# Rigidity, unitary representations of semisimple groups, and fundamental groups of manifolds with rank one transformation group

By Yehuda Shalom

*To Hillel Furstenberg, on the occasion of his 65<sup>th</sup> birthday*

## Contents

## 1. Introduction and discussion of the main results

I. *Introduction.* Throughout the last two or three decades, the theory of rigidity, particularly in relation to semisimple groups and their discrete subgroups, has become an extremely active mathematical field, where tools from diverse areas are simultaneously being employed. In retrospect, it seems that two major achievements in the early 70's have set up somewhat parallel directions in which the theory emerged. These were Mostow's theorem, which, at least in its original form, was concerned with strong rigidity of lattices in rank one Lie groups, and Margulis' theorem, pertaining to the superrigidity of lattices in higher rank groups. Mostow's proof, and the theory which followed his result, are mainly geometric, making a strong use of the hyperbolicity phenomenon in its various forms. The theory of quasi-isometries, negatively curved manifolds (and groups), and the abstract generalization to CAT(-1) spaces, which were deeply pursued by Gromov and his successors, may be viewed as far-reaching developments of Mostow's fundamental ideas. Margulis' theorem, on the other hand, initiated what is now called the (algebraic) ergodic theory of semisimple groups. Zimmer's cocycle superrigidity theorem with its numerous applications, and the work of many others, developed Margulis' ideas to create a striking theory for higher rank groups, and particularly a "nonlinear" finite dimensional rigidity theory. Although drawing on many geometric techniques, this theory relies primarily on measure and ergodic theoretic tools.

A glance at the statement of Mostow-Margulis theorems may hint at the nature of the theories. While an assumption of faithfulness or discreteness on the homomorphism of the lattice is a typical feature of the rank one theory, and owes much to its geometric methods, in the higher rank case one usually aims to get sharp rigidity under extremely weak conditions. Of course, this is hardly a shortcoming of the methods, as it is well-known that the groups $SO(n,1)$ and $SU(n,1)$ fail to have many of the rigidity properties shared by the higher rank groups. It should be mentioned, however, that thanks to the well-known results of Corlette and Gromov-Schoen, the rank one Lie groups $Sp(n,1)$ and $F_{4(-20)}$ do exhibit remarkable rigidity phenomena (see also [Pa]), even if they do not enjoy many of the methods and results the higher rank theory can offer.

In this paper, we attempt to bring closer the rank one and the higher rank theories, in addressing some rigidity questions for the groups $SO(n,1)$, $SU(n,1)$, which are motivated by, and mostly answered in, the higher rank theory. Before elaborating on the precise details, let us describe briefly three examples, which together may help to put our results in a better perspective.

Recall that every isometric action of a higher rank (irreducible) lattice on a tree, fixes a vertex (in most cases this follows from property (T)). For lattices in $SO(n, 1)$ this is well-known *not* be the case. However, we will show that the ambient group does impose a strong limitation on such actions, and even if these lattices need not admit a global fixed vertex, "large" subgroups of them always do. Next, the normal subgroup structure of higher rank lattices, being known to be "trivial" (due to Margulis), is completely analyzed. In particular, the issue of noninjective homomorphisms of them never arises. On the other hand, uniform lattices in rank one groups are word hyperbolic, hence by Gromov's theory, admit many normal subgroups. Actually, some of these lattices are known to surject onto a nonabelian free group, a fact which hardly encourages attempts to study general noninjective homomorphisms. We will nevertheless explore some nonfaithful (discrete) linear representations of these lattices, and show that the algebraic normal subgroup structure is indeed subject to some geometric constraints. Finally, Zimmer applied his cocycle superrigidiy theorem, making use of a deep theorem of Gromov, to show that a fundamental group of a compact manifold on which a higher rank group acts, preserving some geometric structure, must be, roughly speaking, at least as big as the group acting. We will show a similar phenomenon taking place when the group acting is $SO(n, 1)$ or $SU(n, 1)$, relying on a stronger consequence of Gromov's work. In light of the existence of free quotients of lattices, constructions of compact manifolds on which the ambient group acts are quite flexible (using suspension), making such results less expected in general. Our rigidity theorem (see 1.2, 1.3 below) serves here as a bridge, extending both Zimmer's work to the rank one setting, and a result on rank one lattices, due to Yue and Bourdon (see Theorem 1.1 below). While the method we suggest is new, the treatment, using cocycles and the framework of principal bundles, is directly influenced by Zimmer's approach.

The proofs of our results are based heavily, in a new fashion, on unitary representations of semisimple groups, and their cohomology. Although the discussion so far highlights the rigidity aspects, in order to approach these, we prove some results in the representation theory of $SO(n, 1)$ and $SU(n, 1)$, which seem of independent interest. Somewhat surprisingly, it will become evident that the reason for the absence of many rigidity phenomena in these groups, particularly, their failure to have property (T) of Kazhdan, is our main source of strength, being intimately related to the methods we develop. In fact, as far as the subject of asymptotics of matrix coefficients is concerned, we remain with a better understanding of the groups $SO(n, 1)$ and $SU(n, 1)$ than the higher rank groups, and here it is not clear to us whether this has to do only with our approach. (The other rank one simple Lie groups seem to fall, once again, in the middle.) We shall return to elaborate on this matter after describing more fully the main results.

II. *Statement of the main results.* Hereafter, unless specified otherwise, $G$ denotes a simple Lie group with finite center, (locally) isomorphic to $\mathrm{SO}(n,1)$ with $n \geq 3$, or $\mathrm{SU}(n,1)$ with $n \geq 2$ (that is, $\mathrm{SL}_2(\mathbb{R})$ is excluded throughout this section), and $K < G$ denotes a maximal compact subgroup. We fix the $G$-invariant Riemannian metric $d$ on $G/K$, normalized to have constant (maximal) $-1$ sectional curvature when $G = \mathrm{SO}(n,1)$ ($\mathrm{SU}(n,1)$, resp.). This makes the natural equivariant embeddings of the different symmetric spaces, one in the other, isometric. Fix also once and for all the $K$-invariant point $o = \bar{e} \in G/K$ as an "origin" for all the spaces (we may view them all embedded, say, in some $\mathrm{SU}(n,1)/K$ for a large $n$). Recall the notion of *critical exponent* $\delta$, of a discrete subgroup $\Lambda < G$, which will play a central role in the sequel, providing a natural geometric measurement of the size of a subgroup:

1.0 *Definition.*  $\delta(\Lambda) = \inf\{s \mid \sum_{\lambda \in \Lambda} e^{-sd(\lambda o, o)} < \infty\}.$

Note that $\delta(\Lambda)$ does not depend on the ambient group $G$, in the sense that if $\Lambda < G_1 < G_2$, where $G_1 < G_2$ is a natural inclusion, then its value is the same whether we regard $\Lambda$ as a subgroup of $G_1$ or $G_2$. Some well-known properties of $\delta$ to be used, are found in Lemma 2.7 below. We only recall here that if $\Gamma$ is a lattice in $G = \mathrm{SO}(n,1)$ ($n \geq 3$), then $\delta(\Gamma) = n-1$ (a value which we assign, for the time being by formal notation, also to $\delta(G)$), while if $\Gamma$ is a lattice in $G = \mathrm{SU}(n,1)$ ($n \geq 2$), then $\delta(\Gamma) = 2n$ ($= \delta(G)$, as before).

1.1 THEOREM. *Let $\Gamma < G$ be a lattice, and $\varphi : \Gamma \to \Lambda$ be an isomorphism of $\Gamma$ with a discrete subgroup $\Lambda$ of some $H = \mathrm{SO}(m,1)$ or $\mathrm{SU}(m,1)$. Then $\delta(\Gamma) \leq \delta(\Lambda)$.*

Theorem 1.1 extends a result proved in [Yue], where $\Lambda$ is assumed geometrically finite. In [Yue] and [Bou1] the case of equality is studied, which leads to a beautiful superrigidity theorem in this case. After this work was completed, G. Besson, G. Courtois and S. Gallot drew our attention to their recent paper, in which the result of Bourdon and Yue is strengthened, and the assumption on $\Lambda$ removed as well (see [BCG, 1.15]). Yet, in all of the aforementioned papers, the lattice $\Gamma$ is assumed to be uniform, an assumption which is relaxed in Theorem 1.1. Two of our main results present generalizations of this theorem, each in a different direction.

The first direction in which Theorem 1.1 can be extended, concerns uniform lattices $\Gamma$, viewed as (essentially) the fundamental group of the compact $G$-space $G/\Gamma$. Let $G$ be a connected Lie group, and suppose it acts continuously on a compact manifold $M$. Then standard covering space arguments show that the universal covering $\tilde{G}$ of $G$ acts naturally on the universal covering $\tilde{M}$ of $M$. Although we are interested in both families of groups $\mathrm{SO}(n,1)$ and SU$(n, 1)$, the fact that the fundamental group of those in the second family is infinite, makes the formulation of our next result a little more technical for them. For simplicity, we discuss here only $G = \mathrm{SO}(n, 1)$ $(n > 2)$.

1.2 THEOREM. *With the above notation and assumptions on $G$ and $M$, suppose that the $G$-action on $M$ preserves a finite measure, and that the action of $\tilde{G}$ on $\tilde{M}$ is measurably proper with respect to the natural lift of the invariant measure. That is, for almost all $x \in \tilde{M}$ the stabilizer $G_x$ is compact, and the orbital map $g \to gx$ from $\tilde{G}/\tilde{G}_x$ to $\tilde{M}$, is proper. If $\pi_1(M)$ is isomorphic to a discrete infinite subgroup $\Lambda$ of $H = \mathrm{SO}(m, 1)$ or $\mathrm{SU}(m, 1)$ (for some $m$), then $\delta(G) \le \delta(\Lambda)$.*

A simple example of a situation where the condition in Theorem 1.2 holds is when there is at least one element $g \in \tilde{G}$ which acts properly discontinuously on $\tilde{M}$. (This is an easy consequence of the Cartan decomposition, and the rank one property.) Gromov's work on rigid transformation groups, [Gr2], [Zi2, 4.5], provides another general geometric situation in which this condition is satisfied:

1.3 COROLLARY. *Let $G$ be locally isomorphic to $\mathrm{SO}(n, 1)$ for $n \ge 3$.*

(i) *Suppose that the $G$-action on $M$ is real analytic and preserves volume and a rigid (analytic) geometric structure in the sense of Gromov [Gr2] (this includes Cartan's classical notion of a structure of finite type, e.g., a pseudo-Riemannian metric or an affine connection). If $\pi_1(M)$ is infinite, and embeds discretely in $\mathrm{SO}(m, 1)$ or $\mathrm{SU}(m, 1)$, then $\delta(G) \le \delta(\pi_1(M))$.*

(ii) *In particular, in the situation of (i) above, if $\pi_1(M)$ is embedded as an infinite, geometrically finite subgroup, in the same Lie group $G$, then it must be a lattice there.*

Thus, $\mathrm{SO}(n, 1)$ cannot act real analytically on a compact manifold $M$, preserving volume and a connection, if $\pi_1(M)$ is infinite and embeds discretely in $\mathrm{SO}(n-1, 1)$. This type of result holds also when $G = \mathrm{SU}(n, 1)$; however here one has to take into account the possible embedding of the (infinite) center of $\tilde{G}$ in $\pi_1(M)$. We will take up this technicality in Section 7.II, where the proof is presented (see Theorem 7.1 and the preceding discussion).

As mentioned earlier, Theorem 1.2 and Corollary 1.3 were motivated by results of Zimmer in the higher rank case, which, although similar in spirit, are much sharper (see e.g. the account in [Zi2] and the references therein). The superrigidity for cocycles, in conjunction with Ratner's theorem, yields striking consequences, which are simply not valid in our case. We remark that while in Zimmer's work, Gromov's theorem was applied to verify a mild "engaging" condition, we rely on its being stronger, to obtain the fastest possible decay of the matrix coefficients for the $G$-representation on $L^2(\tilde{M})$. It may also be interesting to mention that part (ii) in Corollary 1.3 is now established for *all* simple Lie groups excluding $\mathrm{SL}_2(\mathbb{R})$. The higher rank case, as mentioned, follows from Zimmer's work (see [Zi2, §5]); the case of $\mathrm{Sp}(n,1)$ and $F_{4(-20)}$ from a weaker version of the cocycle superrigidity theorem, established by Zimmer and Corlette via harmonic maps technique [CZ, Theorem 4.2] (and by application of Ratner's theorem in an argument identical to that in [Zi2]). The remaining groups are taken care of above (they are the only ones for which an additional geometric assumption on the image is currently required).

A different way in which Theorem 1.1 may be generalized is by considering also noninjective homomorphisms. There seem to be few known results in this situation.

1.4 THEOREM. *With the assumptions and notation of Theorem* 1.1, *and with $\varphi$ noninjective as well, the following inequality holds*:

$$(1) \qquad \delta(\Gamma) \leq \max\{\delta(\mathrm{Ker}\,\varphi), \frac{\delta(\Gamma)}{2}\} + \max\{\frac{\delta(\varphi(\Gamma))}{2}, 1\}.$$

For instance, if $\delta(\varphi(\Gamma)) < 2$ (e.g., $\varphi(\Gamma)$ is a free group), one necessarily has $\delta(\mathrm{Ker}\,\varphi) \geq \delta(\Gamma) - 1$ (actually, strict inequality must hold; see Theorem 5.6 below). We do not know if inequality (1) is tight in a way more interesting than the fact that equality holds when $\varphi$ extends to $G$. In fact, we believe that equality occurs exactly in that case (see §8.3 below for a further discussion). Indeed, a more symmetric expression on the right-hand side of (1) might be expected, one which resembles better the sum of $\delta$ on the kernel and the image. The appearance of 2 in the denominator has its origins (in the proof), in the special role played by the regular $(=L^2)$ representation. Nevertheless, one should keep in mind the following:

*Observation.* For every $n \geq 2$ there exists a lattice $\Gamma < \mathrm{SO}(n,1)$, and some $\alpha < \delta(\Gamma)$, such that the following holds: for every $\varepsilon > 0$ one can find $\varphi : \Gamma \to \Gamma$ with $\delta(\mathrm{Ker}\,\varphi) < \alpha$ and $\delta(\varphi(\Gamma)) < \varepsilon$.

This is explained after the proof of Theorem 1.4, in Section 7.I. We note that by a result in [BJ], the critical exponent may be replaced (here, as in all of our results) by the Hausdorff dimension of the radial limit set. Notice however, that measuring the full limit set is too crude for our purposes, as the latter is the whole boundary for any infinite normal subgroup of a lattice, whereas the critical exponent of such a subgroup strictly decreases if the quotient by the lattice is nonamenable (see the last paragraph of §7.I).

A result in a more general geometric framework is the following:

1.5 THEOREM. *Let $\Lambda$ be a finitely generated nonamenable group, whose first $\ell^2$-Betti number vanishes. Then for any discrete faithful representation of $\Lambda$ into $\mathrm{SO}(n,1)$ or $\mathrm{SU}(n,1)$, the image has critical exponent at least 2.*

*Moreover, in case of equality, the image must be of divergent type (that is, the series in Definition* 1.0 *diverges at* $s = 2$).

Note that by Sullivan's work [Su1], the additional information in case of equality is of interest. The result itself is tight, as can be seen by lattices in $\mathrm{SL}_2(\mathbb{C})$. There are various families of groups to which Theorem 1.5 applies, such as fundamental groups of closed Kähler manifolds (which are not commensurable to a surface group). See Section 4.III for more examples, applications and references.

1.6 THEOREM. *Let* $\Gamma < G$ *be as in Theorem* 1.1, *and suppose that* $\Gamma$ *acts by isometries on a* (*simplicial*) *tree* $T$, *without a fixed vertex. Then there exists an edge in* $T$, *whose stabilizer* $C < \Gamma$ *satisfies* $\delta(C) \geq \delta(\Gamma) - 1$. *Furthermore, if* $G = \mathrm{SU}(n, 1)$, *or if the tree* $T$ *is locally finite and* $G \neq \mathrm{SO}(3, 1)$, *then for some edge stabilizer* $C$, *strict inequality must hold.*

*In particular, if* $\Gamma = A *_C B$ *is a free product of its subgroups* $A$ *and* $B$, *over the amalgamated subgroup* $C$, *then* $\delta(C) \geq \delta(\Gamma) - 1$, *and strict inequality holds if* $G = \mathrm{SU}(n, 1)$, *or if* $G \neq \mathrm{SO}(3, 1)$ *and* $[A : C] < \infty$, $[B : C] < \infty$.

Notice that we do not assume any kind of geometric condition on the action, or the stabilizers. The bound for $\mathrm{SO}(n, 1)$ is sharp, as is best seen in the constructions of nonarithmetic lattices by Gromov and Piatetski-Shapiro [GPS], which admit a splitting where the amalgamated subgroup $C$ is a lattice in $\mathrm{SO}(n - 1, 1)$. The strengthening to strict inequality relies, in the case of $\mathrm{SU}(n, 1)$, also on a result of Gromov and Schoen [GS] using harmonic maps. In both cases this sharpening depends on an additional detailed spectral analysis, providing support to our belief that equality $\delta(C) = \delta(\Gamma) - 1$ comes only from geometric constructions of splittings of $\Gamma$ (see also §8.2 below). Notice also that some lattices in $\mathrm{SO}(n, 1)$ do admit actions on locally finite trees, e.g., those coming from a homomorphism onto a nonabelian free group.

A different application of our methods concerns the well-known problem of existence of compact quotients for homogeneous spaces. This problem has been addressed using numerous techniques, each providing a solution to some families of groups. Amongst the contributions in the last decade, we mention those of Kobayashi, Labourie, Mozes and Zimmer (see [La] for a survey), and more recently, those of Benoist [Be] and Margulis-Oh ([Mar1], [Oh]). The analysis of Margulis and Oh also uses decay estimates on matrix coefficients, but in a way completely different from the one we suggest. Our approach illustrates another example of a "rank one result", inspired by Zimmer's "higher rank treatment" of the question. Both the formulation and the proof of the theorem are more convenient in the framework of discrete subgroups.

1.7 THEOREM. *Let $G$ be a simple Lie group with finite center, and $\Lambda < G$ be a discrete infinite subgroup which admits a discrete embedding in $\mathrm{SO}(n,1)$ or $\mathrm{SU}(n,1)$ for some $n$. Assume that there exists a nonamenable closed subgroup $L < G$, with noncompact center, which commutes with $\Lambda$. Then $G/\Lambda$ admits no compact quotients; i.e., there is no discrete subgroup $\Gamma < G$ which acts properly discontinuously on $G/\Lambda$, with a bounded fundamental domain.*

For example, a situation, where Theorem 1.7 applies naturally, is the following (compare with Zimmer's questions in [Zi3], and with [LZ]):

COROLLARY. *If $G = \mathrm{SL}_n(\mathbb{R})$, $n \geq 4$, and $H = \mathrm{SL}_2(\mathbb{R}) < G$ is embedded naturally in the upper left corner, then $G/H$ admits no compact quotients. The same is true if $\mathbb{R}$ is replaced by $\mathbb{C}$.*

Replacing $\mathrm{SL}_2(\mathbb{R})$ by a discrete co-compact subgroup $\Lambda$, one indeed reduces the corollary to Theorem 1.7 (observing that $L \cong \mathrm{GL}_2(\mathbb{R})$ can be taken in Theorem 1.7).

III. *The approach.* To describe our approach, we shall first need to discuss further results and notions. Recall that if $G$ is a locally compact group, and $\pi$ is a continuous unitary $G$-representation on the Hilbert space $\mathcal{H}$, then $\pi$ is called *strongly $L^p$*, if for a *dense* subspace $\mathcal{H}_0 \subset \mathcal{H}$, the matrix coefficients $g \to \langle \pi(g)u, v \rangle$ belong to $L^{p+\epsilon}(G)$ for all $\epsilon > 0$, for all $u, v \in \mathcal{H}_0$ (in most examples we shall encounter, the matrix coefficients need not be in $L^p(G)$; see §2.II for further discussion). We note that for notational convenience we allow here any value $0 \leq p \leq \infty$, although standard properties of $L^p$-spaces will be used only for $p \geq 1$ (all unitary representations are strongly $L^\infty$).

Independently, as with all $G$-modules, one has the usual notion of (continuous) first cohomology for a unitary $G$-representation, denoted by $H^1$ (see Definition 3.1 below). Keeping this in mind, we introduce the following invariant, which plays a key role in the present paper:

1.8 *Definition.* For a locally compact group $G$, let $p = p(G)$ denote the value:

$$p(G) = \inf\{0 \leq p \leq \infty \mid \text{there exists some continuous unitary } G\text{-representation}$$
$$\pi \text{ which is strongly } L^p, \text{ and satisfies } H^1(G, \pi) \neq 0\}.$$

If $G$ has no representation $\pi$ with $H^1(G, \pi) \neq 0$ (i.e., $G$ is Kazhdan), set $p(G) = \infty$ (other groups may satisfy $p(G) = \infty$ as well). On the other extreme, we shall see that once $p(G) < 2$, then actually $p(G) = 0$, which is equivalent to having nonvanishing of the first cohomology with coefficients in the regular representation. In the class of finitely generated groups, the latter occurs precisely for groups which are either amenable, or have nonzero first $\ell^2$-Betti number, such as free groups (see Theorem 4.3). For the groups $G = \mathrm{SO}(n,1), \mathrm{SU}(n,1)$ and their lattices, we shall establish the following result:

1.9 THEOREM. *If* $\Gamma < G$ *is a lattice, then* $p(G) = p(\Gamma) = \delta(\Gamma)(= \delta(G))$, *as long as* $G$ *is not locally isomorphic to* $\mathrm{SL}_2(\mathbb{R})$. *In the latter case,* $p(G) = p(\Gamma) = 0$.

See Definition 1.0 for the notation and the numerical values of $\delta$. Notice that this invariant shows, in particular, that lattices in different $\mathrm{SO}(n,1)$'s, say, can be distinguished by their representation theory (one which is far from being understood in any reasonable sense).

The computation of $p(\Gamma)$ when $\Gamma$ is uniform is an easy consequence of the computation of $p(G)$ (to which we shall shortly return), using a restriction-induction argument. The situation is different when $\Gamma$ is nonuniform, as here one does not have in general an isomorphism between $H^*(\Gamma, \pi)$ and $H^*(G, \mathrm{Ind}_\Gamma^G \pi)$. Similar subtle issues, concerning cohomology of nonuniform lattices, were studied mainly by Borel (see [Bor] and the references therein). Here we confine to the first cohomology, but need all unitary (possibly infinite dimensional) $\Gamma$-representations. This question and our approach to it are different, and fortunately we are able to establish the following result, which is essential for the whole treatment of nonuniform lattices in this paper:

1.10 THEOREM. *Assume* $G \cong \mathrm{SO}(n,1)$ *with* $n \geq 4$, *or* $G \cong \mathrm{SU}(n,1)$ *with* $n \geq 2$. *Let* $\Gamma < G$ *be a (nonuniform) lattice, and let* $\pi$ *be some unitary* $\Gamma$-representation. *Then*:

$$(2) \qquad\qquad H^1(\Gamma, \pi) \cong H^1(G, \mathrm{Ind}_\Gamma^G \pi).$$

For the groups excluded in the theorem (which are locally isomorphic to $\mathrm{SL}_2(\mathbb{R})$ or $\mathrm{SL}_2(\mathbb{C})$), we were informed by John Millson that (2) can indeed fail, when $\pi$ is the trivial representation. Our proof of Theorem 1.10 reduces to a question which is independent of $\pi$, concerning the convergence of some integral over $G/\Gamma$. The question is essentially a geometric one, and relies on the structure and precise behavior of the cusps of $G/\Gamma$. The convergence of the integral is crucial in order for a map from the left- to the right-hand side in (2) to be well-defined. We remark that a similar divergence problem invalidates the proof of Raghunathan's local rigidity theorem [Ra] for the nonuniform lattices in $\mathrm{SL}_2(\mathbb{C})$, a result which is indeed false in general. As in [GR, §8], where a substitute for local rigidity was found in this special case, an ad hoc argument is available here, using first $\ell^2$-Betti numbers, to show that $p(\Gamma) = 2$ when $\Gamma < \mathrm{SL}_2(\mathbb{C})$. We also note that in [Sh2], Theorem 1.10 is proved also for irreducible higher rank lattices, using the deep result of [LMR], on the comparison between the word and Riemannian metrics on these groups.

Fix now $G = \mathrm{SO}(n,1)$ or $\mathrm{SU}(n,1)$. The computation of $p(G)$, in Section 4, is carried out in several stages, starting by analysis of the case where $\pi$ is irreducible, and reliance on the classification of those representations with nontrivial cohomology, via Lie algebra cohomology. This proof does not shed light on what seems to be a rather surprising coincidence: the equality in Theorem 1.9. For the groups $\mathrm{SO}(n,1)$ we do indicate a different approach, involving more geometric arguments, proving the result without any classification, and explaining better (though not completely) this phenomenon. At any rate, the above classification yields information only on irreducible representations, and for these, only on the $L^p$ integrability of their $K$-finite matrix coefficients. To implement this classification, we prove the following result, which seems of independent interest:

THEOREM 1.11. *Let $K < G$ be a maximal compact subgroup, and $\pi$ be any unitary $G$-representation. If $\pi$ is strongly $L^p$ for some $2 \leq p < \infty$, then matrix coefficients of all $K$-finite vectors are in $L^{p+\epsilon}(G)$ for all $\epsilon > 0$ (and satisfy an explicit decay estimate, depending on $p$ and not on $\pi$). Furthermore, every irreducible representation in the decomposition of $\pi$ into irreducibles is strongly $L^p$ as well.*

See Theorem 2.1 below for the explicit bound on the $K$-finite matrix coefficients, and for the complete result. Our proof makes essential use of the whole continuous strip of complementary series representations, and tight estimates of the corresponding spherical functions. See also Section 8.4 below for further discussion and several natural questions in this direction which remain unanswered.

We shall see that Theorem 1.11 makes the abstract notion of strongly $L^p$ representations explicit and very amenable to work with. Still, one needs further information to deal with the second ingredient of Definition 1.8. This is provided by the following:

THEOREM 1.12. *Let $G$ be as above, and $\pi$ be any unitary $G$-representation. Let $\Lambda < G$ be any noncompact subgroup. Then the natural restriction map $H^1(G,\pi) \to H^1(\Lambda, \pi|_\Lambda)$ is injective.*

Theorems 1.9–1.12 will be repeatedly used throughout the proofs of the aforementioned rigidity results. Our strategy, roughly, is to show how the situation the group encounters, gives rise to a unitary representation with nonvanishing cohomology (e.g., by use of Theorem 1.12). After analyzing the decay of an appropriately chosen dense subspace of matrix coefficients (not necessarily $K$-finite), we convert the information gained by Theorem 1.9 to the required geometric conclusion.

*Acknowledgments.* This work was initiated while visiting at the University of Chicago, to which we are grateful for its hospitality and support. Especially we thank Bob Zimmer, whose questions during that visit concerning orbit equivalence of rank one lattices (which are not addressed here), led to this paper; his lecture series, given at the NSF/CBMS conference: "Ergodic Theory, Groups and Geometry" in the late stages of this work, motivated and helped us to improve some of the results. We thank the University of Minnesota for its hospitality, and particularly Scot Adams and Dave Witte, for organizing this pleasant and productive workshop. We would also like to acknowledge the contributions of several people to the subject matter in the appendix: Armand Borel, who helped in comparing our result with the known literature; Gregory Margulis, who has drawn our attention to the fact that the criterion we got should in fact be satisfied in dimensions large enough; Shahar Mozes for some helpful discussions; and most of all, Rich Schwartz, who supplied the argument presented there. He deserves full credit for this proof, which replaced the cumbersome tedious arguments we had in mind while writing [Sh1]. We also thank Alain Valette, from whom we learned about the "spaces with walls" in [HP]; Frederic Paulin for communicating that paper; and Bachir Bekka, for indicating an error in a previous version of the paper. Last but not least, we wish to express our deep gratitude and appreciation to the referee, who in a relatively short time examined the paper with great care, pointing out many inaccuracies, and suggesting improvements in various aspects of writing. The final version of the paper owes much to his constructive criticism.

The author is grateful for the NSF support, made through grant DMS-9970974.

## 2. Unitary representations of $SO(n,1)$ and $SU(n,1)$

I. *General preliminaries and notation.* We begin by recalling some of the structure and basic properties of the groups $G = SO(n,1)$, $SU(n,1)$, while introducing some notation.

Let $F$ be the field $\mathbb{R}$ of real numbers, or the field $\mathbb{C}$ of complex numbers. Fix some $n \geq 2$ and let $x_1\bar{x}_1 + \cdots + x_n\bar{x}_n - x_{n+1}\bar{x}_{n+1}$ be a quadratic form of signature $(n,1)$ over $F$, where $x \to \bar{x}$ is the identity in the case $F = \mathbb{R}$, and complex conjugation when $F = \mathbb{C}$. The groups $SO(n,1)$ and $SU(n,1)$ are defined as the subgroups of $SL_{n+1}(F)$ ($F = \mathbb{R}, \mathbb{C}$ resp.) preserving this form. For $F = \mathbb{C}$ this group is connected, but when $F = \mathbb{R}$ the connected component has index 2, and for brevity we keep hereafter the notation $SO(n,1)$ for the connected component.

Fixing the first $m - n$ variables yields a natural embedding of $SO(n,1)$ (resp. $SU(n,1)$) in $SO(m,1)$ (resp. $SU(m,1)$), for $m \geq n$. In addition, $SO(n,1)$ embeds naturally in $SU(n,1)$. These groups are simple Lie groups of real rank one. We fix the embedding of

$$(1) \qquad A \cong \mathrm{SO}(1,1) = \begin{pmatrix} \cosh t & \sinh t \\ \sinh t & \cosh t \end{pmatrix} = a_t$$

in the lower right block, as a Cartan subgroup for all of them. In the case $F = \mathbb{R}$, $K = \mathrm{SO}(n)$ embedded naturally in the upper $n \times n$ submatrix is a maximal compact subgroup. When $F = \mathbb{C}$, $K = S(U(n) \times U(1))$, where $U(n)$ and $U(1)$ ($= \{\lambda \in \mathbb{C} | \ |\lambda| = 1\}$) act only on the first $n$ variables/last variable, resp., is a maximal compact subgroup. Let $M$ be the centralizer of $A$ in $K$. When $F = \mathbb{R}$ we have $M = \mathrm{SO}(n-1)$ (in the upper left $(n-1) \times (n-1)$ submatrix), and when $F = \mathbb{C}$, $M = S(U(n-1) \times \lambda I_2)$ where $\lambda I_2$ denotes the scalar $2 \times 2$ matrices with $|\lambda| = 1$, embedded in the lower right $2 \times 2$ block of $U(n,1)$. Denote by $\tilde{M}$ the normalizer of $A$ in $K$. Then $M \lhd \tilde{M}$ has index 2, and a nontrivial element in the Weyl group $W \cong \tilde{M}/M$ may be represented by the diagonal matrix $w$ whose diagonal has $-1$ in the two entries $n-2$, $n-1$, and 1's elsewhere. We have $wa_t w^{-1} = a_{(-t)}$, and the action of $W$ on $A$ gives rise to two (closures of) Weyl chambers, denoted $A^+ = \{a_t | t \geq 0\}$, $A^- = \{a_t | t \leq 0\}$ (see (1)). The polar (or Cartan) decomposition $G = KA^+K$ holds, where for every $g \in G$ there is a *unique* element $a \in A^+$ with $g = k_1 a k_2$. We may thus define the "Cartan projection" $a : G \to A^+$, by the equality:

$$(2) \qquad g = k_1 a(g) k_2, \qquad k_i \in K, \quad a(g) \in A^+, \quad g \in G.$$

Let $\mathfrak{a} = \mathrm{Lie}\, A$, and let $\exp : \mathfrak{a} \to A$ be the exponential map, with its inverse $\log : A \to \mathfrak{a}$ its inverse. Then $\dim \mathfrak{a} = 1$, and corresponding to our choice of positive Weyl chamber we fix a positive simple root $\beta$. All the positive roots of $\mathfrak{a}$ are either $\beta$ or $2\beta$. For $\mathrm{SO}(n,1)$ the roots are all equal to $\beta$, which has multiplicity $n-1$. Notice that because we identify $A$ as the Cartan subgroup of all the $\mathrm{SO}(n,1)$'s, $\beta$, as functional on $\mathfrak{a}$, does not depend on $n$, and it is the same unique positive simple root for all these groups. The same remark applies to $\mathrm{SU}(n,1)$. Here the root $\beta$ has multiplicity $2n-2$, and $2\beta$ is also a root, with multiplicity 1.

Let $\rho = \rho(G)$ denote half sum of the positive roots of $G$. We define $\delta(G)$ as the (integral) value for which $2\rho = \delta(G)\beta$. That is,

$$(3) \qquad \delta(G) = \left\{ \begin{array}{ll} n-1 & G = \mathrm{SO}(n,1) \\ 2n & G = \mathrm{SU}(n,1) \end{array} \right\}.$$

Since for $\mathrm{SO}(n,1)$, $2\beta$ is not a root, the sum of all the root spaces corresponding to $\beta$ is an abelian subalgebra, which corresponds to an abelian subgroup $N \cong \mathbb{R}^{n-1}$. In the case of $\mathrm{SU}(n,1)$, the appearance of a (one dimensional) root space for $2\beta$, makes $N$ a two-step nilpotent group, the so-called Heisenberg group, of dimension $2n-1$. The following Iwasawa decomposition

holds:

$$(4) \qquad\qquad\qquad G = KAN.$$

The Haar measure of $G$, denoted hereafter simply by $dg$, may be expressed in both polar and in Iwasawa coordinates (2), (4). We shall need both of them in the sequel, but recall now only that of the first one. If $da, dk$ denote the Haar measures of $A, K$ respectively, then the Haar measure of $G$ in terms of the polar decomposition (2) is given by $dg = J(a)dkdadk$, where:

(5)

$$
\begin{aligned}
J(a) &= (e^{\beta(\log a)} - e^{-\beta(\log a)})^{n-1}, & G &= \mathrm{SO}(n,1),\\
J(a) &= (e^{2\beta(\log a)} - e^{-2\beta(\log a)})(e^{\beta(\log a)} - e^{-\beta(\log a)})^{2n-2}, & G &= \mathrm{SU}(n,1),\\
J(a) &\sim e^{\delta(G)\beta(\log a)}, & G &= \mathrm{SO}(n,1), \mathrm{SU}(n,1).
\end{aligned}
$$

II. *Strongly $L^p$ representations*: *The complete formulation of Theorem 1.11*. Before we state the main result of Section 2, Theorem 2.1 below, some remarks on the notion of strongly $L^p$ representation (defined at the beginning of §1.III) are in order. The condition of $L^{p+\varepsilon}$ integrability, rather than simply $L^p$-integrability, may seem peculiar to the reader who comes across this notion for the first time. In fact, in the literature it is the latter condition for strongly $L^p$ representations which is sometimes used. However, it is well-known (see also Theorem 2.1 below), that matrix coefficients for "nice" vectors in irreducible representations of simple Lie groups are, for some $p < \infty$, in $L^{p+\varepsilon}$ for all $\varepsilon > 0$, but not in $L^p$. Consequently, it will turn out to be convenient for us to use this definition. The "nice" vectors alluded to, are those which are $K$-finite. Recall that given a fixed maximal compact subgroup $K$ in a semisimple Lie group $G$, a vector $v \in \mathcal{H}$ is called $K$-*finite* if the dimension of $\mathrm{Span}\{\pi(K)v\}$, denoted henceforth by $\dim_K v$, is finite. Often, the definition of strongly $L^p$ representations for simple Lie groups requires the subspace $\mathcal{H}_0$ to be that of $K$-finite vectors (always dense by the Peter-Weyl theorem), and sometimes this assumption is made implicitly. For our discussion, the distinction between assuming $L^p$-integrability of matrix coefficients for the $K$-finite vectors, and this condition for some dense subspace, is crucial. Our next theorem, the main result of this section, shows in particular that for $\mathrm{SO}(n,1)$, $\mathrm{SU}(n,1)$, these two notions indeed coincide. See Section 8.4 below for more on this issue.

2.1 THEOREM. *Keep the above notation and let $G$ be a simple Lie group with finite center, locally isomorphic to either* $\mathrm{SO}(n,1)$ *or* $\mathrm{SU}(n,1)$ $(n \geq 2)$. *There exists a constant $C = C(G)$, such that for* every *unitary $G$-representation $(\pi, \mathcal{H})$, and $2 \leq p < \infty$, the following conditions are equivalent*:

1. $\pi$ *is strongly $L^p$*.

2. *For every $K$-finite vectors $u, v \in \mathcal{H}$, and every $g \in G$, one has:*

(6)    $|\langle \pi(g) u, v \rangle|$

$$\leq C(1 + \beta(\log a(g))) e^{-\frac{\delta(G)}{p} \beta(\log a(g))} (\dim_K u)^{1/2} (\dim_K v)^{1/2} \cdot \|u\| \cdot \|v\|.$$

*(Recall that $\beta$ is the simple positive root; $a(g)$ and $\delta(G)$ are as in (2) and (3) of §2.I.)*

3. *For every direct integral decomposition $\pi = \int^{\oplus} \pi_x$ (with $\pi_x$'s not necessarily irreducible), the representation $\pi_x$ is strongly $L^p$ for almost every $x$.*

4. *If $\sigma$ is any unitary representation weakly contained in $\pi$, then $\sigma$ is strongly $L^p$ as well.*

   *If, in addition, $\pi$ is irreducible, then 1–4 are also equivalent to:*

5. *There exists one nonzero diagonal matrix coefficient of $\pi$ which is in $L^{p+\varepsilon}(G)$ for all $\varepsilon > 0$.*

Much of the rest of Section 2 will be devoted to proving Theorem 2.1. We remark here that a similar proof works for every rank one simple algebraic group over a non-Archimedean local field, and in fact for "large" subgroups of the automorphism group of a regular tree (e.g., unimodular noncompact subgroups acting transitively on its boundary).

III. *Preliminaries toward the proof of Theorem* 2.1. We begin with a discussion which essentially amounts to showing that (1) implies (2) in Theorem 2.1 for $p = 2$, a result which follows also from [CHH]. However, we present an alternative, simple approach, and establish along the way other results which will be needed in the sequel. The reader interested in a short path to the basics of the pointwise decay estimates (which are fundamental in what follows), may find our treatment helpful (although Theorem 2.2 seems new, and may be of interest also to specialists).

2.2 THEOREM. *Let $G$ a unimodular, second countable, locally compact group, and $K < G$ a compact subgroup. Denote by $\lambda = L^2(G)$ the regular $G$-representation, and let $(\pi, \mathcal{H})$ be some unitary $G$-representation satisfying the following two properties*:

(i) *There exists a nonzero vector $v \in \mathcal{H}$ such that $\langle \pi(g) v, v \rangle \geq 0$ for all $g \in G$ (we call such $v$ a $G$-positive vector).*

(ii) *$\pi^K$, the subspace of $K$-invariant vectors, has dimension 1.*

   *Let $v_K \in \mathcal{H}$ be the unique (up to scalar) nonzero $K$-invariant unit vector in $\mathcal{H}$. Then for every $K$-invariant functions $\varphi, \psi \in L^2(G)$ with $\|\varphi\| = \|\psi\| = 1$, and for every $g \in G$:*

(7)                 $|\langle \lambda(g) \varphi, \psi \rangle| \leq \langle \pi(g) v_K, v_K \rangle.$

*Remark.* If assumption (ii) is removed, we still have an inequality similar to (7), except that on the right-hand side a supremum over the $K$-invariant vectors has to be taken. If (ii) is replaced by the more flexible assumption, $\dim \pi^K < \infty$, this supremum becomes again a maximum (this more general version may be useful in various examples of totally disconnected groups).

*Proof.* For a finite positive measure $\nu$ on $G$, and a unitary $G$-representation $\sigma$, denote by $\sigma(\nu)$ the corresponding convolution operator, defined by $\langle \sigma(\nu)u, v \rangle = \int \langle \sigma(g)u, v \rangle \, d\nu(g)$. Let $m$ be the normalized Haar measure of $K$. If $v$ is a positive $G$-vector, then $\pi(m)v$ is clearly nonzero, $G$-positive, and $K$-invariant; hence by (ii) we deduce that $v_K$ is $G$-positive. Therefore the right-hand side in (7) is nonnegative. By replacing the functions with their absolute values, we may assume that $\varphi$ and $\psi$ are nonnegative, discarding the absolute value on the left-hand side of (7). Suppose now that (7) fails for some $g_0 \in G$. Notice that by $K$-invariance, the set of $g$'s for which it fails is nonempty, open, and bi-$K$-invariant. We may therefore find a probability measure $\mu$ on $G$ (continuous with respect to the Haar measure), which is bi-$K$-invariant, and such that (7) fails in its support. This implies that $\|\lambda(\mu)\| \geq \langle \lambda(\mu)\varphi, \psi \rangle > \langle \pi(\mu)v_K, v_K \rangle$. We now observe that $\langle \pi(\mu)v_K, v_K \rangle = \|\pi(\mu)\|$, which together with the previous inequality yields $\|\lambda(\mu)\| > \|\pi(\mu)\|$. To prove this observation notice that by left $K$-invariance of $\mu$, for every $u$ the vector $\pi(\mu)u$ is $K$-invariant; hence $\pi(\mu)u = t(u)v_K$, for some $t(u) \in \mathbb{C}$. Since $\pi(m)u$ is $K$-invariant, we have $\pi(m)u = s(u)v_K$ where $|s(u)| \leq \|u\|$. Using the right $K$-invariance of $\mu$ we get for every unit vectors $u, w \in \mathcal{H}$:

$$|\langle \pi(\mu)u, w \rangle| = |\langle \pi(\mu)\pi(m)u, w \rangle| = |\langle \pi(\mu)s(u)v_K, w \rangle|$$
$$\leq |\langle \pi(\mu)v_K, w \rangle| = |t(v_K)\langle v_K, w \rangle|$$
$$\leq |t(v_K)\langle v_K, v_K \rangle| = \langle \pi(\mu)v_K, v_K \rangle.$$

Hence $\|\pi(\mu)\| \leq \langle \pi(\mu)v_K, v_K \rangle \leq \|\pi(\mu)\|$, and equality must hold, as claimed. Now, the inequality $\|\lambda(\mu)\| > \|\pi(\mu)\|$ we have established, contradicts the following result, whose proof therefore completes the proof of Theorem 2.2.

**2.3 Lemma.** *Let $\lambda = L^2(G)$ be the regular $G$-representation, and $\pi$ be a unitary $G$-representation with a positive $G$-vector. Then for every positive, continuous, compactly supported, probability measure $\mu$ on $G$, one has*: $\|\pi(\mu)\| \geq \|\lambda(\mu)\|$.

*Proof.* We first reduce the lemma to the case where $\mu$ is symmetric. Indeed, denoting $\hat{\mu}(B) = \mu(B^{-1})$, we see that the convolution $\nu = \hat{\mu} * \mu$ is symmetric, and for every unitary $G$-representation $\sigma$, one has:

$$\|\sigma(\nu)\| = \|\sigma(\hat{\mu})\sigma(\mu)\| = \|\sigma(\mu)^*\sigma(\mu)\| = \|\sigma(\mu)\|^2.$$

Hence, taking $\sigma = \pi, \lambda$ shows that proving the lemma for the symmetric measure $\nu$ suffices. Now, denoting by $\mu^m$ the $m$-th convolution, we have $\pi(\mu^m) = \pi(\mu)^m$. Taking a bounded neighborhood $U$ of the identity such that $\langle \pi(g)v, v \rangle \geq 1/2$ for all $g \in U$, we get for the $G$-positive vector $v$:

$$\langle \pi(\mu)^{2n}v, v \rangle^{1/2n} = \langle \pi(\mu^{2n})v, v \rangle^{1/2n}$$
$$= \left( \int_G \langle \pi(g)v, v \rangle d\mu^{2n}(g) \right)^{1/2n} \geq \left( \frac{1}{2}\mu^{2n}(U) \right)^{1/2n}$$

(for the last inequality, the positivity of $v$ was used). However, it is well-known and easily verified that for every locally compact group $G$, a symmetric continuous (compactly supported) probability measure $\mu$, and a bounded open set $U \subseteq G$, one has $\mu^{2n}(U)^{1/2n} \to \|\lambda(\mu)\|$ (for discrete groups this result is due to Kesten [Ke] and similar arguments apply in general). It follows that $\|\pi(\mu)\| \geq \liminf \langle \pi(\mu)^{2n}v, v \rangle^{\frac{1}{2n}} \geq \|\lambda(\mu)\|$, which completes the proofs of Lemma 2.3 and Theorem 2.2. $\qquad\square$

*Remark.* Notice that Lemma 2.3 yields a lower bound in terms of the regular representation on all the norms of convolution operators for representations with positive vectors. One family of such representations is that of the form $\pi = L^2(X, m)$, where $(X, m)$ is a quasi-invariant $\sigma$-finite $G$-measure space. Another is the family of representations of the form $\pi = \sigma \otimes \bar{\sigma}$, where $\sigma$ is *any* unitary $G$-representation and $\bar{\sigma}$ is its contragredient dual (indeed, any vector $v \otimes \bar{v}$ is then positive). Therefore, $\|\pi \otimes \bar{\pi}(\mu)\| \geq \|\lambda(\mu)\|$ for any unitary $G$-representation $\pi$, where $\lambda$ and $\mu$ are as in Lemma 2.3. This generalizes a result of Pisier [Pis], and when $G$ is amenable, shows that $\pi \otimes \bar{\pi}$ always contains weakly the trivial representation (which is the main result of [Bek1]). After the completion of the paper, we learned from B. Bekka that this notion of "positive vectors" had been studied by R. Godement [Go], who asked, in modern terminology, whether any representation (of an amenable group) with a positive vector, contains weakly the trivial representation. An affirmative answer for some groups was obtained by Bekka in [Bek2]. Lemma 2.3 yields the result for all (amenable) groups.

As we observe below, (7) can be applied to get a sharp decay estimate on $K$-finite matrix coefficients, once an appropriate representation $\pi$ is present. We will be interested, however, not only in decay estimates for the regular representation, but in a more general class. Indeed, it is clear from the proof of Theorem 2.2 that the same result holds if we replace $\lambda$ by any unitary representation $\sigma$ which is weakly contained (see 2.4 below) in $\lambda$. To see that, observe that for any such representation we have $\|\sigma(\mu)\| \leq \|\lambda(\mu)\|$ for every (not necessarily positive) continuous measure $\mu$. All that was proved in the second part of the proof is the inequality $\|\lambda(\mu)\| \le \|\pi(\mu)\|$, which then continues to hold if $\lambda$ is replaced by $\sigma$. As was shown by Eymard [Ey], the condition on norms of convolution operators is in fact equivalent to the weak containment of $\sigma$ in $\lambda$. A more commonly known definition is the following:

2.4 *Definition.* The $G$-representation $\sigma$ is said to be *weakly contained* in $\lambda$, denoted $\sigma \prec \lambda$, if every $n \times n$ submatrix of $\sigma$ can be approximated, uniformly on compact subsets in $G$, by $n \times n$ submatrices of some multiple (perhaps countable) direct sum of $\lambda$.

See [CHH], [Sh3] and the references therein, for more on this notion. We now wish to derive from (7) an estimate for the matrix coefficients of all $K$-finite vectors, not only $K$-fixed ones. For later use, we give this argument in a general framework (compare with [CHH] and [Ho, §7]).

2.5 Proposition. *Keep the notation of Theorem 2.2, but denote now by $\lambda$ the $G$-representation on $L^2(M)$, where $M$ is some $\sigma$-finite, Borel measure $G$-space, on which $G$ acts measurably by measure-preserving transformations (e.g., $M = G$). Suppose that the function $F : G \to \mathbb{R}$ has the property that every $K$-invariant unit vectors $\varphi, \psi \in L^2(M)$ satisfy for all $g \in G$ the inequality $|\langle \lambda(g)\varphi, \psi \rangle| \le F(g)$. Let $\sigma$ be a $G$-representation with $\sigma \prec \lambda$, and let $u, v$ be some $K$-finite vectors in $\mathcal{H}_\sigma$. Then for all $g \in G$:*

$$(8) \qquad |\langle \sigma(g)u, v \rangle| \le F(g) \cdot (\dim_K u)^{1/2} (\dim_K v)^{1/2} \cdot \|u\| \cdot \|v\|.$$

*Proof.* We claim that it is enough to establish (8) when $\sigma = \lambda$. Indeed, first notice that this implies (8) when $\sigma = \infty \cdot \lambda$. Next, fix $K$-finite vectors $u, v \in \mathcal{H}_\sigma$, and apply the fact that every $n \times n$ submatrix of $\sigma$ can be approximated uniformly on compact subsets by submatrices of $\infty \cdot \lambda$ to the submatrix associated with a finite set of vectors including basis elements for both Span$\{Kv\}$ and Span$\{Ku\}$. Restricting our attention to the $K$-action, it follows by its compactness, that the $K$-representations on Span$\{Kv\}$ and Span$\{Ku\}$ are (properly) contained in $\infty \cdot \lambda$ (restricted to $K$), and that moreover the "approximating" vectors in $\infty \cdot \lambda$ can be chosen to span isomorphic $K$-representations (this is a standard exercise related to the discreteness of the so-called Fell topology, on unitary duals of compact groups). In particular, the dimension of the $K$-spans are the same, and since we are assuming now that the pointwise estimate (8) holds for the pointwise approximating matrix coefficients, it must also be satisfied by $\langle \sigma(g)u, v \rangle$.

We are thus reduced to the case $\sigma = \lambda$. To establish (8) in that case, we first fix a countable dense subgroup $K_0 < K$ (recall $G$ is assumed second countable). Given any $K$-finite vectors $u, v \in L^2(M)$, define the functions:

$$\varphi(m) = \sup\{|u(km)| \mid k \in K_0\}, \qquad \psi(m) = \sup\{|v(km)| \mid k \in K_0\}.$$

Observe that since $K_0$ is countable, the functions $\varphi$ and $\psi$ are well-defined as measurable functions, in the sense that they do not depend on the choice of $u$ and $v$ in their class of almost everywhere equal functions. The values of $\varphi(m)$ or $\psi(m)$ may be infinite, but we shall shortly see that $\varphi$ and $\psi$ are in $L^2(M)$, and in particular this can happen only on a null set. Once the latter $L^2$-estimate is established, we deduce from the fact that $\varphi$ and $\psi$ are $K_0$-invariant, and the continuity of the $K$-action on $L^2(M)$ ([Zi1, Appendix B.3]), that they are actually invariant under any $k \in K$ (as measurable functions). Thus, as $\varphi$ and $\psi$ are also positive, and their matrix coefficient dominates (the absolute value of) that of $u$ and $v$ resp., in order to deduce (8) it only remains now to show that their $L^2$-norms satisfy $\|\varphi\| \leq (\dim_K u)^{1/2} \cdot \|u\|$, and $\|\psi\| \leq (\dim_K v)^{1/2} \cdot \|v\|$.

Assume for simplicity $\|u\| = 1$ and set $n = \dim_K u$. Let $u_1, \ldots, u_n$ be an orthonormal basis for $Ku$. Then, for every $k \in K$ there exist unique complex numbers $\theta_1(k), \ldots, \theta_n(k)$, such that $u(km) = \sum \theta_i(k)u_i(m)$ for almost every $m \in M$. In fact, by the countability of $K_0$, almost every $m$ satisfies the latter equality simultaneously for all $k \in K_0$. Obviously, $\sum |\theta_i(k)|^2 \equiv 1$. Now, by the Cauchy-Schwartz inequality, almost every $m \in M$ satisfies for all $k \in K_0$:

$$|u(km)|^2 = |\sum_{i=1}^{n} \theta_i(k)u_i(m)|^2 \leq \sum_{i=1}^{n} |\theta_i(k)|^2 \cdot \sum_{i=1}^{n} |u_i(m)|^2 = \sum_{i=1}^{n} |u_i(m)|^2.$$

Therefore $\varphi^2(m) \leq \sum |u_i(m)|^2$ almost everywhere, and

$$\|\varphi\| = \int_M \varphi^2(m)^{1/2} \leq (\int_M \sum_{i=1}^{n} |u_i(m)|^2)^{1/2} = (\sum_{i=1}^{n} \int_M |u_i(m)|^2)^{1/2} = n^{1/2}$$

as required. A similar analysis applies to the function $\psi$. $\qquad\square$

Theorem 2.2 and Proposition 2.5 yield together a pointwise estimate for all $K$-finite vectors, in all the representations which are weakly contained in the regular representation, in terms of a single matrix coefficient of a special representation. In general, there may be several (even irreducible) representations satisfying the assumptions of Theorem 2.2. Notice that the best estimate will be achieved if a representation $\pi$ as in Theorem 2.2 exists, which is itself weakly contained in the regular representation, for then the pointwise bound that it puts is sharp, being attained by the very same distinguished matrix coefficient. It is perhaps surprising that for any connected topological group (with $K$ a maximal compact subgroup), such a representation can indeed be found (see the discussion after the following result). Although we shall not need it in the sequel, for completeness let us observe its uniqueness.

2.6 PROPOSITION. *Let $G$ be a locally compact, second countable group, and $K$ a compact subgroup. Suppose that there exists a unitary $G$-representa-* *tion which satisfies assumptions* (i) *and* (ii) *of Theorem 2.2, and is also weakly contained in the regular representation. Then there exists such an* irreducible *representation, and it embeds in any other representation with these properties. In particular, there is at most one irreducible representation satisfying all three properties.*

*Proof.* Suppose that $\pi$ satisfies all the properties, and that $\pi_1$ is another such representation. Let $v, v_1$ be $K$-invariant unit vectors in the representation space of $\pi$ and $\pi_1$, respectively. Then by Theorem 2.2 and Proposition 2.5 we have $\langle \pi(g)v, v \rangle \le \langle \pi_1(g)v_1, v_1 \rangle$ for all $g \in G$, and similarly, the opposite inequality holds. Thus, we have equality, which implies that the closures of the $G$-span of $v$ and $v_1$ are isomorphic. Denoting the $G$-representation on these isomorphic subrepresentations by $\sigma$, we see that $\sigma$ satisfies all the properties, and our argument shows that it embeds in any other representation sharing them. We show $\sigma$ is irreducible. Let $\mathcal{H}_\sigma$ be the representation space for $\sigma$, and $v_\sigma$ be the unique (up to scalar) $K$-invariant vector. Suppose we have a $G$-invariant decomposition $\mathcal{H}_\sigma = \mathcal{H}_1 + \mathcal{H}_2$, and write correspondingly $v_\sigma = v_1 + v_2$. As $Gv_\sigma$ spans $\mathcal{H}_\sigma$, if $v_i = 0$ then $\mathcal{H}_i = 0$ ($i = 1, 2$). However, both $v_1, v_2$ are $K$-invariant, hence multiples of $v_\sigma$, and orthogonal to each other. This is impossible unless one of them is zero. $\qquad\square$

Notice that if $G$ is totally disconnected, then a compact open subgroup $K < G$ always exists, and if $P < G$ is any cocompact subgroup, then the $G$-representation on $L^2(G/P)$ has a finite dimensional subspace of $K$-invariant vectors. Thus, (the more general version of) Theorem 2.2 can be applied, and if such amenable $P$ exists, the additional condition in 2.6 is satisfied as well. On the other extreme, this phenomenon occurs when $G$ is connected, and a representation $L^2(G/P)$ of the latter form, satisfying the three conditions in 2.6, always exists (when $K$ is maximal compact). The general case here is reduced by a standard argument to that of a semisimple group, and although the latter may be dealt with generally, we shall discuss henceforth only the simple Lie groups which are of interest to us.

We retain the notations of subsection I. Let $G = KAN$ be as in (4), and $M$ be the centralizer of $A$ in $K$. Then $P = MAN$ is an amenable subgroup, being a compact extension of the solvable group $AN$. We have $G = K \cdot P$. The quasi-regular $G$-representation on $L^2(G/P)$ obviously satisfies the two conditions of Theorem 2.2 (as $K$ is transitive on $G/P$). It is also weakly contained in $L^2(G)$, since $P$ is amenable; thus all three conditions in 2.6 hold. When we choose a quasi-invariant measure on $G/P$ which is invariant under $K$, the Harish-Chandra $\Xi$-function is defined as the diagonal matrix coefficient corresponding to the constant unit function. By 2.2 and 2.5, taking $F(g) = \Xi(g)$ in (8) yields a (tight) bound on matrix coefficients of all $K$-finite vectors in representations $\sigma \prec L^2(G)$. The first estimates of the $\Xi$-function are due to Harish-Chandra. Here we shall use the following (sharp) estimate for the rank one Lie groups (see [GV, 4.6.4]):

$$(9) \qquad \Xi(g) \le C(1 + \beta(\log a(g))) \cdot e^{-\frac{1}{2}\delta(G)\beta(\log a(g))}.$$

Preparing still for the proof of Theorem 2.1, we recall now some basic properties of the class one (also called "spherical") representations of $G = \mathrm{SO}(n,1)$, $\mathrm{SU}(n,1)$.

Recall that a class one (irreducible) $G$-representation, is one that contains a nonzero vector invariant under $K$. Such vector is unique up to scalar. A general method to construct such representations is as follows: Let $\mathfrak{a}_{\mathbb{C}}$ be the complexification of $\mathfrak{a} = \mathrm{Lie}\ A$, and $\lambda \in \mathfrak{a}_{\mathbb{C}}^*$. Then $man \mapsto e^{\lambda(\log a)}$ defines a character on $P = MAN$, which is unitary if and only if $\lambda \in i\mathfrak{a}^*$. Inducing (unitarily) this $P$-representation to $G$, yields a representation $\pi_\lambda$ which is obviously unitary if $\lambda \in i\mathfrak{a}^*$. In that case $\pi_\lambda$ is called a *principal series (spherical) representation*. Denote by $v_\lambda$ a $K$-invariant unit vector, with $\varphi_\lambda(g) = \langle \pi_\lambda(g)v, v \rangle$. Since inducing any unitary representation from an amenable subgroup yields a representation which is weakly contained in the regular representation, from Theorem 2.2 and Proposition 2.5 it follows that $|\varphi_\lambda(g)| \le \Xi(g)$. For $\lambda = 0$ we have, by definition, equality. There are more class one representations $\pi_\lambda$, which are obtained by taking $\lambda \in \mathfrak{a}^*$ with $0 \le \lambda \le \frac{\delta(G)\beta}{2} = \rho$. (Here one has to change the inner product.) These are called the *complementary series (spherical) representations*. For $\lambda = \rho$, $\pi_\lambda$ is the trivial representation. The corresponding spherical functions $\varphi_\lambda(g)$ satisfy the following estimate (see [GV, 5.1 (5.1.18)]), for some constant $C$ depending only on $G$ :

$$(10) \qquad |\varphi_\lambda(g)| \le C(1 + \beta(\log a(g)))e^{(\lambda - \rho)\log a(g)}.$$

As the $\pi_\lambda$'s are irreducible, having one matrix coefficient in $L^p(G)$ implies the same for a dense subset (linear combinations of $G$-translations); hence from (10) and (5) in Section 1 we get:

$(*) \qquad$ For $\lambda = |\lambda|\beta \quad (0 \le |\lambda| \le |\rho| = \frac{1}{2}\delta(G)), \quad \pi_\lambda$ is strongly $L^p$ if
$\qquad\qquad p(\frac{1}{2}\delta(G) - |\lambda|) \ge \delta(G), \quad \text{i.e.,} \quad p \ge 2\delta(G)/\delta(G) - 2|\lambda|.$

Lower bounds on the spherical functions of complementary series representations will be needed as well. It turns out that these spherical functions are positive, and for some constant $C' = C'(G) > 0$ (see [GV, 4.7 (4.7.13)]) satisfy:

$$(11) \qquad C'e^{(\lambda - \rho)\log a(g)} \le \varphi_\lambda(g) \qquad \rho/2 \le \lambda \le \rho.$$

By the reference cited above, (11) holds for every $0 \le \lambda$, with $C'$ depending on $\lambda$ (when $\lambda = 0$, $C' = 1$, see [GV, 4.5]). Actually (11) holds for all $0 \le \lambda \le \rho$, but we have not found an easy reference for that. For $\lambda > 0$, the ratio between the two sides of (11) tends to a positive bounded limit when $g \to \infty$, so (11) reflects the true asymptotics of $\varphi_\lambda$. This fact will be of considerable importance to us later (see (1) in step (ii) of Theorem 4.2 below, and the remark thereafter).

The above principal and complementary series representations exhaust all the class one irreducible representations of $\mathrm{SO}(n,1)$, $\mathrm{SU}(n,1)$. The existence of the whole continuous strip of complementary series representations is special to these two families of groups, and will be essential for the proof of Theorem 2.1, to which we can now turn.

IV. *Proof of Theorem* 2.1.

$1 \implies 2$. Consider first the case $p = 2$. By a general result in [CHH] (which applies to any $G$, not necessarily a simple Lie group), it follows that $\pi$ is weakly contained in the regular representation. The bound (6) then follows from Theorem 2.2, Proposition 2.5 (with $M = G$), and the bound (9) for the Harish-Chandra $\Xi$-function (which is taken as $F(g)$ in Proposition 2.5). The constant $C$ obtained here is the same as in (9). In general we will show that this $C$ can be multiplied by at most $1/C'^2$, where $C'$ is taken from (11).

We may assume that $2 < p < \infty$. Let $q$ be the value for which $\frac{1}{q} + \frac{1}{p} = \frac{1}{2}$. Let us assume first that $p \leq 4$, which means $q \geq 4$. Take $\lambda$ for which the spherical function $\varphi_\lambda$ is in $L^{q+\varepsilon}(G)$ for all $\varepsilon > 0$, but not in $L^q(G)$. This, in light of estimates (10), (11), singles out the $\lambda$ satisfying $q(\rho - \lambda) = 2\rho$. Moreover, $\pi_\lambda$ is strongly $L^q$ (see $(*)$ above). As $q \geq 4$, we have $\lambda \geq \rho/2$. Consider now the tensor product $\pi \otimes \pi_\lambda$. For every $v_1, v_2$ in the assumed dense subspace of $\pi$, and $u_1, u_2$ in the dense subspace of vectors with $L^{q+\varepsilon}$ matrix coefficients in $\pi_\lambda$, we have by the choice of $q$ and the Hölder inequality:

$$\langle (\pi \otimes \pi_\lambda)(g)v_1 \otimes u_1, v_2 \otimes u_2 \rangle = \langle \pi(g)v_1, v_2 \rangle \cdot \langle \pi_\lambda(g)u_1, u_2 \rangle \in L^{2+\varepsilon} \quad \text{for all } \varepsilon > 0.$$

Since this holds also for sums of tensors $v_i \otimes u_i$, and these form a dense subspace of the representation space for $\pi \otimes \pi_\lambda$, we deduce that $\pi \otimes \pi_\lambda$ is strongly $L^2$. The case $p = 2$ at the beginning of the proof, yields now the inequality (6) (with $p = 2$) for all $K$-finite vectors in $\pi \otimes \pi_\lambda$.

Let now $u, w$ be any $K$-finite vectors for the representation $\pi$, and let $v_\lambda$ be the unit $K$-invariant vector in $\pi_\lambda$. Then $u \otimes v_\lambda$, $w \otimes v_\lambda$ are $K$-finite vectors for $\pi \otimes \pi_\lambda$, with $\dim_K(u \otimes v_\lambda) = \dim_K u$ and $\dim_K(w \otimes v_\lambda) = \dim_K w$. By the above discussion we get, assuming $u$ and $w$ are also normalized:

$$|\langle \pi(g)u, w \rangle| \cdot |\langle \pi_\lambda(g)v_\lambda, v_\lambda \rangle| = |\langle (\pi \otimes \pi_\lambda)(g)u \otimes v_\lambda, w \otimes v_\lambda \rangle|$$
$$\leq C \cdot (1 + \beta(\log a(g)))e^{-\rho(\log a(g))} \cdot (\dim_K u)^{1/2} \cdot (\dim_K w)^{1/2}.$$

From this and (11), we now deduce:

$$|\langle \pi(g)u, w \rangle| \leq \frac{C}{C'}(1 + \beta(\log a(g)) \cdot e^{-\lambda(\log a(g))} \cdot (\dim_K u)^{1/2} \cdot (\dim_K w)^{1/2}.$$

To see that this is indeed the required estimate, it only remains to verify that $\frac{\delta(G)\beta}{p}(=\frac{2\rho}{p})=\lambda$. Since $\lambda$ and $\rho$ are (nonzero) multiples one of the other, the ratio $\frac{\lambda}{2\rho}$ makes sense, and we have to show that $\frac{1}{p}=\frac{\lambda}{2\rho}$. As $\lambda$ was chosen such that $\frac{1}{q}=\frac{\rho-\lambda}{2\rho}$, the desired equality $\frac{1}{p}=\frac{\lambda}{2\rho}$ is equivalent to $\frac{1}{p}+\frac{1}{q}=\frac{\lambda}{2\rho}+\frac{\rho-\lambda}{2\rho}=\frac{1}{2}$, which is exactly how $q$ was chosen. This completes the case $p \leq 4$.

If $p > 4$ take again $q$ with $\frac{1}{p}+\frac{1}{q}=\frac{1}{2}$. Now $2 < q < 4$, so $t = 2q > 4$. Let $\lambda'$ be such that the spherical function $\varphi_{\lambda'}$ is "exactly" in $L^{t+\varepsilon}(G)$ (again, see ($*$) above). Then $\rho/2 \leq \lambda'$, and we may repeat the whole previous argument, replacing $\pi_\lambda$ with $\pi_{\lambda'} \otimes \pi_{\lambda'}$ (which is strongly $L^q$ by Hölder), and the vector $v_\lambda$ by $v_{\lambda'} \otimes v_{\lambda'}$. The argument proceeds similarly, with the one observation that the constant may increase up to $C/C'^2$. This completes the proof of $1 \implies 2$.

$2 \implies 4$. An argument identical to that in the first part of the proof of Proposition 2.5, shows that the $K$-finite matrix coefficients of $\sigma$ satisfy the same estimate (6). By (5) it is easy to see that they are in $L^{p+\varepsilon}$ for all $\varepsilon > 0$.

$4 \implies 3$. If $\pi = \int \pi_x$ then, as is well-known, for almost every $x$ the representation $\pi_x$ is weakly contained in $\pi$; so this implication is obvious.

$3 \implies 1$. Taking a trivial decomposition $\pi = \pi$, we see that this claim is again obvious.

Finally, if some nonzero $v \in \mathcal{H}$ satisfies $\langle \pi(g)v, v \rangle \in L^{p+\varepsilon}(G)$, then by invariance of Haar measure it is easy to see that the matrix coefficients of any two vectors in the (algebraic) span of $\pi(G)v$ have this property, a subspace which is dense by irreducibility. This completes the proof of Theorem 2.1. $\square$

V. *Unitary representations and discrete subgroups of* $\mathrm{SO}(n,1), \mathrm{SU}(n,1)$. The notion of critical exponent (see Definition 1.0 above) plays a central role in our discussion. This and other related notions such as the limit set and its Hausdorff dimension are covered, e.g., in [Ni]. Useful information concerning the relation to spectral analysis appears in [Cor], and some extensions to general CAT(-1) spaces can be found in [BM] (and the references therein). For the benefit of the reader, we collect here some basic properties of the critical exponent, which are used in the sequel.

2.7 LEMMA. *Fix* $n \geq 2$. *Retain the notation in Definition* 1.0 *and the discussion preceding it.*

(1) *If* $\Gamma$ *is a lattice in* $G = \mathrm{SO}(n,1)$ *then* $\delta(\Gamma) = n - 1$. *If it is a lattice in* $G = \mathrm{SU}(n,1)$ *then* $\delta(\Gamma) = 2n$.

(2) *If* $\Lambda < G$ *then* $\delta(\Lambda) \leq \delta(G)$. *If equality holds and* $\Lambda$ *is geometrically finite, then* $\Lambda$ *is a lattice.*

(3) *If $G$ is as in* (1) *and* $\Lambda < G$ *is a discrete nonamenable subgroup, then for every* $\epsilon > 0$ *there is a free subgroup* (*say, on two generators*) $\Gamma < \Lambda$, *satisfying* $\delta(\Gamma) < \epsilon$.

*Proof.* We remark first on the initial statement of (2). In [BM, §1] the notion of critical exponent is extended naturally to any (not necessarily discrete) closed subgroup of the isometry group of a locally compact CAT(-1) space. The discussion there, Lemma 2.8 below, and the expression (5) for the Haar measure, show that our formal notation for $\delta(G)$ in (3) above, indeed agrees with the general geometric definition of the critical exponent. In particular, the first claim of (2) follows from [BM, Lemma 1.5]. Next, for a geometrically finite subgroup $\Lambda < G$, the bottom of the Laplacian spectrum on the locally symmetric space $\Lambda \backslash G / K$ is achieved in the discrete $L^2$-spectrum (see [Su2, 2.21], [Cor, §4]), and by its relation with the critical exponent (see the proof of Theorem 2.10 below), we have that $\delta(\Lambda) = \delta(G)$ if and only if $0$ is an ($L^2$-)eigenvalue of the Laplacian, which happens if and only if $\Lambda$ is a lattice. This accounts for (1) and the other direction in (2). For part (3) we pick, using the nonamenability of $\Lambda$, two hyperbolic elements of it, say $a, b$, whose attracting and repelling points on the boundary $G/P$ are all different (their existence follows easily from the nonamenability of $\Lambda$). As is well-known (by Tits' "ping-pong argument"), for $n$ large enough $a^n$ and $b^n$ generate a free group. Using the fact that the metric contraction of $a^n$ and $b^n$ on "most" of the boundary grows linearly with $n$, an elementary argument shows that the Hausdorff dimension of the limit set of the group generated by $a^n$ and $b^n$ is dominated by some constant multiple of $n^{-1}$. Since that dimension always bounds the critical exponent from above (cf. [BJ]), the result follows.  □

2.8 LEMMA. *Let* $a : G \to A^+$ *be as defined in* (2), *and let* $\beta$ *be the positive simple root as before. Then for every* $g \in G$,

$$(12) \qquad\qquad d(go, o) = \beta(\log a(g)).$$

*Proof.* The equality (12) is well-known and easy to verify. Briefly, notice that both sides of (12) are bi-$K$-invariant functions on $G$, and thus depend for their value on $A^+$. On $A^+$ both sides are additive characters, so that there is only an issue of normalization. This may be checked by one example, say, in the two dimensional real hyperbolic space. Another way the reader may be convinced, is to recall that the inverse of the exponential of the two expressions in (12) is in $L^p(G)$ for the same (open) set of $p$'s, the infimum of which is $\delta(G)$ (see (3), (5)).  □

The following consequence illustrates in a simple but useful way, the connection between critical exponents and the notion of strongly $L^p$ representations.

2.9 COROLLARY. *Let* $\Gamma < G = \mathrm{SO}(n,1), \mathrm{SU}(n,1)$ *be any discrete subgroup,* $2 \le p < \infty$, *and* $\pi$ *be a unitary* $G$-*representation which is strongly* $L^p$. *Then* $\pi|_\Gamma$ *is strongly* $L^{\frac{\delta(\Gamma)}{\delta(G)}p}$.

*Proof.* By 2 of Theorem 2.1, the $K$-finite matrix coefficients of $\pi$ satisfy the exponential decay estimate (6); thus by (12) each one is dominated by (a constant multiple of):

$$e^{(-\frac{\delta(G)}{p}+\varepsilon)\beta(\log a(g))} = e^{(-\frac{\delta(G)}{p}+\varepsilon)d(go,o)}.$$

By Definition 1.0, the function on the right-hand side, when restricted to $\Gamma$, lies in $L^q(\Gamma)$ for every $q > \frac{\delta(\Gamma)}{\delta(G)}p$, once $\varepsilon > 0$ is chosen sufficiently small. $\qquad\square$

Another relation between critical exponents and strongly $L^p$ representations is obtained by:

2.10 THEOREM. *Let* $\Gamma < G$ *be any discrete torsion-free subgroup. Then the regular* $G$-*representation on* $L^2(G/\Gamma)$ *is strongly* $L^p$ *for*

$$p = \max\{2, \delta(G)/\delta(G) - \delta(\Gamma)\}.$$

*If* $p > 2$ *this estimate is sharp:* $L^2(G/\Gamma)$ *is not strongly* $L^q$ *for* $q < p$.

*Proof.* Let us first consider the class one spectrum of $\pi = L^2(G/\Gamma)$, particularly, the complementary series representations $\pi_\lambda$, $0 \le \lambda \le \rho$, which are in the support of the spectral measure $\mu$ of the decomposition $\pi = \int \pi_x d\mu(x)$ into irreducibles. Harish-Chandra has shown that the class one representations $\pi_x$ which occur in the decomposition, correspond (through the action of the Casimir operator) to the $L^2$-spectrum of the Laplacian $\Delta$ acting on the locally symmetric space $K\backslash G/\Gamma$. More precisely, denoting by $\pi_\lambda$ ($0 \le \lambda \le |\rho| = \delta(G)/2$) the representation $\pi_{\lambda\beta}$, then this representation lies in the support of $\mu$, if and only if $|\rho|^2 - \lambda^2 (= \delta(G)^2/4 - \lambda^2)$ lies in the spectrum of $\Delta$ (see e.g. [Cor, §4] and the references therein, for this and other related results which are indicated without proof). On the other hand, the bottom of the spectrum of $\Delta$ may be precisely computed in terms of the critical exponent $\delta(\Gamma)$: it is equal to $\delta(\Gamma)(\delta(G) - \delta(\Gamma))$ if $\delta(\Gamma) > \frac{\delta(G)}{2}$, and to $(\frac{\delta(G)}{2})^2$ if $\delta(\Gamma) \le \frac{\delta(G)}{2}$ (see also [Su2, 2.17]). In the second case we get that only principal series spherical representations can appear in $L^2(G/\Gamma)$, whereas in the first, complementary series representations do appear, but only for

$$\lambda \le \sqrt{\frac{\delta(G)^2}{4} - \delta(\Gamma)(\delta(G) - \delta(\Gamma))}\beta = \sqrt{(\delta(\Gamma) - \frac{\delta(G)}{2})^2}\beta = (\delta(\Gamma) - \frac{\delta(G)}{2})\beta.$$

Therefore, all the $K$-invariant unit vectors $v_\lambda$ in these $\pi_\lambda$'s satisfy by (10) the estimate:

$$(13) \qquad |\langle \pi_\lambda(g) v_\lambda, v_\lambda \rangle| \leq C \cdot (1 + \beta(\log a(g)) \cdot e^{(\delta(\Gamma) - \delta(G))\beta \log(a(g))}.$$

(Recall that $2\rho = \delta(G)\beta$.) If $\delta(\Gamma) \leq \delta(G)/2$, by (9) we replace in (13) $\delta(\Gamma)$ by $\delta(G)/2$.

Next, let $v$ be any $K$-invariant unit vector in $L^2(G/\Gamma)$. Then by decomposing $\pi = \int \pi_x d\mu(x)$ we deduce that $v = \int v_x d\mu(x)$, where each $v_x$ is (a multiple of) an appropriate $v_\lambda$, and thus the associated matrix coefficient satisfies the decay estimate (13). A simple direct integral computation then shows that (13) must hold also for $v$. We have thus established a uniform decay estimate for all $K$-invariant vectors in $L^2(G/\Gamma)$, and we may apply Proposition 2.5 (with $M = G/\Gamma$) to get a corresponding estimate for all $K$-finite vectors. Finally, by (5) we see that the right-hand side of (13) is in $L^p(G)$ for every $p$ such that $p(\delta(G) - \delta(\Gamma)) > \delta(G)$; that is, $L^2(G/\Gamma)$ is strongly $L^p$ for $p = \frac{\delta(G)}{\delta(G) - \delta(\Gamma)}$. In case $\delta(\Gamma) \leq \delta(G)/2$, it follows that this representation is strongly $L^2$.

The last claim of Theorem 2.10 follows from the fact that the above value of the bottom of the Laplacian spectrum is sharp, estimate (11) for the spherical functions, and 3 in Theorem 2.1. $\qquad \square$

## 3. Cohomology of unitary representations

I. *Notation and some basic properties.* Let $G$ be a locally compact group and $(\pi, \mathcal{H})$ a (continuous) unitary $G$-representation. Denote by $Z^1(G, \pi)$ the vector space of continuous $G$-*cocycles* with values in $\pi$, namely, the continuous functions $b: G \to \mathcal{H}$ satisfying:

$$(1) \qquad b(gh) = b(g) + \pi(g)b(h) \qquad \text{for all } g, h \in G.$$

The subspace of *co-boundaries*: $b(g) = \pi(g)v - v$ for some $v \in \mathcal{H}$, is denoted by $B^1(G, \pi)$.

3.1 *Definition.* The *first cohomology group* of $G$ with values in $\pi$ is the quotient space: $H^1(G, \pi) = Z^1(G, \pi)/B^1(G, \pi)$.

Given a cocycle $b \in Z^1(G, \pi)$, one can deform the linear action $\pi$ to an affine isometric action of $G$ on $\mathcal{H}$. More precisely, define an action $\rho$ of $G$ on $\mathcal{H}$ by:

$$(2) \qquad \rho(g)v = \pi(g)v + b(g).$$

Conversely, any affine isometric $G$-action on a Hilbert space $\mathcal{H}$ is easily seen to be obtained by a cocycle on the linear part of the action (see e.g. [HV, pp. 45–46] for this and other facts indicated below without proof). The following simple and basic fact will be useful:

3.2 LEMMA. *Let $\pi, b$ and $\rho$ be as above. Then $b \in B^1(G, \pi)$ if an only if the associated $G$-action $\rho$ has a fixed point.*

We shall also need the following well-known fact (see [HV, §4.3]):

3.3 LEMMA. *For any isometric group action on a Hilbert space, the following are equivalent*:

 (i) *The action has a fixed point.*
 (ii) *Every orbit of the action is* (*norm*) *bounded.*
 (iii) *There exists one* (*norm*) *bounded orbit.*

II. *Restriction of cohomology*: *Proof of Theorem* 1.12. We now come to the first main result of this section. In [Sh4] it is proved in a more general framework for any isometric action of a simple algebraic group over a local field, on a CAT(0) space. As it will be essential to us on many occasions in the sequel (see [Sh4] for other applications), we include a short proof here for the more restricted result, sufficient for our purposes.

3.4 THEOREM. *Let $G$ be a simple Lie group with finite center, locally isomorphic to either* $\mathrm{SO}(n, 1)$ *or* $\mathrm{SU}(n, 1)$ $(n \geq 2)$. *Suppose that $G$ acts continuously by isometries on a Hilbert space $\mathcal{H}$, without fixed points. Then the action is proper, namely, for every $v \in \mathcal{H}$ and a sequence $g_i \to \infty$, one has $\|g_i v\| \to \infty$. In particular, in light of 3.2, for any noncompact subgroup $\Lambda < G$ and any unitary $G$-representation $\pi$, the natural restriction map $H^1(G, \pi) \to H^1(\Lambda, \pi|_\Lambda)$ is injective.*

*Proof.* Let us first prove the result in the case $G = \mathrm{SL}_2(\mathbb{R})$. We assume that the action is not proper and show that it admits a global fixed point. Suppose then that for some $v \in \mathcal{H}$, $g_i \to \infty$, and $M < \infty$, one has $\|g_i v\| < M$, or actually, what is more convenient for us, $\|g_i v - v\| < M$. Consider the polar decomposition: $\mathrm{SL}_2(\mathbb{R}) = KA^+K$, where $K = \mathrm{SO}(2)$ and $A^+ = \mathrm{diag}(a, a^{-1}) = \{(a)\}$ with $a \geq 1$. Decomposing $g_i = k_i(a_i)k_i'$, a standard argument using the compactness of $K$ shows that we may assume without loss of generality that $g_i = (a_i) \in A^+$, with $a_i \to \infty$. Let $N$ denote the upper triangular unipotent subgroup, and fix some $n \in N$. Then $(a_i)^{-1}n(a_i) \to e$; hence for all $i$ large enough we have $\|(a_i)^{-1}n(a_i)v - v\| < 1$, and

$$\|n(a_i)v - v\| \leq \|n(a_i)v - (a_i)v\| + \|(a_i)v - v\|$$
$$= \|(a_i)^{-1}n(a_i)v - v\| + \|(a_i)v - v\| < 1 + M.$$

Next, write $\|n(a_i)v - v\| = \|(a_i)v - n^{-1}v\|$ and repeat the same argument, to get $\|v - n^{-1}v\| = \|nv - v\| < 2M + 1$. It follows that the action of $N$ on $\mathcal{H}$ has a bounded orbit and hence, by Lemma 3.3, admits a fixed point $w$.

The last part of the proof (still, for $\mathrm{SL}_2(\mathbb{R})$), is reminiscent of the so called "Mautner phenomenon": one shows that any $N$-fixed point in $\mathcal{H}$ must be fixed by $G$. Indeed, let $w$ be as above and consider the function $\varphi(g) = \|gw - w\|$. Then $\varphi$ is continuous, invariant under $N$, satisfies $\varphi(g) = \varphi(g^{-1})$ for all $g \in G$, and has the property that $\varphi(g) = 0$ if and only if $gw = w$. The argument in [Lub1, Lemma 3.1.13] shows that every such function must be identically zero on $G$; that is, $w$ is indeed fixed by $G$.

Let us return to the general case. Let $G$ be (locally isomorphic to) $\mathrm{SO}(n, 1)$ or $\mathrm{SU}(n, 1)$, and let $H < G$ be a copy of (a finite cover of) $\mathrm{SO}(2, 1) \cong \mathrm{PSL}_2(\mathbb{R})$. The above proof for $\mathrm{SL}_2(\mathbb{R})$ works for any finite cover of it, hence holds also for $H$. Recall that $H$ and $G$ share a mutual Cartan subgroup $A$. Thus, if for some sequence $g_i \to \infty$, and $v \in \mathcal{H}$, we have $\|g_i v - v\| < M$, by writing a polar decomposition $g_i = k_i a_i k_i'$, we deduce that $\|a_i v - v\|$ is bounded, and the $H$-action on $\mathcal{H}$ is not proper. By the first part of the proof, $H$ fixes a point in $\mathcal{H}$. In particular, *every* orbit of $H$ is bounded (see Lemma 3.3). Let $K$ be a maximal compact subgroup of $G$, such that the polar decomposition $G = KAK$ holds, and let $w \in \mathcal{H}$ be a point fixed by $K$. Then, since $\|aw - w\| < M$ for some $M < \infty$ and every $a \in A$, we have $\|k_1 a k_2 w - w\| = \|aw - w\| < M$, so that $\|gw - w\| < M$ for all $g \in G$. By Lemma 3.3 we conclude that $G$ admits a fixed point in $\mathcal{H}$. □

III. *The induction operation on the first cohomology.* An important notion in representation theory is that of induction. We shall now discuss this operation in relation to the first cohomology, and our main concern will be the more subtle issue of *nonuniform* lattices.

To begin however, let $G$ be any locally compact group equipped with Haar measure $\mu$, and let $\Gamma < G$ be a discrete subgroup. Given a unitary $\Gamma$-representation $(\pi, \mathcal{H})$, recall that the (unitarily) induced $G$-representation $\mathrm{Ind}_\Gamma^G \pi$ is defined as follows: The representation space for $\pi$ is the space of measurable functions $f : G \to \mathcal{H}$ satisfying the two properties:

$$f(g\gamma) = \pi(\gamma)^{-1} f(g) \quad \text{for all } \gamma \in \Gamma \text{ and a.e. } g \in G, \qquad \int_{G/\Gamma} \|f(g)\|^2 < \infty.$$

The group $G$ operates from the left: $gf(x) = f(g^{-1}x)$. An equivalent way to define this construction is as follows: Let $X$ be a fundamental domain for $G/\Gamma$, i.e. $G = X\Gamma$. Then we identify $X$ with $G/\Gamma$ and let $\alpha : G \times X \to \Gamma$ denote the associated cocycle:

$$(3) \qquad \alpha(g, x) = \gamma \quad \text{if and only if} \quad gx\gamma \in X.$$

Then the induced representation $\mathrm{Ind}_\Gamma^G \pi$ may be identified with the space $L^2(X, \mathcal{H})$, where the $G$-action is given by:

$$(4) \qquad gf(x) = \pi \circ \alpha(g^{-1}, x) f(g^{-1} \cdot x).$$

We have put the dot notation in $g^{-1} \cdot x$ to emphasize that this multiplication is the $G$-action on $X \cong G/\Gamma$, but we suppress this notation henceforth.

Next, we examine the induction operation on the first cohomology. Let $b \in Z^1(\Gamma, \pi)$ and consider, for every $g \in G$, the function $\tilde{b} \colon G \to \operatorname{Ind}_\Gamma^G \pi$ defined by:

$$(5) \qquad \tilde{b}(g)(x) = b \circ \alpha(g^{-1}, x).$$

Since $\tilde{b}$ is measurable on $G$, it is in fact continuous (this will also become evident in the geometric interpretation we presently give). It is easy to verify that (5) satisfies the formal cocycle identity (1), and we shall return to it shortly; however, one first has to check the integrability condition:

$$(6) \qquad \int_X \|\tilde{b}(g)(x)\|^2 d\mu(x) < \infty \quad \text{for all } g \in G.$$

If we assume that $\Gamma$ is co-compact, then for a choice of bounded $X$, (6) is obvious. Indeed, the cocycle $\alpha$, hence also $\tilde{b}(g)$, is bounded (as $\alpha(g, x)$ takes only finitely many values for every $g \in G$). Condition (6) for nonuniform lattices will occupy our attention shortly.

Using the connection between cocycles and isometric actions (2), we may present the situation in a more geometric framework. Starting with an affine isometric action $\rho$ of $\Gamma$ on a Hilbert space $\mathcal{H}$, whose linear part is $\pi$, consider the induced isometric $G$-action on $L^2(X, \mathcal{H})$, given by:

$$(7) \qquad gf(x) = \rho \circ \alpha(g^{-1}, x) f(g^{-1}x).$$

As in (4), it is easily verified that this is indeed a $G$-action, and recalling that $\rho(\gamma)v = \pi(\gamma)v + b(\gamma)$, we see that:

$$gf(x) = \pi \circ \alpha(g^{-1}, x) f(g^{-1}x) + b \circ \alpha(g^{-1}, x).$$

Thus (5) indeed describes a cocycle $\tilde{b}$ for the induced representation $\operatorname{Ind}_\Gamma^G \pi$. The integrability condition (6) is equivalent to the condition:

$$(8) \qquad \int_X \|\rho \circ \alpha(g^{-1}, x) f(g^{-1}x)\|^2 d\mu(x) < \infty.$$

Notice that (8) does not depend on $f$, as once it holds for some function $f_1 \in L^2(X, \mathcal{H})$; the fact that $\int_X \|gf - gf_1\|^2 d\mu(x) = \int_X \|f - f_1\|^2 d\mu(x) < \infty$ shows that it holds for every other $f \in L^2(X, \mathcal{H})$. We may thus choose any test function $f$, and a natural one would be the constant function $f \equiv o$ (which immediately reduces (8) to (6)). An equivalent, but more convenient way for us to present condition (6) (or (8)) is thus:

$$(9) \qquad \int_X \|\rho \circ \alpha(g^{-1}, x)o\|^2 d\mu(x) < \infty \quad \text{for all } g \in G.$$

Suppose that $\Gamma$ is finitely generated, and let $S < \Gamma$ be a finite generating set. Let $\ell = \ell_S$ be the corresponding word length on $\Gamma$. Denote $M = \max\{\|\rho(\gamma)o\| \mid \gamma \in S\}$. Then successive iterations of the triangular inequality easily show that for any $\gamma \in \Gamma$ one has $\|\rho(\gamma)o\| \leq M \cdot \ell(\gamma)$. It follows that in order to verify condition (9), it is enough to verify:

$$(10) \qquad \int\limits_X \ell(\alpha(g,x))^2 d\mu(x) < \infty \quad \text{for all } g \in G,$$

a condition which depends now only on $\Gamma$ and $G$, and *not* on the representation $\pi$.

3.5. Before discussing (10) for nonuniform lattices, let us assume for the moment that for some $\Gamma < G$ it holds. We shall prove that this implies Theorem 1.10, by showing that the map $b \to \tilde{b}$ induces an isomorphism between the first cohomology groups. That (5) induces a (linear) map on the first cohomology groups $H^1(\Gamma, \pi)$, $H^1(G, \mathrm{Ind}_\Gamma^G \pi)$, namely, that it maps a coboundary to a coboundary, is trivial. We observe first that this map is *injective*. Indeed, by Lemma 3.2 the injectivity is equivalent to the following: If there exists an $\alpha$-invariant $L^2$-function for the $G$-isometric action, i.e., $f \in L^2(X, \mathcal{H})$ with $\rho \circ \alpha(g^{-1}, x) f(g^{-1}x) = f(x)$ for all $g \in G$ and almost every $x \in X$, then there is a $\Gamma$-fixed point for its original action on $\mathcal{H}$. This, however, is exactly what is shown in [Zi1, 4.2.19] (even though in a somewhat different context).

The issue of *surjectivity* of the map is less transparent, and here it is convenient to use the fact that $G$ is a Lie group. Since the argument is standard, and we prove the surjectivity only for completeness (we will not be needing it in the sequel), we shall be brief. Given a unitary $G$-representation, $(\pi, \mathcal{H})$, we denote by $\mathcal{H}^\infty$ the space of smooth ($C^\infty$) vectors, a space which is invariant under the action of $G$, its Lie algebra $\mathfrak{g}$, and the universal enveloping algebra $U(\mathfrak{g})$. Also, $\mathcal{H}^\infty$ is a continuous $G$-module with its $C^\infty$ topology, which is equivalent to the topology defined by the semi-norms $v \to \|Xv\|$, $X \in U(\mathfrak{g})$. It was shown by Pichaud (see [Pi, §4.1-4.2] for further related information, or [BW, p.287] for cohomology in all degrees) that the natural map $H^1(G, \mathcal{H}^\infty) \to H^1(G, \mathcal{H})$ is an isomorphism. Hence, starting with any cocycle $b \in Z^1(G, \mathrm{Ind}_\Gamma^G \pi)$, in order to show that it lies in the image of the map in question, we may assume that it takes values in $(\mathrm{Ind}_\Gamma^G \pi)^\infty$. Assume for the moment that we show that any element in the latter space is a smooth *function* on the manifold $G/\Gamma$ with values in $\mathcal{H}$. Then the map $\gamma \to b(\gamma)(\bar{e})$ makes sense, and it is easy to check that it defines a cocycle $\bar{b} : \Gamma \to \mathcal{H}$ such that $\tilde{(\bar{b})} = b$. Finally, to show that a smooth vector is indeed a smooth function, it is enough to see that it is a weakly smooth function on $G/\Gamma$, which can be shown to follow from Sobolev's Lemma (cf. [Yo, p.174]).

IV. *Induction of cohomology for nonuniform lattices*: *Proof of Theorem* 1.10. The rest of Section 3 is devoted to proving the following result for $G = \mathrm{SO}(n,1)$.

3.6 THEOREM. *If $n \geq 4$ and $\Gamma < G = \mathrm{SO}(n,1)$ is any nonuniform lattice, then condition* (10) *is satisfied for an (appropriately chosen) fundamental domain $X$. Consequently, by* 3.5 *above, for every unitary $\Gamma$-representation $\pi$, the map $b \to \tilde{b}$ (given by* (5)) *defines an isomorphism $H^1(\Gamma, \pi) \cong H^1(G, \mathrm{Ind}_\Gamma^G \pi)$.*

In the appendix, the following is proved:

3.7 THEOREM. *Theorem* 3.6 *holds if $\mathrm{SO}(n,1)$ is replaced by any other rank one simple Lie group.*

Of course, since the groups $\mathrm{Sp}(n,1)$ and $F_{4(-20)}$ have property (T), we will not be interested in them here. Although basically, the same approach used to prove Theorem 3.6 applies to the other cases, the details become more involved, and require somewhat more sophisticated tools. As the proof of Theorem 3.7 is quite technical, we postpone it to the appendix, and concentrate here only on Theorem 3.6. The relative simplification in the proof of the latter comes from the fact that only for the real hyperbolic spaces one knows a "good" *global* description of the Riemannian metric, while in the other cases a more local analysis is required.

We first make some general preliminary reductions, which will serve us in the proofs of both theorems. Replacing $\Gamma$ by a finite index subgroup, we may assume that $\Gamma$ is torsion free. Let $K < G$ be a maximal compact subgroup and choose a fundamental domain $Y$ for the $\Gamma$ (right) action on $G$, which is $K$-invariant. Then for every $k \in K$, one has $\alpha(gk, y) = \alpha(g, ky)$, and since $y \to ky$ is measure preserving, we deduce that the integral expression in (10) is invariant under the map $g \to gk$. Furthermore, if $k \in K$, $g \in G$, $\gamma \in \Gamma$ and $y \in Y$ satisfy $gy\gamma \in Y$, then also $kgy\gamma \in Y$, so that $\alpha(kg, y) = \alpha(g, y)$. It follows that (10) is also invariant under left $K$-multiplication. Hence by polar decomposition, it is enough to establish (10) when $g = a_t \in A$ (see (1) in §2.I).

Let $M$ be the centralizer of $A$ in $K$. Since $A$ commutes with $M$, its (left) action on $G/\Gamma$ induces an action on $M \backslash G / \Gamma$. Recall that the fundamental domain $Y$ is $K$-invariant, and in particular $M$-invariant. For all $y \in Y$, $m \in M, a_t \in A, \gamma \in \Gamma$ we have $a_t my\gamma = ma_t y\gamma$, so that $\alpha(a_t, my) = \alpha(a_t, y)$. Thus, as far as the $A$-action is concerned, we may safely replace $G$ throughout the following by $M \backslash G$ and denote, still by $Y$, a fundamental domain for the (right) $\Gamma$-action on the latter. Recall now that $M \backslash G$ may be identified with the unit tangent bundle of the symmetric space $K \backslash G$. With this identification, the left $A$-action on $M \backslash G$ is no more than the geodesic flow on $T^1(K \backslash G)$. Denoting by $p \colon M \backslash G \to K \backslash G$ the natural projection, and by $d$ the usual $G$-invariant Riemannian metric on $K \backslash G$, we thus have (with an appropriate normalization): $d(p(a_t y), p(y)) = t$. The map $p$ is $\Gamma$-equivariant; hence $F = p(Y)$ is a fundamental domain for the right $\Gamma$-action on $K \backslash G$, and for $\gamma \in \Gamma$, we have $a_t y \gamma \in Y$ if and only if $p(a_t y)\gamma \in F$. Now, notice that the set of $g \in G$, for which condition (10) is satisfied, is closed under multiplication (indeed, with notation as in the beginning of §3.III, write $\|ghf\| \le \|ghf - gf\| + \|gf\| = \|hf - f\| + \|gf\| \le \|f\| + \|hf\| + \|gf\|$). Thus, we may consider only $a_t$ for $-1 \le t \le 1$, and together with the foregoing discussion, we are now reduced to establishing the following result:

3.8 THEOREM. *With the notation above, for every $x \in K \backslash G$ let $w(x)$ denote the word length of the element $\gamma \in \Gamma$ (with respect to some fixed generating set of $\Gamma$), for which $x\gamma \in F$. The function $w$ takes only finitely many values on any compact subset $C$ (since $C$ meets only finitely many $\Gamma$-translates of $F$), and let $f(p) = \sup\{w(x) \mid d(x, p) \le 1\}$. Then $f|_F \in L^2(F)$ with respect to the Riemannian volume on $F \subset K \backslash G$.*

As mentioned, we discuss here only $G = SO(n, 1)$, and in the appendix also $G = \mathrm{SU}(n, 1), \mathrm{Sp}(n, 1)$. First, we need to recall some fundamental results of Garland-Raghunathan [GR] concerning the structure of some fundamental domains for nonuniform lattices in rank one Lie groups.

Let $\Gamma$ be such a lattice. Then a fundamental domain for the $\Gamma$-action on $K \backslash G$ can always be chosen to consist of finitely many cusps. Since their number is finite, it is enough to show the $L^2$-integrability of $f$ on each one, so we confine our discussion hereafter to studying the behavior of $f$ on one such cusp $\Omega$. The direction of the cusp $\Omega$ is determined by one geodesic ray, which without loss of generality may assumed to be given by $oA^+ = o(a_t)$, where $o = \bar{e} \in K \backslash G$, and $a_t$ is as in (1) of Section 2. The cusp group, namely, the group which fixes this geodesic ray, is the group $P = MAN$ (see 2.I). The group $N \cap \Gamma$ is cocompact in $N$.

Let us specialize our discussion henceforth to $G = \mathrm{SO}(n, 1)$. Then we have $N \cong \mathbb{R}^{n-1}$, so $N \cap \Gamma \cong \mathbb{Z}^{n-1}$. Choosing a bounded fundamental domain $T(\cong [0, 1)^{n-1})$ for $N \cap \Gamma$ in $N$, we have $\Omega = o \cdot (a_t)_{t \ge 0} T$ (see [GR, Th. 0.6]). We shall now need to recall the more concrete (and well-known) description of the real hyperbolic space $\mathbb{H}^n \cong K \backslash G$. Recall that in the upper half space model, $\mathbb{H}^n$ is represented by $n$ coordinates $(x_1, \cdots, x_{n-1}, e^t)$, with $x_i, t \in \mathbb{R}$, and the hyperbolic length element is given by:

$$(11) \qquad ds^2 = \frac{dx_1^2 + \cdots + dx_{n-1}^2 + e^{2t}dt^2}{e^{2t}}.$$

These coordinates are induced by the Iwasawa decomposition, $G = KNA$, where $N = \mathbb{R}^{n-1} = \{x_1, \cdots, x_{n+1}, 1\}$ and

$\qquad A \cong \mathbb{R}^* = (0, 0, \cdots, e^{-t}) \;\; (o = (0, ...0, 1)).$

The $G$-Haar measure in the Iwasawa coordinates is given by (see [GV, 2.4.2]): $dg = e^{2\rho(\log a)} dk dn da$, where $dk$ is the $K$-Haar measure, $dn = dx_1 \cdots dx_{n-1}$ is the $\mathbb{R}^{n-1}$-Haar measure, and $da$ is the Haar measure of $A$, which, when written additively through the exponent change of variable, is just the Lebesgue measure $dt$. Since $\rho(\log a) = (n-1)t$, we have for the volume element on $K\backslash G$:

$$(12) \qquad\qquad dv = \frac{dx_1 \cdots dx_{n-1} dt}{e^{(n-1)t}}.$$

For a point $y = (x_1, \cdots, x_{n-1}, e^t) \in \mathbb{H}^n$ call $t = h(y)$ the *height* of $y$.

Let us now return to the proof of Theorem 3.8, concentrating on the cusp $\Omega$ above. It is easy to see that if $p \in \Omega$ has height large enough, then for any point $x$ with $d(p,x) \leq 1$ the element $\gamma \in \Gamma$ such that $x\gamma \in F$ actually satisfies $x\gamma \in \Omega$, and belongs to $N$. Since we may assume that the generating set $S \subset \Gamma$ contains the standard generators for the lattice $\Gamma \cap N \cong \mathbb{Z}^{n-1}$, (11) shows that a bound on the difference between the first $n-1$ coordinates of $p$ and $x$ yields a bound on $l(\gamma)$. This and (11) imply that $f(p)$ (defined in 3.8) is dominated by some constant multiple of the square root of the denominator of (11) at $p$, namely, by $Ce^{h(p)}$. Therefore, using (12) we get:

$$\int_{\Omega} f^2(p) dv \leq C \int_{T \times A^+o} e^{2h(p)} dv \leq \tilde{C} \int_0^{\infty} e^{2t} e^{-(n-1)t} dt = \tilde{C} \int_0^{\infty} e^{(3-n)t} dt.$$

For $n \geq 4$, the last integral converges. This completes the proof of Theorem 3.6. $\qquad\qquad\qquad\qquad\qquad\qquad\qquad\qquad\qquad\qquad\qquad\square$

## 4. The invariants $p(G)$ and $p(\Gamma)$

**I. $p(G)$: *Proof of Theorem* 1.9 *for simple Lie groups.*** Although our interest in the present paper lies in the first cohomology, one may consider higher cohomology groups, which are of interest in their own right, and have some applications in the spirit of those presented here (these will be discussed elsewhere.) We refer the reader to [Gu] for the general definition of the cohomology groups $H^n(G, V)$, for any locally compact group $G$, and a $G$-module $V$. One can thus extend Definition 1.8 of the introduction as follows:

4.1 *Definition.* Let $G$ be a locally compact group and $n \geq 1$. Define

$$p_n(G) = \inf\{1 \leq p \leq \infty | \text{ there exists a unitary G-representation } \pi$$
$$\text{which is strongly } L^p, \text{ and satisfies } H^n(G, \pi) \neq 0\}.$$

As in Definition 1.8, if there is no $\pi$ with $H^n(G, \pi) \neq 0$, set $p_n(G) = \infty$. For consistency with Definition 1.8, we set $p(G) = p_1(G)$. The main result of this subsection, and one which is fundamental in the sequel, is the following:

4.2 THEOREM. *Suppose that $G$ is a connected simple Lie group with finite center, locally isomorphic to* $\mathrm{SO}(n, 1)$ *with* $n \geq 3$, *or to* $\mathrm{SU}(n, 1)$ *with* $n \geq 2$. *Then* $p(G) = \delta(G)$ *(see* (3) *in* §2.I *for the numerical values).*

The proof of Theorem 4.2 is based on the classification of relative Lie algebra cohomology of $(\mathfrak{g}, K)$-modules, and its scheme may be used to compute (or at least estimate) $p_n(G)$ for any simple Lie group $G$ and $n \geq 1$. Following the proof of the theorem we shall remark on the general case, and some fundamental differences which arise when $n > 1$. We shall also indicate a simpler, but more restricted, approach, suitable for the first cohomology and the groups $\mathrm{SO}(n, 1)$ only, avoiding completely the use of $(\mathfrak{g}, K)$-modules and a classification.

*Proof of Theorem* 4.2. We assume that $G$ is not locally isomorphic to

$$\mathrm{SL}_2(\mathbb{R})(\cong \mathrm{SO}(2, 1) \cong \mathrm{SU}(1, 1)),$$

a group which is discussed in Theorem 4.3 below. The proof is carried out in several stages.

(i) *Classifying the irreducible representations $\pi$ with nontrivial cohomology.* The first step is the correspondence between continuous and relative Lie algebra cohomology for the associated $(\mathfrak{g}, K)$ module of $K$-finite vectors (carried out through the smooth cohomology); see [BW] or [Gu]. Once transferred to this framework, one uses the classification of (irreducible) $(\mathfrak{g}, K)$-modules with nontrivial cohomology. It turns out that $\mathrm{SO}(n, 1)$, $n \geq 3$, has exactly one irreducible representation with nontrivial first cohomology, and $\mathrm{SU}(n, 1)$ has two such representations (conjugated under an outer automorphism of the group). This result is due to Delorme (see [De] and the references therein), and Hotta-Wallach [HW], but explicit constructions of the representations were known previously, at the Hilbert space level (see [GGV]). The classification in terms of Langlands parameters is, as we shall see, adequate for our purposes; see [BW, Ch.VI §4] for a complete discussion of $\mathrm{SO}(n, 1), \mathrm{SU}(n, 1)$ (or p.386 therein for a quick survey).

(ii) $L^p$-*integrability of $K$-finite matrix coefficients for the representations in* (i). Once the Langlands parameters of an irreducible representation of a simple Lie group are known, one can compute the exponential decay along the Cartan subgroup $A$, of the $K$-finite matrix coefficients. Since we need to identify precisely the minimal $p$ for which these matrix coefficients are in $L^{p+\epsilon}$, and not just a bound, more should be said.

If $G$ is locally isomorphic to $\mathrm{SL}_2(\mathbb{C})$ then (i) yields exactly one representation (see below), which turns out to be tempered; its $K$-finite vectors have matrix coefficients in $L^{2+\epsilon}(G)$ for all $\epsilon > 0$, but since it is not from the discrete series (which does not exist for $\mathrm{SL}_2(\mathbb{C})$), no matrix coefficient can be square integrable ([HT, Ch. V, 1.2.4]), so that $p = 2$ is indeed minimal. For the other groups, the classification shows that the representations found in (i) are not tempered. The Langlands parameters (cf. [BW, Ch.VI Theorems 3.2, 4.5, 4.11] for $\mathrm{SO}(2k+1,1), \mathrm{SO}(2k,1), \mathrm{SU}(n,1)$, resp.) yield for any such representation $\pi$, and $K$-finite vectors $v, w$, the existence of a limit:

$$(1) \qquad \lim_{a \to \infty, \, a \in A^+} e^{\beta(\log a)} |\langle \pi(a)v, w \rangle| = |L(v, w)| < \infty$$

where $L(v, w)$ is not always zero (it is a matrix coefficient associated with the representation of the compact subgroup $M$ involved in the Langlands parameters of $\pi$); see [BW, p.131] (and pp. 10, 126 therein for the notation), or [Kn, p.198, Lemma 7.23]. Notice also that dim $A = 1$ is helpful here, for the clean expression in (1). By a standard argument (see (2) in §4.III below), one shows, using the $K$-finiteness of $v, w$, that there is $C < \infty$ such that the expression inside the limit (1) is bounded by $C$ if $\pi(a)$ is replaced by $\pi(g)$, and $\log a$ by $\log a(g)$ (defined in (2) of §2.I), for any $g \in G$. This, together with the expression for the Haar measure in terms of polar coordinates (see (5) of §2), shows that for $p = \delta(G)$, the $K$-finite matrix coefficients are all in $L^{p+\epsilon}(G)$ for all $\epsilon > 0$. We should perhaps remind the reader that the representations obtained in the form of Langlands parameters are *not* unitary in general, only (irreducible and) admissible (in our case, they will be unitary exactly for $G = \mathrm{SL}_2(\mathbb{C})$). However, the ones we discuss are infinitesimally equivalent to a unitary (irreducible) representation, and infinitesimally equivalent admissible representations have the same *set* of matrix coefficients, which is all we are interested in here (see [Kn, p.209, p.211 Cor. 8.8]).

Now, to show minimality of $p$, notice that if there exists a $K$-finite matrix coefficient in $L^r(G)$ for some $r < \delta(G)$, then applying Theorem 2.1 (2) and then [Co2, p.157 Cor. 2.2.4], yields a contradiction to the nonidentically vanishing of $L(v, w)$ in (1). This completes step (ii). We only remark that a more subtle question is whether the above matrix coefficients are *not* in $L^{\delta(G)}$, but only in $L^{\delta(G)+\epsilon}$. This is indeed the case, and would refine our geometric applications if it could be shown to hold in the more general situation discussed in step (iv) below (compare with the proof of the second part of Theorem 1.5, in §5.II, where this refinement holds and is used). However for our purposes here, this issue does not arise, as the definition of $p(G)$ considers the infimum of $p$'s, which will be $\delta(G)$ regardless of such a refinement.

(iii) $L^p$-*integrability of all matrix coefficients for the representations in* (i). Recall that since we are still dealing with irreducible representations, having one or a dense set of $L^p$-matrix coefficients is equivalent. By Theorem 2.1 it follows that the same $p$ computed in step (ii) is the minimal value for the $L^p$ integrability of *all* the matrix coefficients (there might be matrix coefficients lying in $L^{p(G)}$ rather that in $L^{p(G)+\epsilon}$ for all $\epsilon > 0$, but as remarked in (ii) above, this has no bearing on the computation of $p(G)$).

(iv) *Passing from irreducible to general representations.* If $G$ is a locally compact group and $\pi = \int^{\oplus} \pi_x$ is a direct integral decomposition into irreducibles, it is in general possible that $\pi$ is strongly $L^p$ for some $p < \infty$, but none of the $\pi_x$ have this property (e.g., $G = \mathbb{R}$, $\pi = L^2(G)$). The situation with simple Lie groups is different, and from Theorem 2.1 we moreover get the sharp result required for the groups of interest to us; as long as $p \geq 2$, almost every $\pi_x$ is strongly $L^p$ for the same value of $p$.

Next, consider a direct integral decomposition $\pi = \int^{\oplus} \pi_x$, and suppose that $H^1(G, \pi) \neq 0$. It is not true in general that for a set $x$ of positive measure one must have $H^1(G, \pi_x) \neq 0$. However this is true, by a result of Guichardet, if we replace $H^1$ with the reduced cohomology $\overline{H^1}$ (cf. [Gu, Ch.3, §2] for the precise notions and results, and also the proof after Theorem 1.7). It is easy to see that $H^1$ coincides with $\overline{H^1}$ if the representation $\pi$ does not weakly contain the trivial representation (cf. [Sh2, 1.6]). However, every strongly $L^p$ representation ($p < \infty$) of a nonamenable group has the latter property. Briefly, this follows by examination of the $n$-fold tensor power $\sigma = \pi^n$ for some integer $n > p/2$. Using Hölder inequality, one shows that $\sigma$ is strongly $L^2$; hence by [CHH] it is weakly contained in the regular representation. But if $\pi$ contains weakly the trivial representation, then so does $\sigma$, which together with the previous argument shows that the same holds for the regular representation, contradicting nonamenability (see e.g. the proof of Theorem C in [Sh3]). We can now deduce from this and the foregoing discussion, that if $\pi = \int^{\oplus} \pi_x$ is a representation of $G = \mathrm{SO}(n, 1)$ or $\mathrm{SU}(n, 1)$, which is both strongly $L^p$ ($2 \leq p < \infty$) and satisfies $H^1(G, \pi) \neq 0$, then there is an *irreducible* $\pi_x$ with the same properties; so by (iii) we have $p \geq \delta(G)$. This completes the proof of Theorem 4.2. $\qquad\square$

Let us say a few words concerning the situation for a general simple Lie group with finite center $G$, and all values of $n$, when estimating the value of $p_n(G)$ in Definition 4.1. The general problem of describing all unitary representations with nontrivial cohomology (in all degrees), for any simple Lie group $G$, is a well-studied one. It is known that for any $G$, only irreducible $\pi$ for which the center of the universal enveloping algebra acts as in the trivial representation, may have $H^n(G, \pi) \neq 0$ for some $n$, and that these are finite in number. Vogan and Zuckerman [VZ] established a general (but explicit) method of describing all such unitary representations, for any $G$ (although the unitarity of all their constructions was proved only later by Vogan). In [VZ] it is also shown how to locate the constructions in the Langlands classifica-tion. Thus, the problem of classification is as much solved as one can expect, in complete generality, hence step (i) may be carried out in general. Next, as mentioned earlier, once the Langlands parameters are known, it is an easy matter to determine the exact decay of $K$-finite matrix coefficients along the Cartan subgroup, thereby getting a sharp value for $p$ in step (ii). Steps (iii) and (iv) reveal the fundamental differences which occur in this more general framework. Due to the lack of (knowledge of) Theorem 2.1 for other $G$'s, we can only assert in general, using Cowling's "Kunze-Stein phenomenon" [Co1], that the value of $p$ may decrease up to the closest *even* integer in the process of passing to all matrix coefficients and all representations, thereby giving only an interval of length (at most) 2 for the value of $p_n(G)$. Moreover, in the case $n \geq 2$, to pass from irreducible to general representations in (iv), one has to work originally with reduced cohomology, rather than the ordinary one (which requires more effort in a later stage, when transferring this property to lattices), or establish a general criterion, as the one we had for $H^1$, relating the two notions. This strategy seems possible at least when $n = \mathrm{Rank}(G)$, however, we shall not pursue here further this direction.

Before presenting a different approach to proving Theorem 4.2, let us consider the general case $p(G) = 0$, and in particular that of $G = \mathrm{SL}_2(\mathbb{R})$.

4.3 THEOREM. *For a locally compact group $G$, the following properties are equivalent*:

(1) $p(G) = 0$;

(2) $p(G) < 2$;

(3) $H^1(G, L^2(G)) \neq 0$.

*If either $G = \mathrm{SL}_2(\mathbb{R})$, or $G$ is amenable and noncompact, then* (1)–(3) *are satisfied. When $G$ is a finitely generated nonamenable group,* (1)–(3) *are equivalent to nonvanishing of the first $\ell^2$-Betti number.*

*Proof.* Clearly $(1) \Rightarrow (2)$ and $(3) \Rightarrow (1)$ (the regular representation has a dense set of matrix coefficients with compact support). The fact that $(2) \Rightarrow (3)$ follows from a well-known result (cf. [HT, Ch.V, 1.2.4]): any unitary representation with a dense subspace of matrix coefficients which are square integrable, can be embedded in a multiple of the regular representation. We only remark that because matrix coefficients are bounded, $L^p$-integrability for $p < 2$ implies square integrability, and that one has $H^1(G, \infty \cdot \pi) \neq 0$ if and only if $H^1(G, \pi) \neq 0$.

Since $\mathrm{SL}_2(\mathbb{R})$ has (two) discrete series representations with nontrivial first cohomology ([BW, Ch.VI]), (3) is satisfied. It is also satisfied by any amenable noncompact group, by a result of Guichardet [HV, p.48] (every representation which contains weakly, but not properly, the trivial representation, has non-vanishing first cohomology). The last assertion is due to Bekka and Valette [BV]. $\qquad\square$

We now turn to discuss a more illuminating approach to the proof of Theorem 4.2 for $G = \mathrm{SO}(n, 1)$, explaining better the connection between $p(G)$ and $\delta(G)$. We shall avoid classification, and give explicit geometric constructions of representations achieving the bound $p(G)$. Some of the results established along the way will also be used in the sequel.

4.4 PROPOSITION. *For $G = \mathrm{SO}(n, 1)$, $p(G) \leq \delta(G)$.*

*Proof.* By Theorem 4.3 we may assume $n \geq 3$. Recall that many lattices (uniform and nonuniform) $\Gamma < \mathrm{SO}(n, 1)$ admit an amalgam decomposition $\Gamma = A *_C B$, where $C$ is a lattice in $\mathrm{SO}(n - 1, 1)$. Constructions of such decompositions go back to Millson [Mi]. Explicit and well-known ones appear in [GPS], where the nonarithmetic lattices share this property virtually by their construction. It can be shown (see the discussion in §5.I below), that for any group $\Gamma$ and a decomposition $\Gamma = A *_C B$, one has $H^1(\Gamma, \ell^2(\Gamma/\tilde{C})) \neq 0$, where $\tilde{C} < C$ has index at most 2. Thus, taking a uniform lattice $\Gamma$, with an amalgam decomposition as above, we may induce the $\Gamma$-representation $\ell^2(\Gamma/C)$ to $G$, and from 3.5 deduce that $H^1(G, L^2(G/C)) \neq 0$. However, recalling that $C$ is a lattice in $\mathrm{SO}(n - 1, 1)$, we have $\delta(C) = n - 1 - 1 = n - 2$. From Theorem 2.10 it follows that $L^2(G/C)$ is strongly $L^{\delta(G)}$; hence by definition, $p(G) \leq \delta(G)$. $\qquad\square$

4.5 THEOREM. *Suppose that $G = \mathrm{SO}(n, 1)$ or $\mathrm{SU}(n, 1)$, and that $H$ is either a group of this type, embedded naturally in $G$ (as described at the beginning of Section 2), or any infinite discrete subgroup. Then $p(H) \leq \frac{\delta(H)}{\delta(G)} p(G)$.*

*Proof.* Assume first that $p(G) \geq 2$ and let $\pi$ be a $G$-representation which is strongly $L^p$ for $p \geq 2$. Then $\pi|_H$ is strongly $L^{\frac{\delta(H)}{\delta(G)} p}$ (if $H$ is discrete use Corollary 2.9; in the other case use Theorem 2.1 and expression (5) for the Haar measure). The result then follows from Theorem 3.4. If $p(G) < 2$, then using Theorem 4.3 and the fact that $L^2(G)|_H$ is a multiple of $L^2(H)$, we get, applying Theorem 3.4 again, that $p(H) = p(G) = 0$, as required. $\qquad\square$

*Remark.* A posteriori, once we have established $p(G) = \delta(G)$ (or at least $p(G) \leq \delta(G)$), it follows from Theorem 4.5 that $p(\Lambda) \leq \delta(\Lambda)$ for every infinite discrete subgroup $\Lambda < \mathrm{SO}(n, 1), \mathrm{SU}(n, 1)$.

4.6 LEMMA. *If $G = \mathrm{SO}(3, 1)(\cong \mathrm{SL}_2(\mathbb{C}))$ and $\Gamma < G$ is a lattice, then $p(G) \geq 2$ and $p(\Gamma) \geq 2$.*

*Proof.* If the conclusion fails either for $G$ or for $\Gamma$, we have by Theorem 4.3: $H^1(G, L^2(G)) \neq 0$, $H^1(\Gamma, \ell^2(\Gamma)) \neq 0$, resp. The first case is reduced to the second by restricting the regular representation (and using, say, Theorem 3.4). However the second is impossible by [BV, Th. D], if we go back to the vanishing of the first $\ell^2$-Betti number for these lattices. $\qquad\square$

4.7 COROLLARY. *For $G = \mathrm{SO}(n, 1)$, $n \geq 3$, $p(G) \geq \delta(G)$. Hence, by Proposition 4.4, $p(G) = \delta(G)$.*

*Proof.* If $n \geq 3$ and $p(G) < \delta(G)$, then with $H = \mathrm{SO}(3, 1)$ in Theorem 4.5 we get $p(H) < \delta(H) = 2$, which contradicts Lemma 4.6. $\qquad\square$

4.8 COROLLARY. *For $G = \mathrm{SU}(n, 1)$, $n \geq 3$, $p(G) \geq \delta(G)$.*

*Proof.* The argument is similar to that in 4.7, with $H = \mathrm{SO}(3, 1) < \mathrm{SU}(n, 1)$. $\qquad\square$

The difference between the families $\mathrm{SO}(n, 1)$ and $\mathrm{SU}(n, 1)$, is that we are not aware of any natural construction of representations of the latter, as in Proposition 4.4, for which one easily verifies $p(G) = \delta(G)$.

II. $p(\Gamma)$: *Proof of Theorem* 1.9 *for lattices.* The proof of the following result illustrates why it was essential to define $p(G)$, when $G$ is a simple Lie group, for an arbitrary dense subspace of matrix coefficients, and not restrict to $K$-finite ones.

4.10 LEMMA. *Let $G$ be a locally compact group and $\Gamma < G$ a lattice. If $\pi$ is a unitary $\Gamma$-representation which is strongly $L^p$ for $p \geq 1$, then $\mathrm{Ind}_\Gamma^G \pi$ is also strongly $L^p$.*

*Proof.* Let $\mu$ be the Haar measure of $G$ and $X$ be a fundamental domain for $\Gamma$ in $G$: $G = X\Gamma$. We normalize $\mu$ such that $\mu(X) = 1$ and denote its restriction to $X$ by $\mu_X$. Let $\alpha \colon G \times X \to \Gamma$ be the associated cocycle (see (3) in Section 3 – we shall use freely the notation of that section). Let $\mathcal{H}$ be the representation space for $\pi$, and $\mathcal{H}_0 \subset \mathcal{H}$ the assumed dense subspace. It suffices to show that if $\varphi, \psi \colon X \to \mathcal{H}_0$ take only finitely many values, then $\langle g\varphi, \psi \rangle \in L^{p+\varepsilon}(G)$ for all $\varepsilon > 0$, as the subspace of those vectors is clearly dense.

Let $F = \{u_i\} \subset H_0$ be the (finite) set of vectors in the image of $\varphi$ and $\psi$, and define for every $\gamma \in \Gamma$: $f(\gamma) = \sum_{i,j} |\langle \pi(\gamma)u_i, u_j \rangle|$. By our assumption on $\pi$ it is clear that $f(\gamma) \in \ell^{p+\varepsilon}(\Gamma)$ for all $\varepsilon > 0$. Next, consider the measure $\lambda$ on $\Gamma$ pushed by $\alpha$ from $G \times X$: $\lambda = \alpha_*(\mu \times \mu_X)$. In other words, $\lambda\{\gamma\}$ is the $\mu \times \mu_X$ measure of all the pairs $(g, x) \in G \times X$ for which $\alpha(g, x) = \gamma$. We claim that $\lambda\{\gamma\} = \alpha_*(\mu \times \mu_X)\{\gamma\} = 1$, namely, that $\lambda$ is just the counting measure on $\Gamma$. To verify this it is enough, by Fubini, to show that for every fixed $x_0 \in X$, the set of $g \in G$ satisfying $\alpha(g, x_0) = \gamma$ has measure 1. Indeed, $\alpha(g, x_0) = \gamma \Leftrightarrow gx_0\gamma \in X \Leftrightarrow g \in X\gamma^{-1}x_0^{-1}$. We can now compute:

$$\int_G |\langle g\varphi, \psi \rangle|^{p+\varepsilon} d\mu(g) = \int_G |\int_X \langle \alpha(g^{-1}, x)\varphi(g^{-1}x), \psi(x) \rangle d\mu_X(x)|^{p+\varepsilon} d\mu(g)$$

$$\leq \int_G (\int_X |\langle \alpha(g^{-1}, x)\varphi(g^{-1}x), \psi(x) \rangle| d\mu_X(x))^{p+\varepsilon} d\mu(g)$$

$$\leq \int_G \int_X |\langle \alpha(g^{-1}, x)\varphi(g^{-1}x), \psi(x) \rangle|^{p+\varepsilon} d\mu_X(x) d\mu(g).$$

Making a change of variables: $\gamma = \alpha(g^{-1}, x)$, and recalling the definition of $f$, we may continue to complete the proof:

$$\leq \sum_{\gamma \in \Gamma} f(\gamma)^{p+\varepsilon} d\lambda(\gamma) < \infty. \qquad \square$$

Notice that for all discrete subgroups $\Gamma < G = \mathrm{SO}(n, 1), \mathrm{SU}(n, 1)$, we have established an inequality $p(\Gamma) \leq \delta(\Gamma) \leq \delta(G) = p(G)$ (by Theorems 4.2 and 4.5). We can now show that for lattices equality always holds.

4.11 THEOREM. *Let $G$ be a locally compact group and $\Gamma < G$ a uniform lattice. Then for every $n \in \mathbb{N}$, $p_n(\Gamma) = p_n(G)$. If $\Gamma$ is nonuniform, then $p(\Gamma) \leq p(G)$, and when $G$ is a connected simple Lie group, $p(\Gamma) = p(G)$ as well.*

*Proof.* For the inequality $p_n(\Gamma) \leq p_n(G)$ ($\Gamma$ uniform), one only needs to verify that the restriction of a strongly $L^p$ representation from $G$ to $\Gamma$ is still strongly $L^p$ (this holds for any closed subgroup $\Gamma$; see [Ho, 6.4]), and that for a representation $\pi$ with $H^n(G, \pi) \neq 0$, one has $H^n(\Gamma, \pi|_\Gamma) \neq 0$. For the second claim, notice that for every uniform lattice $\Gamma$ in a locally compact group $G$, (2) in Theorem 1.10 holds (cf. [Gu, III §4]). Since for every $G$-representation $\pi$, one has $\mathrm{Ind}_\Gamma^G \pi|_\Gamma \cong \pi \otimes L^2(G/\Gamma)$, and $1_G \subset L^2(G/\Gamma)$, we have $\pi \cong \pi \otimes 1_G \subset \mathrm{Ind}_\Gamma^G \pi|_\Gamma$; hence by (2) we cannot have $H^n(\Gamma, \pi|_\Gamma) = 0$. To show the opposite inequality $p_n(G) \leq p_n(\Gamma)$, use again (2) from 1.10, and Lemma 4.10.

If $n = 1$ and $\Gamma$ is nonuniform, we can repeat the argument above to show that $p(\Gamma) \leq p(G)$. Just replace the use of (3) by [HV, 3.c.19] to show nonvanishing. Finally, if $\Gamma$ is nonuniform and $G$ is a connected simple Lie group, the cases $\mathrm{SO}(n,1)$ with $n \geq 4$ and $\mathrm{SU}(n,1)$ with $n \geq 2$, follow by inducing the cohomological representation as above, by Theorems 3.6, 3.7. The case of $\mathrm{SO}(3,1) \cong \mathrm{SL}_2(\mathbb{C})$ follows from Lemma 4.6, and that of $\mathrm{SO}(2,1) \cong \mathrm{SL}_2(\mathbb{R})$ is trivial, as by Theorem 4.3 and the above discussion, we have $0 \leq p(\Gamma) \leq p(\mathrm{SL}_2(\mathbb{R})) = 0$. Since all other simple Lie groups have property (T), we have for them $p(\Gamma) = p(G) = \infty$. $\qquad\blacksquare$

III. *First applications*: *Proofs of Theorems* 1.1 *and* 1.5. We are now able to pick up the first fruits of our consideration of unitary representations and their cohomology.

*Proof of Theorem* 1.1. As $\Lambda$ and $\Gamma$ are isomorphic, $p(\Lambda) = p(\Gamma) = \delta(\Gamma)$ (by Theorems 4.11 and 4.2). Since by Theorem 4.5 we have $p(\Lambda) \leq \delta(\Lambda)$, the result follows. $\qquad\blacksquare$

Next, we prove Theorem 1.5 stated in the introduction, and discuss some further applications.

*Proof of Theorem* 1.5. Assume first that $\Lambda < G = \mathrm{SO}(n,1), \mathrm{SU}(n,1)$ is a nonamenable subgroup with $\delta(\Lambda) < 2$. From the remark after Theorem 4.5, it follows that $p(\Lambda) < 2$; hence by Theorem 4.3 the result follows.

Assume now that $\delta(\Lambda) = 2$. As in the previous case, it is enough to show that $H^1(\Lambda, \ell^2(\Lambda)) \neq 0$, when we assume the Poincaré series in Definition 1.0 converges at $s = 2$. By enlarging $G$ if necessary, we may assume that it is not $\mathrm{SO}(k,1)$ for $k = 2$ or $k = 3$. Recall now from the proof of Theorem 4.2 (see step (ii)) that $G$ admits an irreducible unitary representation $\pi$ with $H^1(G,\pi) \neq 0$, such that $\pi$ is *not* tempered. Hence $K$-finite matrix coefficients in $\pi$ satisfy along $A^+$ the asymptotic behavior described in (1) there. Let us extend the decay estimate given by (1) to all $G$ in terms of the Cartan decomposition (this will become an upper bound, rather than precise asymptotics).

Let $v, w$ be $K$-finite unit vectors, and let $v_1, ..., v_n$ $(w_1, ..., w_m)$ be an orthonormal basis for $\mathrm{Span}\{Kv\}$ $(\mathrm{Span}\{Kw\}$, resp.). As in the proof of Theorem 2.5, for every $k \in K$ there exist unique complex numbers $\theta_1(k), ..., \theta_n(k)$, such that for all $k \in K$, $\pi(k)v = \sum \theta_i(k)v_i$. Obviously, for all $i$, $|\theta_i(k)| \leq 1$. Similarly, let $\tau_1(k), ..., \tau_m(k)$ be the complex functions corresponding to $w$ and the basis $w_1, ..., w_m$. By (1) in step (ii) of Theorem 4.2 we know that there is a constant $C < \infty$ such that for all $a \in A^+$ one has $|\langle \pi(a)v_i, w_j \rangle| < Ce^{-\beta(\log a)}$, for all $i, j$. Given $g \in G$, write $g = k_1 a k_2$ where $a = a(g) \in A^+$ and $k_1, k_2 \in K$,

and evaluate:

$$
\begin{aligned}
(2) \qquad |\langle \pi(g)v, w \rangle| &= |\langle \pi(a)\pi(k_2)v, \pi(k_1^{-1})w \rangle| \\
&= |\langle \pi(a) \sum \theta_i(k_2)v_i, \sum \tau_j(k_1^{-1})w_j \rangle| \\
&= |\sum_{i,j} \theta_i(k_2)\overline{\tau_j(k_1^{-1})}\langle \pi(a)v_i, w_j \rangle| \\
&< Ce^{-\beta(\log a)} \sum_{i,j} |\theta_i(k_2)\overline{\tau_j(k_1^{-1})}| \\
&\leq Ce^{-\beta(\log a)} n \cdot m.
\end{aligned}
$$

From Lemma 2.8 and (2) above we conclude that if the Poincaré series in Definition 1.0 converges for $s = 2$, then all $K$-finite matrix coefficients of $\pi|_\Lambda$ are in $\ell^2(\Lambda)$. Since this is a dense subspace, it follows ([HT, Ch.V, 1.2.4]) that $\pi|_\Lambda$ is embeddable in a multiple of $\ell^2(\Lambda)$. However, from Theorem 3.4 we have $H^1(\Lambda, \pi|_\Lambda) \neq 0$; thus $\infty \cdot \ell^2(\Lambda)$ has nonvanishing first cohomology, and the same holds for $\ell^2(\Lambda)$ (see Theorem 4.3), as required. $\qquad \square$

As mentioned in the introduction, there are various known results on vanishing of the first $\ell^2$-Betti number. The case of Kähler groups, mentioned in the introduction, is due to Gromov [Gr1] (see [ABR] and [ABC, Ch.4] for a detailed discussion). Various other examples, including all known (and conjectured) fundamental groups of closed prime 3-manifolds, were established by Lück [Lu]. Very recently, Gaboriau [Ga] showed (sharpening a result by Lück [Lu]), that the first $\ell^2$-Betti number vanishes for every group which admits an infinite index, infinite, finitely generated normal subgroup. Together with Theorem 1.5, this yields the following purely group theoretic consequence (which is a classical theorem of Schreier in the case of free groups):

4.12 COROLLARY. *Let* $\Lambda < \mathrm{SO}(n, 1)$, $\mathrm{SU}(n, 1)$ *be a discrete, nonamenable, finitely generated subgroup, which is either of critical exponent less than* 2, *or a convergence group of critical exponent* 2. *Then every infinite, infinite index normal subgroup of* $\Lambda$, *is* not *finitely generated.*

Note that the result is tight, and fails in general for divergence groups with critical exponent 2. Indeed, many lattices $\Lambda < \mathrm{SL}_2(\mathbb{C})$ are known to surject onto $\mathbb{Z}$, with a finitely generated kernel.

Finally, notice that from [BV] it follows that the family of discrete subgroups of $\mathrm{SO}(n, 1)$ or $\mathrm{SU}(n, 1)$ with critical exponent less than 2 (or groups of convergence type with critical exponent 2), provide, by Theorem 1.5, a new class of groups admitting nonconstant Dirichlet harmonic functions (with respect to any choice of finite generating set). For instance, by Lemma 2.7 (2) every geometrically finite subgroup $\Lambda < \mathrm{SL}_2(\mathbb{C})$, which is not a lattice, satis-fies a strict inequality $\delta(\Gamma) < \delta(\mathrm{SL}_2(\mathbb{C})) = 2$, and hence the above discussion applies (see e.g. [KAG, Ch.III] for many such constructions). Fundamental groups of normal $Z^d$-covers, $d \geq 3$, of closed hyperbolic 3-manifolds, are examples of convergence groups with critical exponent 2 (see [Re]).

## 5. Isometric actions on trees

I. *Trees, spaces with walls, and cohomology.* Although our main interest lies in actions on trees, we begin by indicating how our methods apply in a more general framework. The main result, as well as the approach taken in this subsection, is not original, but follows the work of Haglund, Paulin and Valette [HPV]. Since it is of considerable importance in what follows, for completeness we include a brief introductory discussion.

The following notion was introduced by Haglund and Paulin [HP] as a tool in studying rich automorphism groups of certain CAT(0) spaces. It also unifies virtually all the known constructions of negative definite functions on groups, which arise from geometric actions.

5.1 *Definition* [HP]. Let $X$ be a (countable) set. A *wall* is a partition of $X$ into two disjoint subsets, called *half spaces*. Let $W$ be a set of walls of $X$. Then the pair $(X, W)$ is called a *space with walls*, if for all $x, y \in X$, the set $W(x, y)$ of walls separating $x$ and $y$ is finite.

Given a space with walls $(X, W)$, let $H$ be the set of half spaces determined by $W$. For $x, y \in X$ let $E_x, E_y \subset H$ be the set of half spaces containing $x, y$, resp. Then, by our assumption, the function on $H$ defined by $1_{E_x} - 1_{E_y}$ has finite support. In particular, it is in $\ell^2(H)$, the space of square integrable functions on the (discrete) set $H$. Next, any action of a group $\Gamma$ by automorphisms on $(X, W)$ induces an action on $H$, and thus a unitary representation on $\ell^2(H)$. Choosing some fixed base point $x_0 \in X$, it is easy to verify that the function $b : \Gamma \to \ell^2(H)$ defined by:

$$(2) \qquad\qquad b(\gamma) = 1_{E_{\gamma x_0}} - 1_{E_{x_0}} \in \ell^2(H)$$

belongs to $Z^1(\Gamma, \ell^2(H))$. Since $b \in B^1(\Gamma, \ell^2(H))$ if and only if $b$ is bounded (Lemma 3.3), we find that if for some point $x_0 \in X$ the function

$$(3) \qquad\qquad \gamma \to \|1_{E_{\gamma x_0}} - 1_{E_{x_0}}\|^2 = 2 \# W(\gamma x_0, x_0)$$

is unbounded on $\Gamma$, then $H^1(\Gamma, \ell^2(H)) \neq 0$. The latter condition is satisfied in most natural examples of such actions, which we illustrate with the simplest one, namely, an action on a tree.

Suppose that $X = T$ is a (simplicial) tree. Denote by $V$ (resp. $E$) the set of vertices (resp. edges). Then "cutting in the middle" an edge $e \in E$ defines in an obvious way a partition of $V$ into two half spaces, and the resulting space $H$ of half spaces may be identified with the set $\tilde{E}$ of *oriented* edges. It is clear from (3), that (2) defines an unbounded cocycle, unless the orbit $\Gamma x_0$ is bounded, in which case, as is well-known, a subgroup of index 2 in $\Gamma$ must fix some $v \in V$. Summarizing our discussion, we have the following:

5.2 PROPOSITION [HPV]. *Suppose that a group $\Gamma$ acts isometrically on a tree $T = T(V, E)$. Let $\tilde{E}$ denote the set of oriented edges. Then either $\Gamma$ fixes some vertex or edge, or the natural unitary $\Gamma$-representation on $\ell^2(\tilde{E})$ has nonvanishing first cohomology.*

Before discussing Theorem 1.6, we mention that an analysis similar to the one below for trees can be performed with other known examples of spaces with walls. These walls have typically a natural geometric interpretation in terms of the structure of the space, as in the case of trees. Amongst the known constructions we mention Coxeter complexes [BJS], CAT(0) cubical complexes [NR], constructions in [BaSa] and [Bou2], and various examples of CAT(-1) spaces in [HP] (admitting a discrete co-compact action, and a simple, locally compact, nonlinear, nondiscrete group of automorphisms). Consequently, as we presently illustrate in the case of trees, unbounded isometric actions of lattices in $\mathrm{SO}(n, 1)$ and $\mathrm{SU}(n, 1)$ on these spaces, must admit "co-dimension one" subgroups (measured by the critical exponent), which stabilize walls (themselves identified with co-dimension one geodesic hyper-planes).

II. *Proof of Theorem* 1.6. We start by proving the first assertion on the weak inequality. This result will then be used in the more refined analysis needed for the strict inequality, which will occupy most of our attention. Notice first that if $G = \mathrm{SO}(2, 1)$ then $\delta(\Gamma) = 1$ and there is nothing to prove. In case $G = \mathrm{SO}(3, 1) \cong \mathrm{SL}_2(\mathbb{C})$ our method of proof fails (at least for nonuniform lattices; this is due to the failure of Theorem 3.6 in that case), but Z. Sela has informed us that then the theorem is known, based on 3-manifolds techniques (which do not seem to generalize). Thus, we assume henceforth that $\delta(G) > 2$.

Suppose that there is no fixed vertex for the $\Gamma$-action, and consider the unitary $\Gamma$-representation on $\ell^2(\tilde{E})$, as in Proposition 5.2. Since $\tilde{E}$ is a discrete set on which $\Gamma$ acts, we have:

$$\ell^2(\tilde{E}) \cong \oplus \ell^2(\Gamma/\Gamma_e)$$

where the sum is taken over the orbits, and $\Gamma_e$ denotes a stabilizer of an (oriented) edge in the orbit. Notice that a stabilizer of an oriented edge has index at most 2 in the edge stabilizer; hence for the discussion of critical exponents we may fail to distinguish between them. Since we assume that $\Gamma$ does not fix a vertex, and the result is trivial if it fixes an edge, we deduce from Proposition 5.2 that $H^1(\Gamma, \ell^2(\tilde{E})) \neq 0$. Therefore, induction from $\Gamma$ to $G$

(and application of Theorem 1.10 when $\Gamma$ is nonuniform) yields:

$$H^1(G, \operatorname{Ind}_\Gamma^G \ell^2(\tilde{E})) = H^1(G, \oplus L^2(G/\Gamma_e)) \neq 0.$$

We distinguish now between two cases: if the above $G$-representation contains weakly the trivial representation, this means that spherical complementary series $G$-representations from any neighborhood of the trivial representation must occur there weakly, and from the discussion in the proof of Theorem 2.10 it follows that there is a sequence of edges with $\delta(\Gamma_{e_i}) \to \delta(\Gamma)$, thereby proving even more than we claim. Otherwise, it follows that for some $\Gamma_e$ one has $H^1(G, L^2(G/\Gamma_e)) \neq 0$ (see the discussion at the second part of the proof of step (iv) in Theorem 4.2). From Theorem 4.2 it then follows that $L^2(G/\Gamma_e)$ can be strongly $L^p$ only for $p \geq p(G) = \delta(G) = \delta(\Gamma)$, which by Theorem 2.10, implies that $\delta(\Gamma_e) \geq \delta(\Gamma) - 1$, as required.

The first part of Theorem 1.6 is now established. We continue our discussion through a series of results, aiming at the second part, and assume that $G \neq \mathrm{SO}(2,1), \mathrm{SO}(3,1)$.

5.3 PROPOSITION. *Let $(M, m)$ be a Borel $\sigma$-finite measure space on which $G$ acts measurably, preserving the measure $m$. Assume that the natural $G$-representation on $L^2(M)$ is strongly $L^p$ for some $2 < p < \infty$, and that it admits an irreducible subrepresentation $\sigma$ which is strongly $L^p$ for the same $p$, but* not *strongly $L^q$ for any $q < p$. Then the unique spherical complementary series representation $\pi_\lambda$ which is "exactly" strongly $L^p$, namely, the representation $\pi_\lambda$ with $p(\rho - \lambda) = 2\rho$, also embeds as a subrepresentation of $L^2(M)$.*

See Section 2.III for the notation and discussion of complementary series representations, and the relation between $\lambda$ and $p$. The relevance of Proposition 5.3 to our discussion comes via the following:

5.4 COROLLARY. *Let $G$ be as above and $C < G$ be a discrete, torsion-free subgroup, which satisfies both $\delta(C) = \delta(G) - 1$, and $H^1(G, L^2(G/C)) \neq 0$. Let $\tilde{\lambda}_0$ denote the bottom of the $L^2$-spectrum of the Laplacian $\Delta$ on the locally symmetric space $C \backslash G/K$. Then $\tilde{\lambda}_0$ appears* discretely *in the spectrum of $\Delta$; that is, there exists a nonzero $L^2$-eigenfunction of $\Delta$ with eigenvalue $\tilde{\lambda}_0$.*

*Proof of 5.4 based on 5.3.* We use freely the discussion in Sections 2.III and 2.V, regarding the connections between the critical exponent, the class one spectrum with the associated values for $L^p$-integrability, and the bottom of the $L^2$-spectrum of $\Delta$, acting in the corresponding locally symmetric space. Since $\delta(C) = \delta(\Gamma) - 1$, the representation $L^2(G/C)$ is strongly $L^p$ for $p = \delta(G)$. On the other hand, as $H^1(G, L^2(G/C)) \neq 0$, and $G$ has only finitely many irreducible representations with nonvanishing first cohomology, we deduce, using the fact that $L^2(G/C)$ does not admit almost invariant vectors, that one of these irreducible representations must appear *discretely* in $L^2(G/C)$; see the second part of the proof of step (iv) in Theorem 4.2. However, in step (ii) of the proof of 4.2 we saw that these irreducible representations are strongly $L^{\delta(G)}$ as well, and that $p = \delta(G)$ is the minimal value for which they are strongly $L^p$. Therefore from Proposition 5.3 it follows that the complementary series representation $\pi_{\lambda_0}$, with $\lambda_0 = \rho - \beta$, appears discretely in $L^2(G/C)$. By the connection between the class one spectrum for the $G$-representation, and the Laplacian spectrum for the corresponding locally symmetric space, Corollary 5.4 follows. $\qquad\square$

Let us postpone further the proof of Proposition 5.3, and continue discussing the proof of Theorem 1.6. Recall the following purely group theoretic notion: For a subgroup $C < G$, the *commensurator* of $C$ in $G$ is the subgroup $\mathrm{Comm}_C G = \{g \in G \mid gCg^{-1} \cap C \text{ has finite index in both } C \text{ and } gCg^{-1}\}$. Notice that even if $C$ is discrete, $\mathrm{Comm}_C G$ need not be a closed subgroup. We shall need the following general observation:

5.5 LEMMA. *Let $G$ be a rank one simple Lie group with finite center, and $C < G$ a torsion free discrete subgroup, for which the bottom $\tilde{\lambda}_0$ of the $L^2$-spectrum of the Laplacian on the locally symmetric space $C\backslash G/K$ is attained by a nonzero $L^2$-eigenfunction $f$. Then $C$ has finite co-volume in the closure of $\mathrm{Comm}_C G$ in $G$. In particular, if $\Gamma < G$ is a discrete subgroup containing $C$, then $C$ has finite index in $\mathrm{Comm}_C \Gamma$.*

*Proof.* It is a general fact about complete Riemannian manifolds, that if the bottom $\tilde{\lambda}_0$ of the spectrum of $\Delta$ occurs discretely, then up to a scalar multiple it admits a *unique* ($L^2$-)eigenfunction, which is in addition continuous and positive everywhere (up to a choice of sign); cf. [Su2, Theorem 2.8]. Consider the function $\tilde{f}$ on $G/K$ which lifts $f$ from $C\backslash G/K$, and fix $g \in \mathrm{Comm}_C G$. Since $\tilde{f}$ is $C$-invariant, the function $g\tilde{f}(x) = \tilde{f}(g^{-1}x)$ is $gCg^{-1}$-invariant, and in particular it is $C_0$-invariant, where $C_0 = gCg^{-1} \cap C$. By the finite index property of $C_0$, the projection to $C_0\backslash G/K$ of both $\tilde{f}$ and $g\tilde{f}$, is an $L^2$-eigenfunction of $\Delta$ with eigenvalue $\tilde{\lambda}_0$ (which is still the bottom of the spectrum). By the uniqueness property it follows that the projections of the two functions have constant ratio, denoted $\chi(g)$, where $\chi: \mathrm{Comm}_C G \to \mathbb{R}^*$ is a multiplicative character. For every $g \in \mathrm{Ker}\,\chi$ we then have $g\tilde{f} = \tilde{f}$, and as the elements in $G$ which fix $\tilde{f}$ form a closed subgroup, $\tilde{f}$ is fixed by $H = \overline{\mathrm{Ker}\,\chi}$ (of course, $C < H$). However, recalling that the projection of $\tilde{f}$ to $C\backslash G/K$ is square-integrable, we see easily that the $H$-invariance of $\tilde{f}$ is possible only if $C$ has finite co-volume in $H$.

Thus, it remains only to show that $\chi \equiv 1$. Indeed, every $g \in \mathrm{Comm}_G C$ normalizes $H$ (which, admitting the lattice $C$, is a unimodular subgroup), hence its left action on $H\backslash G$ is well-defined, and preserves the $G$-invariant measure on the latter. On the other hand, since $C$ is a lattice in $H$, and $f \in L^2(C\backslash G/K)$, the $H$-invariance of $\tilde{f}$ and the compactness of $K$ show that $\tilde{f}$ lifts (and then descends) to an $L^2$-function $F$ on $H\backslash G$. Integrating over $H\backslash G$, we get for all $g \in \mathrm{Comm}_G C$ the equality $\chi(g)^2 \int F^2(x)\ dx = \int F^2(gx)\ dx = \int F^2(x)\ dx$; hence $\chi(g) = 1$, as required. $\qquad\square$

We are now able to complete the proof of the second part of Theorem 1.6 (modulo Proposition 5.3). Replacing $\Gamma$ by a finite index subgroup, we may assume that it is torsion free. Recall that the proof of the first part of Theorem 1.6 (i.e., the weak inequality), showed that it is enough to study an (oriented) edge stabilizer $C < \Gamma$ such that $H^1(G, L^2(G/C)) \neq 0$, a fact which implies $\delta(C) \geq \delta(\Gamma) - 1$. We further need to rule out a possibility of equality in either one of the two cases in the second part of the theorem. Notice that without any additional assumption, by Corollary 5.4 it follows that if equality holds, then the bottom of the Laplacian spectrum on $C\backslash G/K$ is attained.

Consider the case of an isometric action on a locally finite tree $T$. It is easy to see that stabilizers of any two adjacent edges are commensurable, and hence, by connectivity, all edge stabilizers are commensurable one to the other, and the group $N < \Gamma$ they generate, commensurates each one of them. It follows from the above discussion and the last statement of 5.5, that $C$ has finite index in $N$. Since $N$ is normal in $\Gamma$, $C$ must be commensurated by $\Gamma$. Again by 5.5 it follows that $C$ has finite index in $\Gamma$, which is of course impossible under our assumption $\delta(C) = \delta(\Gamma) - 1$. This completes the proof in the case of actions on locally finite trees.

To establish the strict inequality for lattice $\Gamma < \mathrm{SU}(n,1)$ ($n \geq 2$), we call on a result of Gromov and Schoen [GS, §9], asserting that any action of $\Gamma$ on a tree, factors on some nonempty invariant subtree, through a discrete homomorphism to $\mathrm{PSL}_2(\mathbb{R})$. (That result assumes nonelementarity of the action, which is an interesting case for us as well, since otherwise it is easy to see that there are edge stabilizers with critical exponents tending to $\delta(\Gamma)$.) Since the image of such a homomorphism has positive first $\ell^2$-Betti number (e.g., by Theorem 1.5), the asserted strict inequality follows from the next general result, which may be of independent interest (and also sharpens Theorem 1.4):

5.6 THEOREM. *Let $\Gamma < G$ be a lattice, with $G$ not locally isomorphic to $\mathrm{SO}(k,1)$, $k = 2,3$. Let $N < \Gamma$ be the kernel of an epimorphism of $\Gamma$ onto a group with positive first $\ell^2$-Betti number. Then $\delta(N) > \delta(\Gamma) - 1$.*

*Proof.* By [BV] (see Theorem 4.3) the regular representation of the group $\Gamma/N$ has nonvanishing first cohomology. Lifting the representation and the 1-cocycle to $\Gamma$, and then inducing to $G$, yield $H^1(G, L^2(G/N)) \neq 0$ (by Theorem 1.10 and the isomorphism $\mathrm{Ind}_\Gamma^G \ell^2(\Gamma/N) \cong L^2(G/N)$). The foregoing discussion and the proof of the first part of Theorem 1.6 show first that this implies $\delta(N) \geq \delta(\Gamma) - 1$, and secondly, if equality holds, then $[\mathrm{Comm}_N(\Gamma) : N] < \infty$ (by Corollary 5.4 and Lemma 5.5). However $N$ is normal in $\Gamma$, so it must have finite index, which again contradicts the possibility $\delta(N) = \delta(\Gamma) - 1$. $\qquad\square$

Finally, we turn to the proof of Proposition 5.3, which will then complete the proof of Theorem 1.6, and with it our whole discussion of actions on trees.

*Proof of Proposition* 5.3. We shall assume, only to avoid introducing further notation, that $p = \delta(G)$ (which is the value needed in the proof of Theorem 1.6). The proof of Proposition 5.3 for a general value of $p$ with $2 < p < \infty$ is identical.

We shall need to return to the discussion and notation in the proof of Proposition 2.5, taking $G$ to be one of the simple Lie groups of interest to us, and $K < G$ a maximal compact subgroup. Let $u = v$ be a $K$-finite vector in the assumed irreducible subrepresentation $\sigma \subset L^2(M) = \tau$, and let $\varphi = \psi$ denote the $K$-invariant nonnegative function defined in the proof of Proposition 2.5 via $u$. It is obvious by the definition of $\varphi$, that we have for every $g \in G$, $|\langle \sigma(g)u, u\rangle| \leq \langle \tau(g)\varphi, \varphi\rangle$. While in Proposition 2.5 we have used this kind of inequality to obtain an upper bound on the left side in term of the right one, here we shall use it to get a lower bound on the right side in terms of the left one. More precisely, as discussed in the proof of step (ii) of Theorem 4.2, the Langlands classification and the fact that $\sigma$ is assumed to be strongly $L^p$ for $p = \delta(G) > 2$ (but not strongly $L^q$ for $q < \delta(G)$), imply that we may have chosen the $K$-finite vector $u$ to satisfy

$$(1) \qquad \lim_{a \to \infty,\, a \in A^+} e^{\beta(\log a)} |\langle \sigma(a)u, u\rangle| = C > 0.$$

Therefore, in view of the relation between the matrix coefficients of $u$ and $\varphi$, (1) implies

$$(2) \qquad \liminf_{a \to \infty,\, a \in A^+} e^{\beta(\log a)} \langle \tau(a)\varphi, \varphi\rangle \neq 0.$$

Note that while in (1) we have $C < \infty$, we do not claim here (although it is true) that the lim inf in (2) is finite as well.

Now, the $K$-invariance of $\varphi$ yields, by a standard direct integral computation, that for every $g \in G$ one has: $\langle \tau(g)\varphi, \varphi\rangle = \int \varphi_\lambda(g)\, d\mu(\lambda)$, where $\varphi_\lambda$ is the spherical function of the class one representation $\pi_{\lambda\beta}$, with $\lambda \in [0, \frac{1}{2}\delta(G)] \cup i\mathbb{R}$, and the positive Borel measure $\mu$ has total mass $\|\varphi\|^2 < \infty$ (see §2.III for the notation). Notice that by our assumption that $L^2(M)$ is strongly $L^{\delta(G)}$, Theorem 2.1, and the computations in Section 2.I of $L^p$-integrability of complementary series representations, it follows that $\mu$ is supported on $[0, \frac{1}{2}\delta(G) - 1] \cup i\mathbb{R}$.

The assertion of Proposition 5.3 would follow by showing that $\mu$ has an atom at $\frac{1}{2}\delta(G) - 1$, as the complementary series representation at that value of $\lambda$ is precisely the one claimed to embed in $L^2(M)$. Our strategy is to assume that $\mu(\{\frac{1}{2}\delta(G) - 1\}) = 0$, and derive a contradiction to the nonvanishing in (2), using the faster decay of the other $\varphi_\lambda$'s, and Lebesgue's dominated convergence theorem.

Indeed, recall that for $t \in \mathbb{R}$ we have $|\varphi_{it}| \leq \varphi_0$ (cf. Theorem 2.2), and for $0 \leq \lambda_1 \leq \lambda_2 \leq \frac{1}{2}\delta(G)$, we have $0 < \varphi_{\lambda_1} \leq \varphi_{\lambda_2} \leq 1$ (this well-known fact follows, e.g., from the integral representation of the spherical functions; cf. [GV, p.104 (3.1.15)]). Recall also that when $\lambda$ is *positive*, there exists a *finite* limit:

$$(3) \qquad \lim_{a \to \infty, \, a \in A^+} e^{(\rho - \lambda)(\log a)} \varphi_\lambda(a) = L_\lambda$$

(this well-known fact is actually covered by the discussion in step (ii) of Theorem 4.2; see also [GV, p. 172, Th. 4.7.4 (a)]). In particular, it follows from (3) for $\lambda = (\frac{1}{2}\delta(G) - 1)\beta = \rho - \beta$ and the monotonicity of the $\varphi_\lambda$'s over the positive $\lambda$'s, that for some $L < \infty$, all the $\varphi_\lambda$'s with $\mathrm{Re}\lambda \leq \frac{1}{2}\delta(G) - 1$ satisfy for $a \in A^+$ the uniform bound $e^{\beta(\log a)}|\varphi_\lambda(a)| \leq L$. Therefore, Lebesgue's convergence theorem may be applied to show that when $a \to \infty$ in $A^+$, one has:

$$(4) \qquad \int e^{\beta(\log a)} \varphi_\lambda(a) \, d\mu(\lambda) \to 0,$$

once we show that for $\mu$ almost every $\lambda$, the expression in the integrand tends to zero. This, however, follows from (3) and our assumption that $\mu(\{\frac{1}{2}\delta(G) - 1\}) = 0$, when $\lambda \in (0, \frac{1}{2}\delta(G) - 1]$, and from the bound on $a \in A^+$,

$$|\varphi_\lambda(a)| \leq \varphi_0(a) = \Xi(a) \leq C(1 + \beta(\log a)) \cdot e^{-\frac{1}{2}\delta(G)\beta(\log a)}$$

for all $\lambda \in i\mathbb{R}$ (see (9) in §2.III). As the expressions for which the limits in (4) and (2) are taken, identify, this establishes the required contradiction.  $\square$

## 6. Principal bundles, induction, and cohomology

This section is devoted to additional preliminaries toward the proofs of the remaining results, in Section 7. The reader may prefer to skip to that section on first reading, after looking at Definition 6.1 below, which describes the setting in which we shall work in the sequel.

I. *Principal bundles and induction of unitary representations.* The framework of principal bundles will be convenient, and we now define the exact notion to be used.

6.1 *Definition.* Suppose that $\tilde{M}$ is a locally compact space, equipped with a $\sigma$-finite Borel measure $\mu$, on which a discrete group $\Lambda$ acts properly discontinuously, preserving $\mu$. Denote the $\Lambda$-action from the right, and let $M = \tilde{M}/\Lambda$ be the quotient space. We refer to this structure as a *covering (principal)* $\Lambda$*-bundle* $\tilde{M} \to M$. We say that it has a compact/finite measure base, if $M$ has this property. If $\tilde{M}$ has no topology, and the $\Lambda$-action is measurably proper with a fundamental domain $M \cong \tilde{M}/\Lambda$, we call $\tilde{M} \to M$ a *measurable covering* $\Lambda$*-bundle.* If a locally compact group $G$ acts continuously (measurably) on $M$ and $\tilde{M}$, preserving $\mu$, so that its action on $\tilde{M}$ commutes with that of $\Lambda$ (and hence may safely be denoted from the left), we say that $G$ acts by (measurable) *covering* $\Lambda$*-bundle automorphisms.*

Choosing a measurable section $\varphi$ for a measurable covering $\Lambda$-bundle, namely, a fundamental domain $\varphi(M) = X \subset \tilde{M}$ for the $\Lambda$-action, yields a cocycle $\alpha$ for any action of a group $G$ by measurable covering $\Lambda$-bundle automorphisms. More precisely, $\alpha : G \times M \to \Lambda$ is defined by:

(1)    $\alpha(g, m) = \lambda$    if and only if    $g\varphi(m)\lambda \in X$.

The outstanding example is of course that of a discrete subgroup $\Lambda$ of a locally compact group $G$, where $\tilde{M} = G$ and $M = G/\Lambda$. Our notation then reduces to that in Section 3.III (compare (3) there with (1) above). Exactly as in Section 3.III, one may then induce a unitary representation $\pi$ from $\Lambda$ to $G$, keeping the same notation $\mathrm{Ind}_\Lambda^G \pi$ (or sometimes just $\mathrm{Ind}\pi$, where the context is clear). The representation space is $L^2(M, \mathcal{H})$, and the $G$-action is given as in (4) there. Although *a priori* this construction depends on the choice of $X$, it is easily verified that a different section (i.e., fundamental domain) yields a cohomologous cocycle, and hence an isomorphic representation. Some basic properties that we shall need, following standard properties of Mackey's unitary induction operation, are given in the following:

6.2 LEMMA. *Suppose that $G$ acts measurably on a measurable covering $\Lambda$-bundle $\tilde{M} \to M$ with finite measure base. Then for any unitary $\Lambda$-representations, $\pi, \sigma$, with $\mathrm{Ind}$ standing for $\mathrm{Ind}_\Lambda^G$, one has:*

(i) $\mathrm{Ind}(\pi \oplus \sigma) = \mathrm{Ind}\pi \oplus \mathrm{Ind}\sigma$.

(ii) *If $\pi$ is a mixing $\Lambda$-representation, and the natural $G$-representation on $L^2(\tilde{M})$ is mixing as well, then so is the $G$-representation $\mathrm{Ind}\pi$ (recall that a unitary representation is called mixing if all the associated matrix coefficients vanish at infinity).*

(iii) *If $\pi$ is weakly contained in $\sigma$ then $\mathrm{Ind}\pi$ is weakly contained in $\mathrm{Ind}\sigma$.*

(iv) $\mathrm{Ind}\ell^2(\Lambda) \cong L^2(\tilde{M}, \mu)$, *where the latter is the natural $G$-representation coming from its measure-preserving action on $\tilde{M}$, and $\ell^2(\Lambda)$ is the left regular $\Lambda$-representation.*

*Proof.* (i) is trivial. To prove (ii), recall that it is enough to exhibit a dense subspace of Ind$\pi$ for which all diagonal matrix coefficients vanish at infinity. The subspace considered is that of the functions taking finitely many values (in $\mathcal{H}_\pi$). Let $f \in L^2(M, \mathcal{H}_\pi)$ be such a function, $\{v_i\}$ the finite set of vectors in its image, and $C < \infty$ a bound on the norms of all the $v_i$'s. Fix some $\epsilon > 0$. Let $F \subset \Lambda$ be the set of $\lambda$'s for which $|\langle \pi(\lambda) v_i, v_j \rangle| > \epsilon$ for some $i, j$. By the assumption on $\pi$, $F$ is finite. Let $X \subset \tilde{M}$ and $\alpha \colon G \times M \to \Lambda$ be as in (1) above. Denote $Y = XF^{-1} = \cup X\lambda^{-1}$, $\lambda \in F$. Notice that $Y \subset \tilde{M}$ has finite measure. Now, given any $g \in G$, dividing $X$ into the subsets where $\alpha(g^{-1}, x)$ belongs, or does not belong, to $F$, we can estimate:

$$|\langle gf, f \rangle| = |\int_X \langle \pi(\alpha(g^{-1}, x)) f(g^{-1}x), f(x) \rangle d\mu(x)| \leq \epsilon\mu(X) + C^2\mu(X \cap gY).$$

However, applying the assumption that $L^2(\tilde{M})$ is mixing to the characteristic functions of $X$ and $Y$, yields that $\mu(X \cap gY) \to 0$ when $g \to \infty$. Since $\epsilon$ was arbitrary, we are done.

Claim (iii) is again standard; see the proof of [Zi1, 7.3.7]. Assertion (iv) and its proof are analogous to the well-known statement $\text{Ind}_\Lambda^G \ell^2(\Lambda) \cong L^2(G)$ in the group theoretic setting. For the proof we distinguish between the $G$-action on $M$ and $\tilde{M}$, by writing for the former $g \cdot m$. Notice that by the definition (1) of $\alpha$, and the section $\varphi$, we have:

$$(2) \qquad\qquad \varphi(g \cdot m) = g\varphi(m)\alpha(g, m).$$

The space $L^2(M, \ell^2(\Lambda))$ may naturally be identified with $L^2(M \times \Lambda)$ with the product measure, by assigning $f(m)(\lambda) \to f(m, \lambda)$. (Note that $\int_{M \times \Lambda} |f(m, \lambda)|^2 = \int_M \int_\Lambda |f(m)(\lambda)|^2 = \int_M \|f(m)\|_{\ell^2(\Lambda)}^2$.) The $G$-action on the former space takes the form:

$$[gf(m)](\lambda) = [\alpha(g^{-1}, m) f(g^{-1} \cdot m)](\lambda) = [f(g^{-1} \cdot m)](\alpha(g^{-1}, m)^{-1}\lambda),$$

and thus its action on the latter is:

$$(3) \qquad\qquad gf(m, \lambda) = f(g^{-1} \cdot m, \alpha(g^{-1}, m)^{-1}\lambda).$$

Next, $M \times \Lambda$ may naturally be identified with $X \times \Lambda$ by $(m, \lambda) \to (\varphi(m), \lambda)$. This identifies also the spaces $L^2(M \times \Lambda)$ and $L^2(X \times \Lambda)$, and (3) now has the form:

$$\begin{aligned} gf(x, \lambda) &= f(\varphi(g^{-1} \cdot \varphi^{-1}(x)), \alpha(g^{-1}, \varphi^{-1}(x))^{-1}\lambda) \\ &= f(g^{-1}x\alpha(g^{-1}, \varphi^{-1}(x)), \alpha(g^{-1}, \varphi^{-1}(x))^{-1}\lambda) \end{aligned}$$

(the last equality follows by substituting $g = g^{-1}$ and $m = \varphi^{-1}(x)$ in (2)).

It is now clear that the map $(x, \lambda) \to x\lambda$ from $X \times \Lambda$ to $\tilde{M}$ is a measure-preserving isomorphism, which takes the latter $G$-action to the natural $G$-action on $L^2(\tilde{M})$: $gf(m) = f(g^{-1}m)$, as required. $\qquad \square$

II. *A "transference" theorem.* The following technical result will be essential in the proofs of Theorems 1.2, 1.4 and 1.7.

6.3 THEOREM. *Let $\tilde{M} \to M$ be a measurable covering $\Lambda$-bundle with finite measure base, acted upon by $G = \mathrm{SO}(n,1)$ or $\mathrm{SU}(n,1)$, by bundle automorphisms. Assume that the $G$-representation on $L^2(\tilde{M})$ is strongly $L^p$ for some $2 \le p < \infty$. Suppose that the group $\Lambda$ is torsion free, and is discretely embedded in $H = \mathrm{SO}(m,1)$ or $\mathrm{SU}(m,1)$ for some $m \in \mathbb{N}$. Identify $\Lambda$ with its image in $H$. Let $\pi$ be a unitary $H$-representation which is strongly $L^q$ for some $2 \le q < \infty$ (as an $H$-representation). Then $\mathrm{Ind}_\Lambda^G \pi|_\Lambda$ is strongly $L^r$ for $r = \max\{\frac{\delta(\Lambda)}{\delta(H)} \cdot \frac{pq}{2}, p\}$.*

*Proof.* The structure of the proof is as follows: In the first part we show how to convert the information given by the assumption into an $L^p$-integrability property of the associated cocycle $\alpha$; see (8) below. It is only in the second part that the representation $\pi$ is considered, where we show how to deduce the required estimate from (8).

Let $\rho = \rho(H)$ be the half sum of the positive roots of $H$, namely, $2\rho = \delta(H)\beta$ (we retain the notation of Sections 2 and 3, and in particular that for the complementary series representations $\pi_\lambda$, discussed in 2.II). Recall that for $\lambda = |\lambda|\beta \in \mathfrak{a}^*$, $0 \le |\lambda| \le |\rho| = \frac{1}{2}\delta(H)$, the representation $\pi_\lambda$ is strongly $L^t$ for $t = 2\delta(H)/(\delta(H) - 2|\lambda|)$ (see (∗) in 2.III). By Corollary 2.9, $\pi_\lambda|_\Lambda$ is strongly $L^t$ for $t = 2\delta(\Lambda)/(\delta(H) - 2|\lambda|)$. Assume for the moment that $\delta(\Lambda) > 0$. Let $\lambda_0$ be defined by the equality $t = 2$, or more precisely:

$$(5) \qquad 2\delta(\Lambda)/(\delta(H) - 2|\lambda_0|) = 2, \quad \text{equivalently,} \quad |\lambda_0| = \frac{1}{2}(\delta(H) - \delta(\Lambda)).$$

If $\delta(\Lambda) = 0$ replace $\delta(\Lambda)$ by any positive value strictly less than 2, say, by 1. Notice that in both cases $\lambda_0 < \rho$ and $\pi_{\lambda_0}$ is nontrivial. Therefore, by the choice of $\lambda_0$, $\pi_{\lambda_0}|_\Lambda$ is strongly $L^2$ (when $\delta(\Lambda) = 0$ we get that it is strongly $L^1$, hence also $L^2$). By [CHH] it is weakly contained in $\ell^2(\Lambda)$. From Lemma 6.2 (iii), (iv), and the assumption, it now follows that $\mathrm{Ind}_\Lambda^G \pi_{\lambda_0}|_\Lambda$ is weakly contained in a strongly $L^p$ representation, hence by Theorem 2.1 (4), it is itself strongly $L^p$.

We claim that by the torsion freeness assumption, one can always find a $K$-equivariant section ($K < G$ maximal compact) $\varphi : M \to \tilde{M}$ for the natural projection, or equivalently, a fundamental domain $X \subset \tilde{M}$ for the (right) $\Lambda$-action, which is (left) $K$-invariant. Let us postpone the proof of this fact to the end of the proof (see 6.4 below), and assume for the moment that $X \subset \tilde{M}$ has this property. Let $\alpha : G \times M \to \Lambda$ be the cocycle as in (1), corresponding to this $X$. Let $v_0$ be the (unique up to scalar) $K$-invariant unit vector in $\pi_{\lambda_0}$, and consider the constant function $f : X \to \mathcal{H}_{\pi_{\lambda_0}}$: $f(x) \equiv v_0$, in the representation space $L^2(X, \mathcal{H}_{\pi_{\lambda_0}})$ of $\mathrm{Ind}_\Lambda^G \pi_{\lambda_0}|_\Lambda$. Then $f$ is $K$-invariant for the $G$-action (this follows only from $K$-invariance of $X$), so by Theorem 2.1 and the fact (established above) that this representation is strongly $L^p$, we get for every $\varepsilon > 0$:

$$(6) \quad \int_G |\int_X \langle gf, f \rangle dx|^{p+\varepsilon} dg = \int_G |\int_X \langle \pi_{\lambda_0} \alpha(g^{-1}, x) v_0, v_0 \rangle dx|^{p+\varepsilon} dg < \infty.$$

However, recall from (11) in Section 2 (and the remark thereafter), that for some constant $C > 0$:

$$(7) \quad 0 < C \cdot e^{(\lambda_0 - \rho) \log a(\alpha(g^{-1}, x))} \leq \langle \pi_{\lambda_0} \alpha(g^{-1}, x) v_0, v_0 \rangle.$$

Thus, substituting (7) in (6) we deduce that for all $\varepsilon > 0$:

$$\int_G (\int_X e^{(\lambda_0 - \rho) \log a(\alpha(g^{-1}, x))} dx)^{p+\varepsilon} dg < \infty$$

or equivalently, by (5):

$$(8) \quad \int_G (\int_X e^{-\frac{1}{2} \delta(\Lambda) \beta \log a(\alpha(g^{-1}, x))} dx)^{p+\varepsilon} dg < \infty.$$

Moving on to the second part of the proof, let $\pi$ be some unitary $H$-representation which is strongly $L^q$. If $q\delta(\Lambda)/\delta(H) \leq 2$ then from Corollary 2.9, $\pi|_\lambda$ is strongly $L^2$. By the argument above, $\mathrm{Ind}_\Lambda^G \pi|_\Lambda$ is strongly $L^p$, so we have nothing to prove. Otherwise we have $r = \frac{\delta(\Lambda)}{\delta(H)} \cdot \frac{pq}{2} > p$. Let $\mathcal{H}_0 \subset \mathcal{H}$ be the (dense) subspace of $K$-finite vectors. Then, as in the proof of 4.10, the subspace of functions $f : X \to \mathcal{H}_0$ which take finitely many values is dense, and it suffices to prove that the matrix coefficient associated with any two such functions is in $L^{r+\varepsilon}(G)$, for all $\varepsilon > 0$ (where $r$ is as in the theorem).

Let then $\varphi, \psi : X \to \mathcal{H}_0$ be two such functions, and $F = \{u_i\} \subset \mathcal{H}_0$ be the finite set of ($K$-finite) vectors in their image. Then, for every $\theta > 0$ we can find, by Theorem 2.1, a constant $C_\theta < \infty$, such that for every $h \in H$:

$$(9) \quad \sum_{i,j} |\langle \pi(h) u_i, u_j \rangle| < C_\theta \cdot e^{(-\frac{\delta(H)}{q} + \theta) \beta \log a(h)}.$$

Fix some "small" $\theta$ for the moment. Choose any $r_0 > r$, and take $p_0 > p$ such that $t = \frac{r_0}{p_0} > \frac{r}{p}$. As we assume now that $r > p$, we have $t > 1$. By Hölder inequality and the assumption that $X$ has finite total measure (which may assumed to be one), we have for every positive function $f : X \to \mathbb{R}$, the inequality $(\int f\,dx)^{r_0} \le (\int f^t dx)^{p_0}$. Applying this and (9) yields (arguing as in the proof of Lemma 4.10):

$$\int_G |\langle g\varphi, \psi\rangle|^{r_0} dg \le \int_G (\int_X C_\theta e^{(-\frac{\delta(H)}{q}+\theta)\beta \log a(\alpha(g^{-1},x))} dx)^{r_0} dg$$

$$\le C_\theta^{r_0} \int_G (\int_X e^{(-\frac{t\delta(H)}{q}+t\theta)\beta \log a(\alpha(g^{-1},x))} dx)^{p_0} dg.$$

In light of (8), the last integral converges once $-t\frac{\delta(H)}{q}+t\theta < -\frac{1}{2}\delta(\Lambda)$, or $\theta < \frac{\delta(H)}{q} - \frac{1}{2}\frac{\delta(\Lambda)}{t}$. The right-hand side in the latter inequality would be zero for $t = r/p$, by definition of (and our current assumption on) $r$. Since $t > r/p$, the integral would indeed converge had $\theta > 0$ been taken small enough (according to the choice of $r_0$ and $p_0$, which has to be made first). Since $r_0 > r$ is arbitrary, the theorem follows. We have only to complete the following unfinished point in the proof:

6.4 LEMMA. *With the notation of Theorem 6.3, one can always find a measurable $K$-equivariant section $\varphi : M \to \tilde{M}$.*

The result is standard, but we failed to find a reference. For completeness we sketch the proof (note that the torsion freeness assumption is necessary).

*Proof.* Let $\psi : M \to \tilde{M}$ be a section and $\alpha : K \times M \to \tilde{M}$ the associate cocycle. Consider the induced representation $\operatorname{Ind}_\Lambda^K \ell^2(\Lambda)$, which as in the proof of (iv) in Lemma 6.2, is isomorphic to the $K$-representation on $L^2(M \times \Lambda)$, where $K$ acts on $M \times \Lambda$ by $k(m, \lambda) = (km, \alpha(k, m)^{-1}\lambda)$. Since $K$ is compact, such a representation must have a nonzero invariant vector. Suppose for the moment that $K$ acts ergodically on $M$. Then, as stabilizers of nonzero vectors in $\ell^2(\Lambda)$ (for the $\Lambda$-regular representation) must be trivial (by torsion freeness), we get from Zimmer's cocycle reduction lemma [Zi1, 5.2.11], that $\alpha$ is cohomologous to a cocycle taking values in $e \in \Lambda$. Equivalently, for some $f : M \to \Lambda$ we have $f(m)^{-1}\alpha(k, m)f(km) = e$ for all $k \in K$ and (almost all) $m$. The section $\varphi(m) = \psi(m)f(m)$ then satisfies $\varphi(km) = k\varphi(m)$ as required. In general, of course, the $K$-action on $M$ is not ergodic, but one can use the above argument for each ergodic component (which is just an orbit). The proof of Lemma 6.4, and hence that of Theorem 6.3, is now complete. □

III. *The bundle-induction operation on the first cohomology.* Our goal here is to show that under suitable conditions, one has an injection of $H^1(\Lambda, \pi)$ into $H^1(G, \operatorname{Ind}_\Lambda^G \pi)$. As in Section 3.III, we first need to know that an appropriate map between the cohomology groups exists, and then study its properties. The question of existence of the map is completely analogous to that in Section 3.III. The notation and arguments used there apply here verbatim, and in particular condition (10) there ensures that the map $b \to \tilde{b}$ defined in (5) yields a well-defined linear map $H^1(\Lambda, \pi) \to H^1(G, \operatorname{Ind}_\Lambda^G \pi)$. To summarize these facts, we have:

6.5 PROPOSITION. *Let $\tilde{M} \to M$ be a measurable covering $\Lambda$-bundle, acted upon by $G$ by bundle automorphisms. If there exists a fundamental domain $X$ for which the corresponding cocycle $\alpha : G \times X \to \Lambda$ satisfies condition (10) in Section 3, then the map $b \to \tilde{b}$ from $H^1(\Lambda, \pi)$ to $H^1(G, \operatorname{Ind}_\Lambda^G \pi)$ in (5) is well-defined. In particular, this is the case if $\alpha(g, \cdot)$ takes finitely many values for every $g \in G$, and $X$ (equivalently $M$) has finite measure. Consequently, if $\tilde{M} \to M$ is a (continuous) covering $\Lambda$-bundle with compact base, this condition is satisfied.*

Once the map between the first cohomology groups is well-defined, we may inquire as to when it is injective. This will not always be the case, and our next purpose is to describe a situation where it is. To this end we have the following:

6.6 LEMMA. *Let $G$ be a locally compact group, $(\pi, \mathcal{H})$ a unitary $G$-representation, and $b \in Z^1(G, \pi)$ a cocycle. Then the following conditions are equivalent:*

(i) *If $g \to \infty$ in $G$ then $\|b(g)\| \to \infty$ in $\mathcal{H}$.*

(ii) *The associated isometric $G$-action $\rho$ on $\mathcal{H}$ is proper (see Section 3.I).*

(iii) *There is one vector $v \in \mathcal{H}$ such that $\|\rho(g)v\| \to \infty$ whenever $g \to \infty$.*

*When one (and hence all) of these conditions hold, we call the cocycle $b$ proper. If $G$ is not compact, and $b$ is proper, then $b$ is not a coboundary.*

The verification of the lemma is simple and we omit it. The following result is a little more general than the one we shall actually need here, but it will also be used elsewhere.

6.7 LEMMA. *Let $\tilde{M} \to M$ be a measurable $\Lambda$-bundle with a base of finite measure (denoted $\mu$), on which a locally compact group $G$ acts by bundle automorphisms. Suppose that for the associated cocycle $\alpha : G \times M \to \Lambda$, the following two conditions are satisfied:*

(i) *$\alpha(g, \cdot)$ takes finitely many values in $\Lambda$, for any $g \in G$.*

(ii) *For every sequence $g_i \to \infty$ and $\lambda \in \Lambda$, $\mu\{m \in M | \alpha(g_i, m) = \lambda\} \to 0$.*

*Then for every unitary $\Lambda$-representation $\pi$, condition (10) in Section 3 is satisfied, and for every proper cocycle $b \in Z^1(\Lambda, \pi)$, the cocycle $\tilde{b} \in Z^1(G, \operatorname{Ind}_\Lambda^G \pi)$ is proper as well.*

*Furthermore, conditions* (i) *and* (ii) *are satisfied if* $\tilde{M} \to M$ *is a covering* $\Lambda$*-bundle with compact base* (*Def.* 6.1), *where the cocycle* $\alpha$ *corresponds to a bounded section, and the natural* $G$*-representation on* $L^2(\tilde{M})$ *is mixing.*

*Proof.* For the first assertion see Proposition 6.5. To show the second, we switch to the framework of isometric actions and use condition (iii) in Lemma 6.6. Denote by $\rho$ the isometric action of $\Lambda$ on $\mathcal{H}$ which corresponds to $b$, and let $v \in \mathcal{H}$ be any vector. Consider the constant function $f(m) \equiv v$ in $L^2(M, \mathcal{H})$. By the assumption on $b$, we may find for every $C < \infty$ a finite set $F \subset \Lambda$, such that for every $\lambda$ outside $F$ one has $\|\rho(\lambda)v\| > C$. Given a sequence $g_i \to \infty$, by assumption (ii) we have for all $i$ large enough

$$\mu\{m \in M | \alpha(g_i^{-1}, m) \in F\} < \frac{1}{2}\mu(M),$$

which gives

$$\int\limits_M \|\rho \circ \alpha(g_i^{-1}, m) f(g_i^{-1}m)\|^2 d\mu(m) > C^2 \cdot \frac{1}{2}\mu(M).$$

Finally, in the assumed topological covering bundle situation, it is clear that condition (i) is satisfied. To establish (ii) denote by $X \subset \tilde{M}$ the image of a bounded section of the bundle, and fix $\lambda \in \Lambda$. Then the mixing assumption implies that $\mu(X \cap g_i(X\lambda^{-1})) \to 0$, which is exactly the condition asserted in (ii). $\qquad\square$

## 7. Proofs of Theorems 1.2, 1.4 and 1.7

I. *Proof of Theorem* 1.4.   Retain the notation in the statement of the theorem. Denote $N = \text{Ker}\varphi$, $\Lambda = \text{Im}\varphi$, $\Lambda \cong \Gamma/N$. Since $\Lambda$ is linear, by passing to a subgroup of finite index in $\Lambda$ (and then in $\Gamma$), we may and shall assume that it is torsion free. Consider the natural projection (covering) map $\tilde{M} = G/N \to M = G/\Gamma$, which is naturally a $\Lambda$-bundle (Definition 6.1), using the (well-defined) $\Lambda$-action on $G/N$ from the right. Of course, $G$ acts from the left by $\Lambda$-bundle automorphisms. Suppose now that $\pi$ is some unitary $\Lambda$-representation. Then there are two ways to induce $\pi$ to a $G$-representation. The first is as in the procedure described in Section 6.I, with $G$ acting by $\Lambda$-bundle automorphisms. The second is by consideration of $\pi$ as a $\Gamma$-representation factoring through its projection to $\Lambda$, and then by induction as in Section 3.III, according to the usual Mackey's construction. It is an easy exercise, which we leave for the reader to verify, that these two ways yield (naturally) isomorphic representations. We shall need in the sequel both constructions: the first to control $L^p$-integrability of the induced representation (making use of Theorem 6.3), and the second to deduce that the first cohomology injects, if $\Gamma$ is nonuniform as well, applying Theorem 1.10.

First, notice that the statement of the theorem is trivial if $G = \mathrm{SO}(3,1)$, for the first summand is at least $\frac{1}{2}\delta(\Gamma)$ and the second is at least 1. If $\delta(\Gamma) \le 2$ we clearly have the inequality asserted there. Thus, we may exclude that group, and deduce from the foregoing discussion and Theorem 1.10, that for any lattice $\Gamma < G$, uniform or not, we have for every unitary $\Lambda$-representation $\sigma$, an isomorphism $H^1(\Lambda, \sigma) \cong H^1(G, \mathrm{Ind}_\Lambda^G \sigma)$.

The case $\frac{\delta(\Lambda)}{2} < 1$, namely, $\delta(\Lambda) < 2$, is dealt with by Theorem 5.6 and Theorem 1.5, which in fact yield a strict inequality. We may therefore assume that $\delta(\Lambda) \ge 2$. Let $\pi$ be an irreducible unitary $H$-representation with $H^1(H, \pi) \ne 0$. By (the proof of) Theorem 4.2, $\pi$ is strongly $L^q$ for $q = \delta(H)$. Recall that we have $\tilde{M} \cong G/N$. By Theorem 2.10 the $G$-representation on $L^2(\tilde{M}) \cong L^2(G/N)$ is strongly $L^p$ for $p = p_0 = \max\{2, \delta(G)/\delta(G) - \delta(N)\}$. Thus, we may substitute this $p = p_0$, and the above value of $q$, in Theorem 6.3. It follows that $\mathrm{Ind}_\Lambda^G \pi|_\Lambda$ is strongly $L^r$ for $r = \max\{\frac{\delta(\Lambda)}{\delta(H)} \cdot \frac{p_0 q}{2}, p_0\} = \delta(\Lambda)\frac{p_0}{2}$ (as we are assuming $\delta(\Lambda) \ge 2$). On the other hand, by Theorem 3.4 we have $H^1(\Lambda, \pi|_\Lambda) \ne 0$, and hence, by the discussion at the beginning of the proof, $H^1(G, \mathrm{Ind}_\Lambda^G \pi|_\Lambda) \ne 0$. Theorem 4.2 now implies that $\mathrm{Ind}_\Lambda^G \pi|_\Lambda$ is strongly $L^r$ only if $r \ge \delta(G) = \delta(\Gamma)$. Comparing this with the above value of $r$, we see that we must have $\delta(\Lambda)\frac{p_0}{2} \ge \delta(\Gamma)$. Substituting for $p_0$ we deduce:

$$(1) \qquad \delta(\Gamma) \le \frac{\delta(\Lambda)}{2}\max\{2, \delta(\Gamma)/\delta(\Gamma) - \delta(N)\}.$$

Now, if $\delta(\Gamma)/\delta(\Gamma) - \delta(N) \le 2$ then (1) implies that $\delta(\Gamma) \le \delta(\Lambda)$, which is just what the theorem claims. Otherwise, we conclude:

$$\delta(\Gamma) \le \frac{1}{2}\delta(\Lambda)\delta(\Gamma)/\delta(\Gamma) - \delta(N) \qquad \text{and hence} \qquad \delta(\Gamma) - \delta(N) \le \frac{1}{2}\delta(\Lambda),$$

which completes the proof of the theorem. $\qquad\square$

Finally, let us discuss the observation made in the introduction, following Theorem 1.4. The starting point of this construction is the existence of homomorphisms of certain lattices in $\mathrm{SO}(n,1)$ (for all $n \ge 2$), onto nonabelian free groups; see e.g. [Lub2]. Suppose that $\varphi : \Gamma \to F$ is such a homomorphism, and denote by $N$ its kernel. Notice the following general claim: If $\Gamma/N$ is nonamenable (which is the case for us), then $\delta(N) < \delta(\Gamma)$ (the "only if" part holds as well, and is easier to establish). To show this claim observe that the conclusion is equivalent to the statement that the bottom of the Laplacian spectrum on $N\backslash\mathbb{H}^n$ is positive, which is itself equivalent to saying that the $G$-representation on $L^2(G/N)$ does not contain $1_G$ weakly (see the discussion in Section 2.V above). However by the assumption, $\ell^2(\Gamma/N)$ does not contain $1_\Gamma$ weakly, and we also have $\mathrm{Ind}_\Gamma^G \ell^2(\Gamma/N) \cong L^2(G/N)$. The claim now follows from [Mar2, III.1.11 (b)] and the well-known fact that $L_0^2(G/\Gamma)$ does not contain $1_G$ weakly. Now, returning to the observation, denote $\alpha = \delta(N) < \delta(\Gamma)$. Then, given any $\epsilon > 0$, one can embed the free group $F$ inside $\Gamma$, with a critical exponent less than $\epsilon$ (Lemma 2.7 (3)). Taking the composition as a homomorphism of $\Gamma$ into itself establishes now the observation. $\qquad\square$

II. *Proof of Theorem* 1.2. We consider first the case $G \cong \mathrm{SO}(n, 1)$ $n \geq 3$, and then discuss the modifications required for $\mathrm{SU}(n, 1)$.

Assume $G \cong \mathrm{SO}(n, 1)$, $n \geq 3$. Recall that the universal cover of $G$, $\tilde{G} \cong \mathrm{Spin}(n, 1)$, is a two-fold covering, and hence has a finite center. We may thus replace $G$ by $\tilde{G}$ and assume that $G$ is simply connected. Denote $\Lambda = \pi_1(M)$. Then by standard covering theory we may realize $\Lambda$ as acting discontinuously on $\tilde{M}$ (from the right), commuting with the lifted $G$-action (from the left). In other words, $\tilde{M} \to M$ is a covering $\Lambda$-bundle with compact base, on which $G$ acts by bundle automorphisms (see Definition 6.1).

The assumption of measurable properness of the action implies that the $G$-representation on $L^2(\tilde{M})$ may be embedded in a multiple of the regular $G$-representation. Indeed, by Mackey's well-known orbit theorem, $L^2(\tilde{M}) \cong \int L^2(G/G_x) d\mu(x)$, where the integration is taken over the (measurable by our assumption) set of orbits, and $G_x$ is a stabilizer of a point $x$ in the section (the measure $\mu$ is the one induced by $\tilde{M}$ on this space). However, by our assumption, for almost every $x$ the subgroup $G_x$ is compact, hence $L^2(G/G_x) \subseteq L^2(G)$. It follows that $L^2(\tilde{M})$ may be embedded in $\infty \cdot L^2(G)$, and in particular it is strongly $L^2$ (and mixing).

Let $\Lambda$ be identified with its discrete image in some $H \cong \mathrm{SO}(m, 1)$ or $\mathrm{SU}(m, 1)$. Let $\pi$ be an irreducible unitary $H$-representation with $H^1(H, \pi) \neq 0$. By (the proof of) Theorem 4.2, $\pi$ is strongly $L^{\delta(H)}$ (as an $H$-representation). We now wish to invoke Theorem 6.3 for the representation $\pi$. With notation as in that theorem, we may take by the above discussion $p = 2$, and $q = \delta(H)$. It then follows from Theorem 6.3 that the induced representation $\mathrm{Ind}_\Lambda^G \pi$ is strongly $L^r$ for $r = \max\{\delta(\Lambda), 2\}$.

Let us assume first that $n \geq 4$; that is, $G \neq \mathrm{SO}(3, 1) \cong \mathrm{SL}_2(\mathbb{C})$. Then $\delta(G) = n - 1 \geq 3$. Thus, if $\delta(\Lambda) < \delta(G)$ we have $r < \delta(G)$. On the other hand, by Proposition 6.5 the map $b \to \tilde{b}$ from $Z^1(\Lambda, \pi|_\Lambda)$ to $Z^1(G, \mathrm{Ind}_\Lambda^G \pi|_\Lambda)$ is well-defined. Choosing $b$ as the restriction to $\Lambda$ of a non-coboundary $H$-cocycle, we have by Theorem 3.4 and Lemma 6.6, that $b$ is proper. Hence, from Lemma 6.7 it follows that $H^1(G, \mathrm{Ind}_\Lambda^G \pi|_\Lambda) \neq 0$ (by the last part of Lemma 6.7, conditions (i) and (ii) there are indeed satisfied in our case). However from Theorem 4.2 we now deduce that $\mathrm{Ind}_\Lambda^G \pi|_\Lambda$ cannot be strongly $L^r$ for $r < \delta(G)$, a contradiction which proves the result.

We are left with the case $G = \mathrm{SO}(3, 1)$. Here we have to show that $\delta(\Lambda) \geq 2$. Otherwise it follows from Theorem 1.5 that $H^1(\Lambda, \ell^2(\Lambda)) \neq 0$, and that moreover, there is a proper $\Lambda$-cocycle for the regular $\Lambda$-representation (the latter follows from the proof of Theorem 1.5, going back to Theorem 3.4). As before we deduce that $H^1(G, \operatorname{Ind}_\Lambda^G \ell^2(\Lambda)) \neq 0$; hence by Lemma 6.2 (iv) $H^1(G, L^2(\tilde{M})) \neq 0$. However as we observed at the beginning of the proof, $L^2(\tilde{M}) \subseteq \infty \cdot L^2(G)$. If the latter representation has nonvanishing first cohomology, then so does $L^2(G)$; hence by Theorem 4.3 we would have $p(G) = 0$, contradicting Theorem 4.2.

*Remark.* Notice that in the case of $\operatorname{SO}(3,1) \cong \operatorname{SL}_2(\mathbb{C})$, the argument above, together with that used in the proof of Theorem 1.5, show that if equality $\delta(\pi_1(M)) = 2$ holds, then $\pi_1(M)$ must be a divergence group.

Let us discuss now the case $G = \operatorname{SU}(n,1)$. In the higher rank results of Zimmer, an assumption of finite fundamental group is always made on $G$, which is a minor restriction in view of the relatively few examples of simple groups without this property (see [Zi2]). As we cannot afford to lose "half of our clients", we discuss the issue of how these groups may be treated as well, which is not simple enough to be left untouched.

To illustrate the necessity of a modification of Theorem 1.2 for these groups, consider the outstanding example $M = G/\Gamma$, which motivates the general theorem. Here $\tilde{M} \cong \tilde{G} = \widetilde{\operatorname{SU}(n,1)}$, and $\pi_1(M) \cong \tilde{\Gamma}$ is just the "lift" of $\Gamma$ to $\tilde{G}$, which is an infinite central extension of $\Gamma$. Having an infinite center, $\tilde{\Gamma}$ admits no discrete embedding in a rank one simple Lie group with finite center. Consequently, our result regarding injective homomorphisms of $\pi_1(M)$ into such groups, although strictly speaking is valid, becomes of little interest. Analyzing the problem in the more general perspective of Theorem 1.2, observe that lifting the $G$-action on $M$, to an action of $\tilde{G}$ on $\tilde{M}$, one easily verifies by the construction of the lift that $p^{-1}(e)$, where $p : \tilde{G} \to G$ is the covering homomorphism, acts on $\tilde{M}$ by deck transformations over $M$, and hence maps naturally into $\pi_1(M)$ (possibly with kernel). Since the $\tilde{G}$-action on $\tilde{M}$ commutes with that of $\pi_1(M)$, the image of $p^{-1}(e)$ lies in the center of $\pi_1(M)$. As in the preceding discussion, of "real" interest to us are those homomorphisms of $\pi_1(M)$ which may have kernel contained in the image of $p^{-1}(e)$. We shall call such a homomorphism an *isogeny*, and prove Theorem 1.2 for this class of homomorphisms. In the exemplary case of $M = G/\Gamma$, $\pi_1(M) = \tilde{\Gamma}$, we see that any embedding of $\Gamma$, when viewed as coming from $\tilde{\Gamma}$ through the natural projection, defines an isogeny of $\pi_1(M)$, and will therefore be covered by our analysis.

Yet, when $G = \operatorname{SU}(n,1)$, the assumption that the lifted $\tilde{G}$ action on $\tilde{M}$ is proper (topologically or measurably) is too restrictive, as often the $G$-action lifts to an action of a finite cover of $G$, so the $\tilde{G}$-action will admit an infinite kernel. In addition, Gromov's theorem, which is so essential in our main application of Theorem 1.2, is not quite the same when $G$ has an infinite fundamental group. Rather than stating some technical general assumptions on the $G$-action, which resolve the above problems, and then verify them in our main application (Corollary 1.3), let us aim directly at this application, relying upon Gromov's theorem [Gr2]. The proof of the Corollary lends itself to the appropriate generalization, and we leave this matter to the reader.

7.1 THEOREM (Corollary 1.3 for SU($n,1$)). *Suppose that the connected Lie group $G$ has finite center, and is locally isomorphic to* SU($n,1$) ($n \geq 2$). *Assume that $G$ acts on a compact manifold $M$, satisfying all the conditions in Corollary 1.3. Let $\varphi : \pi_1(M) \to \Lambda$ be an isogeny (see the above discussion) onto an infinite discrete subgroup $\Lambda$ of $H \cong$ SO($m,1$) or SU($m,1$), for some $m \in \mathbb{N}$. Then $\delta(\Lambda) \geq \delta(G) = 2n$.*

*Proof.* Consider the lifted $\tilde{G}$ action on $\tilde{M}$. Denote by $p : \tilde{G} \to G$ the covering map. The proof proceeds now in the following steps:

1. *We may assume each element in $p^{-1}(e)$ acts properly on $\tilde{M}$, hence the $Z(\tilde{G})$-action is proper.* Indeed, by standard covering theory, every $z \in p^{-1}(e)$ acts as a deck transformation over $M$. Denoting by $\mathcal{N} \subseteq p^{-1}(e)$ the kernel of the action, the $\tilde{G}$-action on $\tilde{M}$ factors through an action of $\tilde{G}/\mathcal{N}$ (the latter still covers $G$, and we keep the notation $p$ for the covering homomorphism). Replacing $\tilde{G}$ by $\tilde{G}/\mathcal{N}$ throughout the rest of the proof, we reduce to the situation of 1. The last claim of 1 follows readily from the fact that $G$ has finite center, so that $p^{-1}(e)$ has finite index in $Z(\tilde{G})$.

2. *The $\tilde{G}$ action on $\tilde{M}$ is (measurably) proper.* Indeed, notice that Gromov's theorem in case $\pi_1(G)$ is infinite ([Gr2, 6.1.B]; see also Remark 6.1.c therein), asserts that the lifted action on $\tilde{M}$ is proper modulo the center; that is, for almost every $m \in M$, for any sequence $g_i \in \tilde{G}$ with $\mathrm{Ad}(g_i) \to \infty$ in $\mathrm{Ad}(\tilde{G})$, one has $g_i m \to \infty$ (we thank R. Zimmer for drawing our attention to this crucial fact). Suppose now that for such $m \in \tilde{M}$, there exists a sequence $g_i \in \tilde{G}$, $g_i \to \infty$, with $g_i m$ bounded. Discarding a null set, by the above we must then have $g_i = z_i k_i$, where $z_i \in Z(\tilde{G})$, $z_i \to \infty$, and the $k_i$'s are bounded. However then the points $k_i m$ are all contained in some compact subset, and since by (1) the $Z(\tilde{G})$-action is proper (moreover, a finite index subgroup acts effectively by deck transformations over $M$), this is impossible.

3. Denote $N = \mathrm{Ker}\, \varphi$. The assumption that $\varphi$ is an isogeny means that $N \subset p^{-1}(e)$. By 1 every element of $N$ acts properly on $\tilde{M}$, and we may consider the quotient $\tilde{M}/N$ or $N \backslash \tilde{M}$ (both of which make sense and coincide, because $N \subset \pi_1(M)$). We have a normal covering $\tilde{M}/N \to M$, and the group of deck transformations of $\tilde{M}/N$ over $M$ is naturally isomorphic to $\Lambda$. The $\tilde{G}$-action on $N \backslash \tilde{M}$ factors through an action of $G_1 = \tilde{G}/N$, which is (measurably) proper by 2.

4. *The index of* $N$ *in* $p^{-1}(e)$ *is finite. Consequently,* $G_1$ (*defined in* 3)
*has finite center.* Indeed, if $[p^{-1}(e): N] = \infty$, it follows that $\Lambda$ has infinite
center. Since we assume that it is embedded discretely in a rank one simple
Lie group with finite center, $\Lambda$ must be amenable, hence $\pi_1(M)$ is amenable as
well: $1 \prec \ell^2(\pi_1(M))$. Inducing the latter weak containment of representations
to $G_1$, and applying Lemma 6.2 (iii), (iv), we get $1 \subseteq L^2(M) \prec L^2(\tilde{M})$.
However, by step (3) above, $G_1$ acts measurably properly on $\tilde{M}$; hence, as in
the proof of Theorem 1.6, Mackey's orbit theorem implies that $L^2(\tilde{M})$ embeds
in a multiple of $L^2(G_1)$. But this and the previous conclusion imply together
that $1 \prec \infty \cdot L^2(G_1)$, contradicting the nonamenability of $G_1$.

*Conclusion of the argument.* By 3 we see that $\tilde{M}/N \to M$ is a covering $\Lambda$-
bundle with compact base, on which $G_1$ acts by bundle automorphisms, acting
properly on $\tilde{M}/N$. From 4 it follows that $G_1$ is commensurable with $\mathrm{SU}(n, 1)$,
and in particular $p(G_1) = p(G)(= 2n)$. At this point, exactly the same proof
as in the case of $\mathrm{SO}(n, 1)$ applies.                                         $\square$

III. *Proof of Theorem* 1.7. Reversing sides for convenience, suppose by
contradiction that $M = \Lambda \backslash G/\Gamma$ is compact. Dividing out the finite center, we
may assume that $G$ is (real) algebraic. Since the Zariski closure of $L$ commutes
with $\Lambda$ and satisfies the same hypothesis as in the theorem, we may also assume
that $L$ is algebraic. In particular, by nonamenability of $L$, it contains a copy
of $(P)\mathrm{SL}_2(\mathbb{R})$, denoted henceforth by $S$. As $L$ commutes with $\Lambda$, it acts on $M$
by left multiplication. Denote $\tilde{M} = G/\Gamma$. Then both $L$ and $\Lambda$ act from the left
on $\tilde{M}$, in a commuting way, so $\tilde{M} \to M$ is a covering $\Lambda$-bundle with compact
base, on which $L$ acts by bundle automorphisms (Def. 6.1).

Assume now that $\Lambda$ is embedded discretely in some $\mathrm{SO}(n, 1)$ or $\mathrm{SU}(n, 1)$.
Let $\pi$ be a unitary representation of that Lie group which is strongly $L^p$
for some $p < \infty$, and has nontrivial first cohomology (see the proof of Theorem
4.2). Then by Theorem 3.4 and Corollary 2.9, restricting $\pi$ to $\Lambda$ yields a
strongly $L^q$ representation for some $q < \infty$, which we continue to denote
by $\pi$, with a *proper* cocycle $b \in H^1(\Lambda, \pi)$ (Lemma 6.6). Consider the induced
$L$-representation $\sigma = \mathrm{Ind}_\Lambda^L \pi$. We claim that $\sigma|_S$ is strongly $L^r$ for some $r < \infty$.
Indeed, to invoke Theorem 6.3 for that purpose, we only need to explain why
the $S$-representation on $L^2(\tilde{M}) = L^2(G/\Gamma)$ is strongly $L^t$ for some $t < \infty$.
For this, we use the (necessary) assumption that $\Lambda$ is infinite, which implies
that $\Gamma$ is not a lattice in $G$ (indeed, by Poincaré recurrence, the $\Lambda$-action on
$G/\Gamma$ would not then be proper). Therefore $L^2(G/\Gamma)$ does not admit a nonzero
$G$-invariant vector. By a well-known result of Cowling [Co2], together with
the fact that $G$ must have real rank at least 2 (by the existence of a simple
noncompact subgroup $S$ which commutes with a noncompact subgroup $\Lambda$),
we deduce that $L^2(G/\Gamma)$ is strongly $L^t$ for some $t < \infty$. The latter holds as a $G$-representation, and hence also when restricted to any closed subgroup, particularly, to $S$ (cf. [Ho, 6.4]). We have thus shown that $\sigma|_S$ is strongly $L^r$ for some $r < \infty$, and since $S = (P)\mathrm{SL}_2(\mathbb{R})$ is not amenable, we conclude that $\sigma|_S$ does not contain weakly the trivial representation ([Sh3, Th. C]; see also step (iv) in the proof of Theorem 4.2 above). In particular, this conclusion holds for $\sigma$ as an $L$-representation.

Recall that the original $\Lambda$-representation $\pi$ came with an associated proper 1-cocycle; hence by Lemma 6.7 we have $H^1(L, \sigma) \neq 0$. (By the last part of Lemma 6.7, conditions (i) and (ii) there are indeed satisfied in our case; see below for the second one.) We now call on a result in [Sh2], concerning the so-called first *reduced* cohomology group $\overline{H^1}$ (see the discussion preceding Theorem 9; briefly, with the notation of Section 3.I above, $\overline{H^1}$ is the quotient of $B^1$ by the *closure* of $Z^1$ in the topology of uniform convergence on compact subsets). On the one hand, an easy and well-known fact is that for representations without almost invariant vectors, the two cohomology groups $\overline{H^1}$ and $H^1$ coincide (cf. [Sh2, 1.6]). On the other hand, it is shown in [Sh2, 3.6] that for any unitary representation $\sigma$ of a locally compact group $L$, satisfying $\overline{H^1}(L, \sigma) \neq 0$, the center $Z(L)$ admits a nonzero invariant vector. We deduce that the induced $L$-representation admits a nonzero vector, invariant under $Z(L)$.

Finally, we observe that the $Z(L)$-action in $\sigma$ must be mixing, a contradiction which will complete the proof of the theorem. Indeed, this follows from Lemma 6.2 (ii), whose two assumptions are verified using Howe-Moore's theorem [Zi1, 2.2.20]; that theorem applies to the original $\Lambda$-representation $\pi$, being a restriction of a representation of a simple Lie group without invariant vectors, and to the representation on $L^2(\tilde{M}) = L^2(G/\Gamma)$, using our previous observation that $\Gamma$ cannot be a lattice in $G$. This completes the proof of Theorem 1.7. □

## 8. Further remarks and related questions

8.1. The results of Section 4 suggest the following question: Suppose that $\pi$ is an irreducible unitary $\Gamma$-representation ($\Gamma < G = \mathrm{SO}(n,1), \mathrm{SU}(n,1)$, a lattice), such that $H^1(\Gamma, \pi) \neq 0$, and $\pi$ is strongly $L^{p(\Gamma)}$ (i.e., strongly $L^{\delta(\Gamma)}$). Is $\pi$ necessarily a restriction of a $G$-representation? In other words, is there only one such $\pi$ for $\Gamma < \mathrm{SO}(n,1)$ and two for $\Gamma < \mathrm{SU}(n,1)$? Notice that by [CS] the above $\Gamma$-representations are indeed irreducible (and in the case of $\mathrm{SU}(n,1)$, the two are non isomorphic). Of course, one has to exclude here $\mathrm{SO}(2,1) \cong \mathrm{PSL}_2(\mathbb{R})$. An affirmative answer to this question would be quite striking, and come close to a representation theoretic proof of Mostow's rigidity theorem. Indeed, if $\varphi : \Gamma_1 \to \Gamma_2$ is an isomorphism between two lattices of $G$, then the $\Gamma_1$-representation $\pi \circ \varphi$, when $\pi$ is the restriction to $\Gamma_2$ of an irreducible $G$-representation with $H^1 \neq 0$, is also a strongly $L^p$ cohomological representation, with 'minimal' $p$. Knowing that $\pi \circ \varphi$ extends to $G$ yields interesting information on $\varphi$, and the missing piece to deduce that $\varphi$ itself extends to $G$, is *exactly* the main result of [BiSt], proved for $\mathrm{PSL}_2(\mathbb{R})$. This beautiful theorem in [BiSt] is trivial for other rank one lattices, once one uses Mostow's rigidity; however the fact that it holds for $\mathrm{PSL}_2(\mathbb{R})$, in which rigidity is typically less seen, suggests that it might be proved independently, for other groups.

8.2. In view of the simplest, geometric constructions, of amalgam splittings $\Gamma = A *_C B$ for lattices $\Gamma < \mathrm{SO}(n,1)$, where $C$ is the fundamental group of a totally geodesic, co-dimension one separating hyperplane, one is naturally led to ask whether this must be the case whenever equality in Theorem 1.6 is satisfied, namely, when $\delta(C) = \delta(\Gamma) - 1$. Such a result would be analogous to Bourdon-Yue's superrigidity theorem mentioned after Theorem 1.1. The analysis in §5.II provides evidence to an affirmative answer, and especially to our (apparently weaker) conjecture that in case of equality, $C$ must be a lattice in $\mathrm{SO}(n-1,1)$. Recall that it was shown in §5.II that in case of equality, the bottom of the Laplacian spectrum in $C\backslash G/K$ is attained, which yields significant geometric information (see [Su2]). The question of existence of nonelementary actions on trees, of lattices in $\mathrm{SU}(n,1)$, seems open in general.

8.3. Separate discussion of the different possibilities which can occur in inequality (1) of Theorem 1.4 leads to a more illuminating result, putting forward an intriguing question. Again, we shall not discuss the group $\mathrm{SL}_2(\mathbb{R})$, and let us leave aside, for the moment, also the group $\mathrm{SL}_2(\mathbb{C})$. Retain the notation in Theorem 1.4. By dividing out the finite center, we assume hereafter that $G$ is center free. Consider now the case where $\varphi$ is *not* injective. Then $N = \mathrm{Ker}\,\varphi$ is a normal, *nonamenable* subgroup of $\Gamma$, which seems to imply in general that $\lambda_0$, the bottom of the spectrum of $\Delta$ for $N\backslash G/K$, is *strictly smaller* than its value $(\delta(G))^2/4$ for the symmetric space. This is the Riemannian analogue of Kesten's well-known theorem [Ke], asserting that for any generating, symmetric probability measure on a finitely generated group $\Gamma$, the norm of the corresponding convolution operator in $\ell^2(\Gamma)$ is strictly smaller than its norm in $\ell^2(\Gamma/N)$, for any *normal* nonamenable subgroup $N$ (the normality assumption is crucial). We have failed to find reference for the result in the Riemannian context, however, at our request it was recently proved in the co-compact case by Vadim Kaimanovich [Ka], to whom we extend our gratitude. Of course, if the bottom of the spectrum of $N\backslash G/K$ is strictly smaller than $\delta(G)^2/4$, then $\delta(N) > \delta(G)/2$ (see §2.V and the references therein).

Once this claim is established, one may break inequality (1) in Theorem 1.4 into a more refined version. Indeed, with the notation of Theorem 1.4, and denoting $\varphi(\Gamma) = \operatorname{Im} \varphi$, the following possibilities occur:

(i) Either $\varphi$ is *injective*, implying $\delta(\Gamma) \leq \delta(\operatorname{Im} \varphi)$   (Theorem 1.1); or

(ii) $\delta(\operatorname{Im} \varphi) < 2$, implying $\delta(\Gamma) < \delta(\operatorname{Ker} \varphi) + 1$   (Theorem 5.6); or

(iii) $\varphi$ is *noninjective* and $\delta(\operatorname{Im} \varphi) \geq 2$, implying $\delta(\Gamma) \leq \delta(\operatorname{Ker} \varphi) + \frac{1}{2}\delta(\operatorname{Im} \varphi)$ (Theorem 1.4).

The superrigidity theorem of Bourdon-Yue (and its strengthening in [BCG]), mentioned in the introduction, yields a complete understanding of case (i) (for uniform lattices), and it seems that not much more can be said in case (ii), which already provides a strict inequality. Thus, if one could prove that in case (iii), strict inequality always holds as well (which we believe to be the case), this would show that equality in Theorem 1.4 occurs if and only if the homomorphism $\varphi$ extends to the ambient group. Returning now to the case of $\operatorname{SL}_2(\mathbb{C})$, we see that, curiously enough, this superrigidity-type result indeed holds, at least for uniform lattices, simply because $\delta(\Gamma) = 2$ (relying only on Kaimanovich's result mentioned above). Illustration of the wealth of situations which can occur for noninjective homomorphisms can be found in an example of Mostow [Mo, §22], where an arithmetic lattice $\Gamma < \operatorname{SU}(2,1)$ is shown to map surjectively, with infinite kernel, onto a different arithmetic lattice of that group.

8.4. Our proof of Theorem 2.1 made essential use of the full continuous strip of complementary series representations. While for the other rank one groups, results of Cowling regarding existence of appropriate uniformly bounded representations may serve as a substitute, for higher rank groups some new idea seems to be required. By Cowling's general "Kunze-Stein phenomenon" [Co1], for any simple Lie group and *even p*, once a representation is strongly $L^p$, then *all* its matrix coefficients are in $L^{p+\epsilon}$ for all $\epsilon > 0$. However, even for $\operatorname{SO}(n,1)$ and $\operatorname{SU}(n,1)$, it is not known whether for any irreducible representation, all matrix coefficients lie in the same $L^p$. We remark that by applying a general result of [KS] in addition to Theorem 2.1, it follows that for any irreducible representation $\pi$ of $G = \operatorname{SO}(n,1), \operatorname{SU}(n,1)$, if there is one matrix coefficient in $L^{p+\epsilon}(G)$, then all matrix coefficients of the *smooth* vectors satisfy a decay estimate similar to (6), essentially replacing the Hilbert norm by a Sobolev type norm. Another related question, considered by Howe [Ho] (see the remarks preceding the proof of Corollary 7.3 there), is whether for every $p$, not only an even integer, the set of strongly $L^p$ representations is closed in the Fell topology. Theorem 2.1 (4) yields a positive answer for $\operatorname{SO}(n,1)$ and $\operatorname{SU}(n,1)$.

8.5. It is natural to look for more classes of groups $\Gamma$, for which $0 < p(\Gamma)$ $< \infty$, especially with a precise estimate. By Theorem 4.5, this would yield a lower bound $\delta(\varphi(\Gamma)) \geq p(\Gamma)$ for any discrete faithful representation $\varphi$ into $SO(n,1)$ or $SU(n,1)$. Discrete subgroups of the latter groups are themselves a rich source of examples, and especially intriguing would be to find there non-lattices which satisfy equality $\delta(\Gamma) = p(\Gamma) > 0$. Can one find a group $\Gamma$ for which $p(\Gamma)$ is not an integer? What can be said about the numerous known hyperbolic groups acting discretely co-compactly on CAT(-1) spaces (see [HP] and the references cited in §5.I) ?

Finally, as the finitely generated nonamenable groups with $p(\Gamma) < 2$ (equivalently, $p(\Gamma) = 0$), are exactly the ones having nonvanishing first $\ell^2$-Betti number, one can consider two properties which may be viewed as "almost having" positive first $\ell^2$-Betti number: satisfying $p(\Gamma) = 2$, and admitting a representation with $H^1 \neq 0$, which *weakly* embeds in the regular representation. By [CHH], the first implies the second, but they are not equivalent. Fundamental groups of closed (or finite volume) hyperbolic 3-manifolds satisfy the former, and groups which split as an amalgamated product over an amenable subgroup, the latter.

# Appendix: Proof of Theorem 3.7

Although in the present paper we are interested only in $SO(n,1)$ and $SU(n,1)$, our treatment of the problem of inducing the first cohomology enables one to induce isometric actions on general CAT(0) spaces. In that framework, the other rank one groups are of relevance as well. The exceptional group $F_{4(-20)}$ can be dealt with by similar methods, and is sufficiently specialized that we will not describe it here. Thus, in terms of the corresponding symmetric space, we will work here with either the complex hyperbolic space $\mathbb{CH}^n$ ($G = SU(n,1)$), or the quaternionic hyperbolic space $\mathbb{QH}^n$ ($G = Sp(n,1)$), where $n \geq 2$.

1. *Statement of the problem and some notation.* Let $X = G/K$ be either $\mathbb{CH}^n$ or $\mathbb{QH}^n$. Let $\Gamma < G$ be a nonuniform lattice, which we may assume to be torsion free. (Hereafter, for convenience, we let $\Gamma$ act from the left on $G$ and $G/K$.) Let $F \subset G/K$ be a fundamental domain for the action of $\Gamma$. As mentioned in Section 3.IV, such a fundamental domain consists of finitely many cusps, and as in there, it will suffice to study the behavior of one such cusp $\Omega$, which we fix once and for all in what follows. Let $\infty \in \partial X$ be the geodesic ray determined by $\Omega$. Let $\{H_t \mid t \in \mathbb{R}\}$ be the usual foliation of horospheres which are based at $\infty$. We have chosen our parameter so that the distance from $H_s$ to $H_t$ is $|s - t|$.

Let $\Gamma_\infty < \Gamma$ be the intersection of the cusp group $P$ (the subgroup of $G$ fixing $\infty$), with $\Gamma$. The subgroup $\Gamma_\infty$ preserves each $H_t$, and acts on it co-compactly (cf. [GR]). Let $F_0 \subset H_0$ be a (compact) fundamental domain for $\Gamma_\infty$. We define $F_t \subset H_t$ to be the set of points which lie on geodesics determined by $\infty$ and by points of $F_0$. By symmetry, $F_t$ is a fundamental domain for the action of $\Gamma_\infty$ on $H_t$. Furthermore, $\Omega$ is the union of a compact part, irrelevant to our discussion (as our function is locally bounded), and, say, $\cup_{s \geq 0} F_s$ (denoted $F^0$ below).

As is well-known by now, $\Gamma$ is finitely generated (again, see [GR]). We fix a finite generating set once and for all, which may be assumed to contain the generators of $\Gamma_\infty$, and equip $\Gamma$ with the corresponding word length. Given any point $x \in X$, let $w(x)$ be the word length of the element $\gamma \in \Gamma$ for which $\gamma x \in F$. Finally, define the function $f(p)$ on $F$, and the regions $F^s \subset F$ by:

$$f(p) = \sup\{w(x) \mid d(x, p) \leq 1\} \qquad F^s = \cup_{t \geq s} F_t.$$

Our goal is to show that $f$ is square integrable on $\Omega$, or equivalently, on $F^0$ (see Theorem 3.8 and the discussion preceding it). $\qquad \square$

2. *Horospherical projection.* Given $t > 0$, there is a canonical projection $\pi_t : H_t \to H_0$. This map is defined so that the geodesic determined by $x$ and $\pi_t(x)$ intersects $\partial X$ at $\infty$. Notice that $\pi_t$ commutes with the action of $\Gamma_\infty$. Since $F^0$ is foliated by horospheres, we may combine the maps $\pi_t$ into a single projection $\pi : F^0 \to H_0$.

We will equip $H_0$ with the Riemannian metric induced from the Riemannian metric on $X$. Let $B(m)$ be the portion of $F^0$ which is bounded by $H_m$ and $H_{m+1}$. In other words, let $B(m) = F^m - F^{m+1}$. Let $B_1(m)$ be the unit tubular neighborhood of $B(m)$ in $X$. Let $g(m)$ be the diameter of $\pi(B_1(m))$, and $v(m)$ be the volume of $B(m)$.

2.1 LEMMA. *The function $f$ is square integrable on $F^0$ if the sequence $\{g(m)^2 v(m) \mid m \in \mathbb{N}\}$ is summable.*

*Proof.* Given $p \in B(m)$ we define $\hat{g}(p) = g(m)$. It suffices to show that $f(p) < C\hat{g}(p)$ for some universal constant $C$ and for all $p \in F^0$. If $x \in X$ is such that $d(x, p) \leq 1$, then $x \in B_1(m)$ and $\pi(x) \in \pi(B_1(m))$. Thus, $\pi(x)$ has distance at most $g(m)$ from any point in $F_0$.

Since $F_0$ is equipped with a path metric, and $H_0/\Gamma_\infty$ is compact, $H_0$ is quasi-isometric to the set $\Gamma_\infty$, equipped with its word metric. This is the well-known criterion of Milnor-Svarc (see [Sc] for details). Thus, there is a universal constant $C$ such that a word of length at most $Cg(m)$ maps $\pi(x)$ to a point in $F_0$. Since $\Gamma_\infty$ commutes with $\pi$, this same word maps $x$ into $F^0$. $\qquad \square$

It remains to estimate $g(m)$ and $v(m)$.

3. *Estimates.* As a preliminary step, we describe the Carnot metric on the horospheres $H_t$ (for details see [Gr3] or [Sc]). $H_t$ has a canonical nonintegrable distribution. This distribution has codimension one in case $X = \mathbb{C}\mathbb{H}^n$, and has codimension three when $X = \mathbb{Q}\mathbb{H}^n$. The Riemannian metric on $X$ induces a Riemannian metric on the planes of the distribution. One defines lengths of curves, integral to the distribution, via integration. It is well-known that every two points of $H_t$ can be connected by an integral path, so that the infimal length of integral paths from $x$ to $y$, called the *Carnot distance*, is finite.

3.1 LEMMA. *There is a constant $C_0$ such that $g(m) < C_0 \exp(m)$.*

*Proof.* The map $\pi_t$ preserves the relevant distributions and $T(x) = \exp(t)x$, for every vector $x$ tangent to such a plane in the distribution on $H_t$. If follows immediately that

$$d_C(\pi_t(x), \pi_t(y)) = \exp(t)d_C(x, y),$$

for all $x, y \in H_t$. Here $d_C$ is the Carnot distance.

Every two points $x, y \in B_1(m)$ can be connected by a path of the form $\alpha \cup \beta$, where $\pi(\alpha)$ is a point, and $\beta \in H_t$ is integral, with $t \in [m, m+1]$. Since every $B_1(m)$ can be isometrically embedded in $B_1(k)$ for $m > k$, we can find a constant $C_0$, independent of $m$, such that $\beta$ has Carnot length at most $C_0$. By construction, $\pi(\alpha \cup \beta)$ has Carnot length at most $C_0 \exp(m)$. Hence, the path metric distance, in $H_0$, from the endpoints of this projected path, is at most $C_0 \exp(m)$.                                                                        □

3.2 LEMMA. *For $X = \mathbb{C}\mathbb{H}^n$, there is a constant $\tilde{C}_1$ such that*

$$v(m) < \tilde{C}_1 \exp(-2nm).$$

*Here $\tilde{C}_1$ depends on the group $\Gamma$.*

*Proof.* The linear differential $d\pi_t$ expands distances by $\exp(t)$ in all directions tangent to the nonintegral distribution on $H_t$. It is a well-known fact that $d\pi_t$ expands distances by $\exp(2t)$ along one of the tangent vectors to $H_t$, which does not lie in the distribution. More precisely, this "extra-expanding vector" is generated by the Lie brackets of vectors tangent to the distribution. Thus, $d\pi_t$ multiplies the volume form on $H_t$ by $\exp((2n-2)t+2t) = \exp(2nt)$. Put the other way, $d\pi_t^{-1}$ multiplies volumes on $H_0$ by $\exp(-2nt)$. The lemma now follows readily from integration.                                                                        □

3.3 LEMMA. *For $X = \mathbb{Q}\mathbb{H}^n$, there is a constant $\tilde{C}_2$ such that*

$$v(m) < \tilde{C}_2 \exp(-4nm - 2m).$$

*Here $\tilde{C}_2$ depends on the group $\Gamma$.*

*Proof.* The proof is the same as above, except that there is now a set of three linearly independent vectors tangent to $H_t$, which $d\pi_t$ expands by $\exp(2t)$. $\qquad\square$

*Conclusion of the argument.* Note that $n \geq 2$. Therefore, by 3.2 and 3.3 we have, for the two families, $v(m) < C \exp(-4m)$. Together with Lemma 3.1, we get

$$g(m)^2 v(m) < \tilde{C} \exp(2m) \exp(-4m) < \tilde{C} \exp(-2m),$$

a sequence which is obviously summable. By Lemma 2.1 we are done. $\qquad\square$

Princeton University, Princeton, NJ
*Current Address:* Yale University, New Haven, CT, and Hebrew University, Jerusalem, Israel
*E-mail address:* yshalom@math.yale.edu, yehuda@math.huji.ac.il